%% file: egpaper_for_review.tex
\documentclass[10pt,twocolumn,letterpaper]{article}

\usepackage{wacv}
\usepackage{times}
\usepackage{epsfig}
\usepackage{graphicx}
\usepackage{amsmath}
\usepackage{amssymb}

\makeatletter
\@namedef{ver@everyshi.sty}{}
\makeatother
\usepackage{amsfonts}
\usepackage{subfigure}
\usepackage{tikz}
\usepackage{booktabs}
\usepackage{multirow}
\usepackage{pgfplots}
\usepackage{overpic}
\usepackage{tabularx}
\newlength\figureheight 
\newlength\figurewidth 

\def\eqref#1{(\ref{eqn:#1})}
\def\eqlabel#1{\label{eqn:#1}}

\def\be{\begin{equation}}
\def\ee{\end{equation}}

\def\code#1{{\tt #1}}
\def\vcat#1{\left[\begin{matrix}#1\end{matrix}\right]}
\def\hcat#1#2{\left[#1 ~\middle|~ #2 \right]}

\usepackage{microtype}

\wacvfinalcopy 

\def\wacvPaperID{413} 

\ifwacvfinal\fi
\begin{document}

\title{QRkit: Sparse, Composable QR Decompositions for\\Efficient and Stable Solutions to Problems in Computer Vision}

\author{Jan Svoboda \\
USI Lugano, Switzerland\\
{\tt\small jan.svoboda@usi.ch}
\and
Thomas Cashman \hspace{2cm} Andrew Fitzgibbon \\
Microsoft HoloLens, Cambridge, UK\\
{\tt\small {tcashman,awf}@microsoft.com}
}

\maketitle
\ifwacvfinal\thispagestyle{empty}\fi

\begin{abstract}
Embedded computer vision applications increasingly require the speed and power benefits of single-precision (32 bit) floating point.  
However, applications which make use of Levenberg-like optimization can lose significant accuracy when reducing to single precision, sometimes unrecoverably so.
This accuracy can be regained using solvers based on QR rather than Cholesky decomposition, but the absence of sparse QR solvers for common sparsity patterns found in computer vision means that many applications cannot benefit.
We introduce an open-source suite of solvers for Eigen, which efficiently compute the QR decomposition for matrices with some common sparsity patterns (block diagonal, horizontal and vertical concatenation, and banded).
For problems with very particular sparsity structures, these elements can be composed together in `kit' form, hence the name {\em QRkit}. We apply our methods to several computer vision problems, showing competitive performance and suitability especially in single precision arithmetic.
\end{abstract}


\section{Introduction}

\def\sz#1#2{#1\times#2}
\noindent
Computer vision applications are increasingly required to run on low-power architectures, where single precision floating point is significantly more efficient than double~\cite{Baboulin2009,Langou2006,RoldaoLopes2009}.  However, where such applications depend on Levenberg-like solvers (e.g. bundle adjustment~\cite{Agarwal2010, Triggs2000, Zach2014}, SLAM~\cite{Cadena2016}, 3D reconstruction~\cite{Tomasi1992}, surface fitting~\cite{Cashman2013,Tagliassacchi2015,Taylor2016}), single precision operation can negatively impact accuracy, and may therefore require more solver iterations, or may simply never be accurate enough.

\begin{figure}[t]
\centering
\subfigure[]
{
\includegraphics[width=17mm,height=32mm]{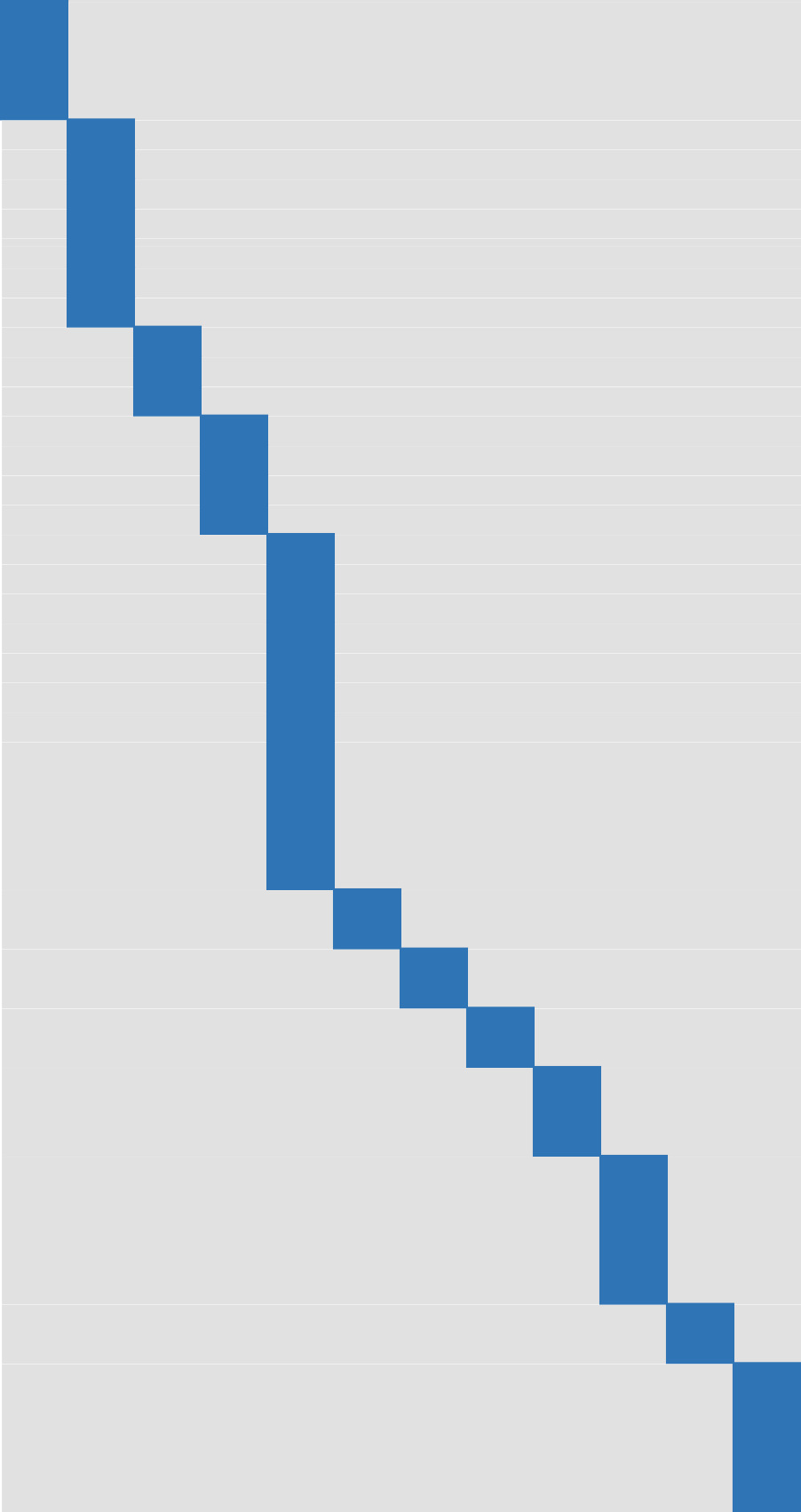}
\label{fig:patBlockDiag}
}
\subfigure[]
{
\includegraphics[width=19mm,height=32mm]{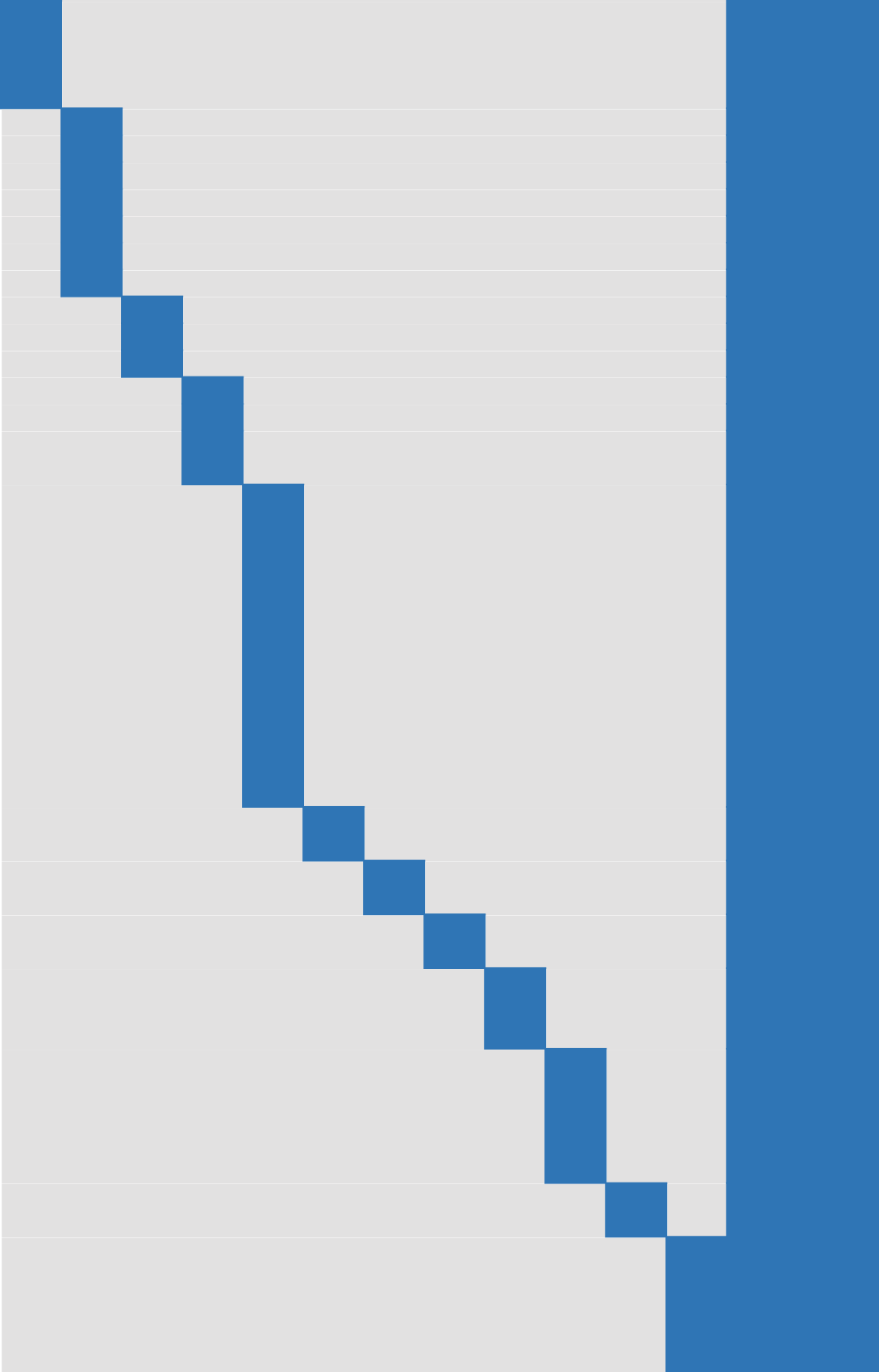}
\label{fig:patHorzcat}
}
\subfigure[]
{
\includegraphics[width=17mm,height=32mm]{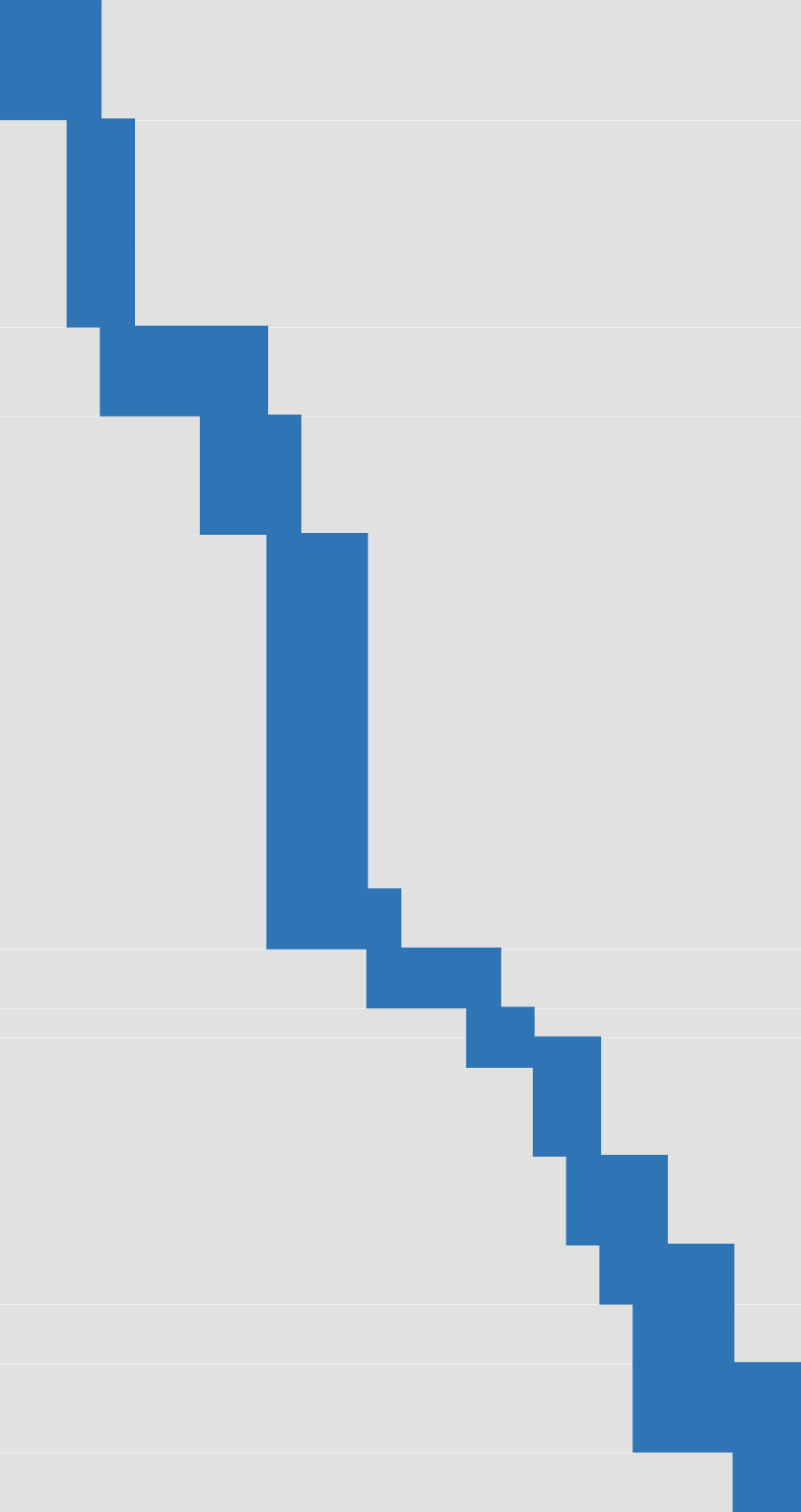}
\label{fig:patBlockBand}
}
\subfigure[]
{
\includegraphics[width=17mm,height=32mm]{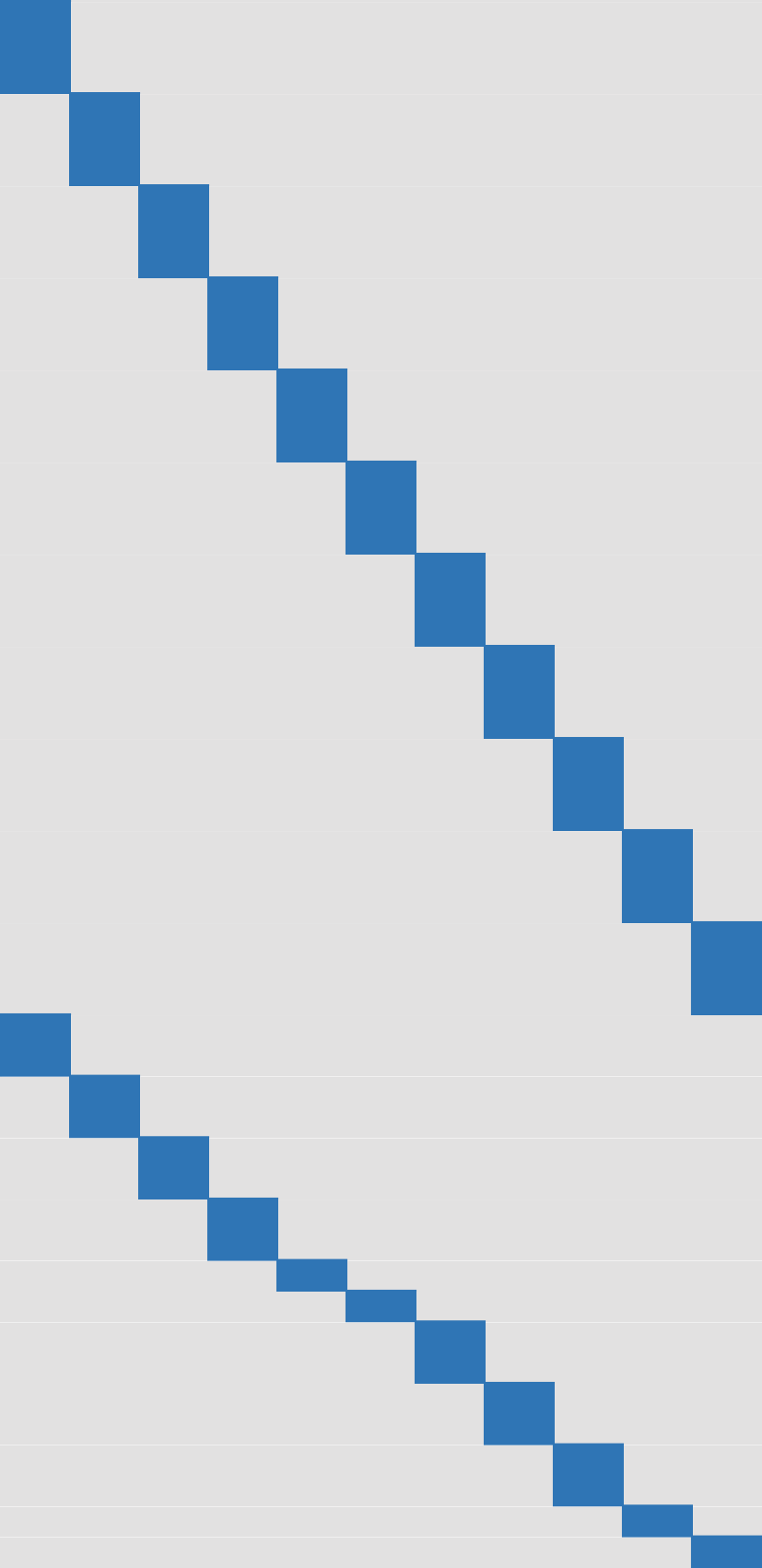}
\label{fig:patVertcat}
}
\caption{Sparsity patterns. (a) Block diagonal, (b) Horizontal concatenation of {\tt [Easy, Dense]}, including Angular, (c) Banded, (d) Row-permuted banded (includes vertical concatenation of some easier cases). }
\label{fig:sparsityPatterns}
\end{figure}

In this paper, we identify an important contributor to such inaccuracy: the use of Cholesky decomposition to solve equations of the form
\be 
(J^\top J + \lambda D^2) p = J^\top b
\label{eq:LMstep}
\ee
for unknown vector $p$ given: vector~$b$, diagonal matrix~$D$, and matrix~$J$.
Typically $J$ is rectangular and sparse with particular sparsity structures (some are illustrated in Figure~\ref{fig:sparsityPatterns}) that can be exploited in Cholesky decomposition.  However, $J^\top J$ has a condition number which is the square of~$J$, which can adversely affect numerical precision and robustness. 

The squaring can be avoided by solving the equivalent least squares system
\be 
\begin{pmatrix} J\\\lambda^{\frac12} D\end{pmatrix} p = 
\begin{pmatrix} b \\ 0 \end{pmatrix}
\label{eq:LevMarqQR}
\ee
using QR decomposition, as is in the classic Levenberg--Marquardt implementation of Mor\'e~\cite{More1983}.
In the dense case, this is typically twice as many floating-point operations per iteration than the Cholesky solution, but the reduction in iteration count for a given accuracy can be significant.

In the sparse case, the story is less rosy.
There are general-purpose sparse QR implementations, e.g.\ the SPQR routines in SuiteSparse~\cite{Davis2011}, which make the method competitive with Cholesky for many problems.  However, there are no libraries which offer QR solvers which can exploit the particular sparsity structures in typical computer vision problems, so the QR method is unusably slow, even though it may be ultimately more accurate.

The contribution of this paper is to collate from the numerical analysis literature a small collection of special-purpose solvers, with the property that they can be composed in `kit' form to easily build fast solvers for a wide range of computer vision applications.  The wider contribution is to allow the use of efficient single-precision optimization routines without catastrophic loss of accuracy.  

We structure the paper by first reviewing the properties and use of the QR decomposition.
Second, we explain our strategy for building efficient solvers, and present a collection of such solvers.  
Third, we show results on a range of computer vision benchmarks, for QR-based and Cholesky-based algorithms in both single and double precision.

\def\e#1{\mathbf{e}_{#1}}
\def\qr#1{\code{qr}(#1)}
\def\blk#1#2{\code{blk_{#1}}{#2}}
\def\kron#1#2{#1 \otimes #2}
\def\blkdiag#1{\mathtt{blkdiag}(#1)}
\def\zeros#1{\mathtt{zeros}(#1)}

\paragraph{Notation}
Let $\e i$ denote the $i$th column of the $\sz n n$ identity matrix. The Kronecker product $\kron AB$ puts a scaled copy of $B$ at every entry in $A$, for example if $B$ is $\sz pq$,
then $\kron{\e i}{B}$ is a matrix of size $\sz{np}q$ with structure $[0_{\sz{(i-1)p}q} ; B ; 0_{\sz{(n-i+1)p}q} ]$.
To flag where a QR decomposition is computed numerically, we write $Q R = \qr A$.

\def\Qobj{\code{Q}}
\def\Qfull{\mathbf{Q}}
\def\Rfull{\mathbf{R}}
\def\Qperp{Q^{\perp}}
\def\Qt{Q^\top}
\section{Background: QR decomposition}
\label{sec:BackgroundQR}
Given a matrix $A$, the QR decomposition finds matrices $Q$ and $R$ such that $Q$ is orthonormal ($Q^\top Q = I$), $R$ is upper triangular, and $A = QR$.
We generally deal with matrices $A$ in `portrait' orientation, so $m > n$, and with the `economy-size' decomposition where $Q$ has size $\sz m n$ and $R$ is $\sz n n$.  However, we will need access to the orthogonal complement of $Q$, written $Q^\perp$, size $\sz{(n-m)}{n}$, and to the `full size' matrices
\be 
\Qfull = \hcat{Q}{\Qperp},  \qquad{\Rfull} = \vcat{R\\0}
\eqlabel{Qfull}
\ee 
which will often be too large to store explicitly.  Therefore, as is common in QR decomposition routines, we may store not $Q$ itself, but some C++ object $\Qobj$ which behaves like $Q$; that is, it implements the operations of matrix multiplication.
Then, to look at $Q$ itself (for example, when debugging small problems), we multiply $\Qobj$ by a representation of the $\sz n n$ identity matrix.

The rich history of QR decomposition dates back to the early 1950s, when it was introduced independently by 
Francis~\cite{Francis1961Part1,Francis1961Part2} and 
Kublanovskaya~\cite{Kublanovskaya1962}, who proposed the use of QR decomposition as a solution to the eigenvalue problem. 
QR decomposition can be carried out 
using Gram-Schmidt orthogonalization~\cite{Bjorck1994}, Givens rotations~\cite{Givens1958} or Householder reflections~\cite{Householder1958}; 
all the approaches were summarized later by Gander~\cite{Gander1980}. 
The method became popular at the time, exploited by many showing its use for singular value decomposition~\cite{Golub1970} and solution to least squares problems~\cite{Golub1970}. 

Many important methods are based on Householder matrices of the form $H = I - 2vv^T$, where $v$ is a Householder vector~\cite{Householder1958} having unit 2-norm. Householder matrices are orthogonal and can be used to zero-out selected columns of a matrix. Householder QR decomposition of $\sz{n}{m}$ matrix $A$ can be therefore expressed as multiplication by sequence of Householder matrices $H_{m-1}, \cdots, H_0$, i.e.
\be
\Rfull = H_{m-1}\cdots H_0 A = \Qfull^T A,
\ee
where each $H_k, k \in \{0, \cdots, m-1\}$ is of size $\sz{m}{m}$.

Bischof \etal~\cite{Bischof1987} and Scheiber \etal~\cite{Schreiber1988} introduced the blocked versions of Householder transformation, in order to reduce the burden of inefficient matrix-vector operations on the supercomputing architectures of the time. Given matrices $B,C \in \mathbb{R}^{\sz{n}{r}}$ and $B_1,C_1 \in \mathbb{R}^{\sz{r}{r}}$ we want to find an orthogonal $\Qfull \in \mathbb{R}^{\sz{n}{n}}$ such that
\be
B = \vcat{B_1 \\ B_2}, \qquad\Qfull^T B = C = \vcat{C_1 \\ 0}.
\ee
Using the `WY' representation by Bischof \etal~\cite{Bischof1987}, the solution to the problem is represented as
\be
\Qfull = I + WY^T, \qquad W,Y \in \mathbb{R}^{\sz{n}{r}},
\label{eq:blockedHouseholderWY}
\ee
where $Y$ is lower trapezoidal, i.e., $y_{ij} = 0 \mathrm{~if~} i < j$. The submatrix $C_1$ is upper triangular and $\Qfull$ is a rank-1 correction to the identity, and so it can be regarded as a generalization of the Householder matrix~\cite{Schreiber1989}.

Schreiber \etal~\cite{Schreiber1989} show how to modify the WY representation so that only $nr$ storage is required. The matrix $\Qfull$ from (\ref{eq:blockedHouseholderWY}) can be expressed as
\be
\Qfull = I + YTY^T, \qquad Y \in \mathbb{R}^{\sz{n}{r}}, T \in \mathbb{R}^{\sz{r}{r}},
\ee
where $Y$ is lower trapezoidal and $T$ upper triangular. This is usually referred to as a compressed WY representation. For details on computation of blocked Householder representations, we refer the reader to \cite{Schreiber1989}.

Later interest in parallel computing motivated researchers to explore parallel QR decomposition algorithms~\cite{Cosnard1986, Oleary1990}, in a thread that continues as a subject of active research~\cite{Buttari2008,Demmel2008,Gunter2005}.

Alongside other methods, QR decomposition has been implemented as a part of the LAPACK package~\cite{Anderson1990} in 1990 and much more recently Davis~\cite{Davis2011} developed a multifrontal rank-revealing version (SPQR) in the SuiteSparse library.

A very appealing application of QR factorization arises in the field of numerical optimization, in particular solving non-linear least squares problems. It is of particular interest in use together with the Levenberg~\cite{Levenberg1944} Marquardt~\cite{Marquardt1963} algorithm. 
This is a popular variant of the Gauss--Newton method for finding the minimum of a function $F(x)$ represented as a sum of squares of generally nonlinear functions
\be
F(x) = \| f(x) \|^2 = \frac{1}{2}\sum_{i=1}^{m} [f_i(x)^2],
\ee
from which comes the above-mentioned instance of (\ref{eq:LMstep}), with the matrix $J$ being the Jacobian of $f : \mathbb R^n \mapsto \mathbb R^m$, 
and the vector $b = -f(x)$.

There are various strategies for updating the damping parameter~$\lambda$. Mor{\'e}'s implementation employs a trust-region method proposed by Hebden~\cite{Hebden1973}. 



To reduce the number of matrix factorizations, Lourakis \etal\cite{Lourakis2009} presented an alternative approach to updating the damping factor, better suited for computer vision latent variable problems such as bundle adjustment. 
Instead of seeking a nearly exact solution for $\lambda$ using Newton's algorithm in a trust-region framework (as proposed by Conn \etal~\cite{Conn2000} and Mor\'e~\cite{More1983}), it directly controls the damping factor $\lambda$ with a line-search algorithm~\cite{Lourakis2009}. 




\section{Building sparse QR solvers}
The general strategy to build efficient solvers~\cite{Bjorck2014,Demmel2008} is to express the matrix $A$ as some combination of smaller matrices $A_{1..K}$, for whose shape it is easy (i.e.\ efficient) to store and compute the QR decomposition.  Then manipulations of these easy QRs leads to the decomposition of $A$.  The composition of these solvers in code is a compile-time declaration in terms of the data types of component solvers.  For example, a matrix with four blocks (see Figure~\ref{fig:exampleA}) might have the top-left block defined as block diagonal, with $p$ component blocks of size $\sz23$, the top-right as general dense of size $\sz p q$, and the bottom-left block expressed as `block-banded'.  Consulting Figure~\ref{fig:exampleA}, the best way to encode this is as a vertical concatenation of $A_1$ and $A_2$ followed by horizontal concatenation with $A_3$ and a $0$ block, expressed as
\be
\hcat{\vcat{A_1 \\ A_2}}{\vcat{A_3 \\ 0}}.
\ee

\section{QRkit}
QRkit supports efficient factorization of the sparsity patterns and their compositions that will be described in the following sections.
\label{sec:qrkit}
\subsection{Block diagonal}
\label{sec:BlockDiag}
\begin{figure}[t]
\centering
\subfigure[$A$]
{
\includegraphics[width=26mm,height=32mm]{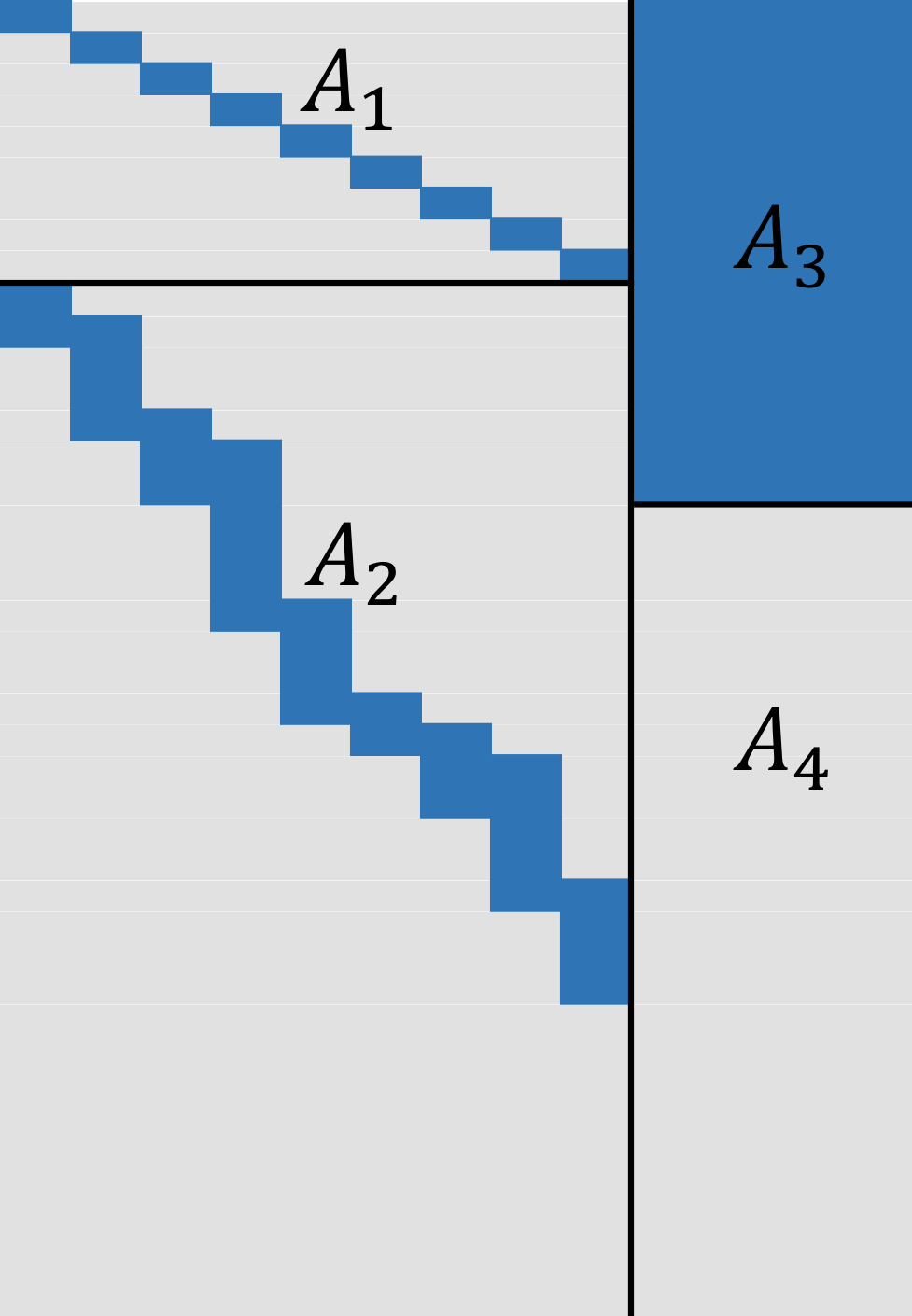}
\label{fig:exampleA}
}
\subfigure[$\Qfull$]
{
\includegraphics[width=32mm,height=32mm]{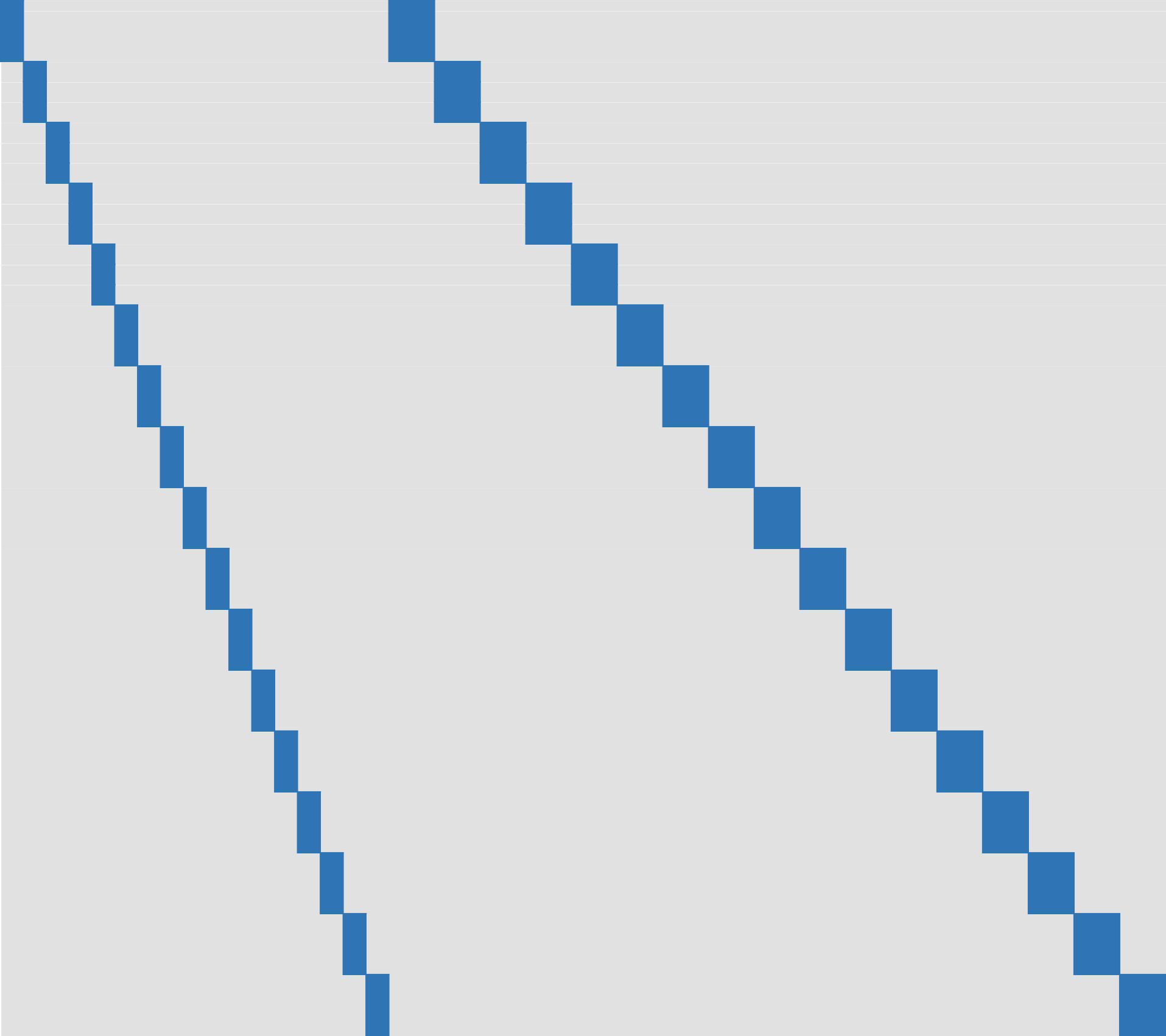}
\label{fig:blockDiagonalQ}
}
\subfigure[$R$]
{
\includegraphics[width=17mm,height=32mm]{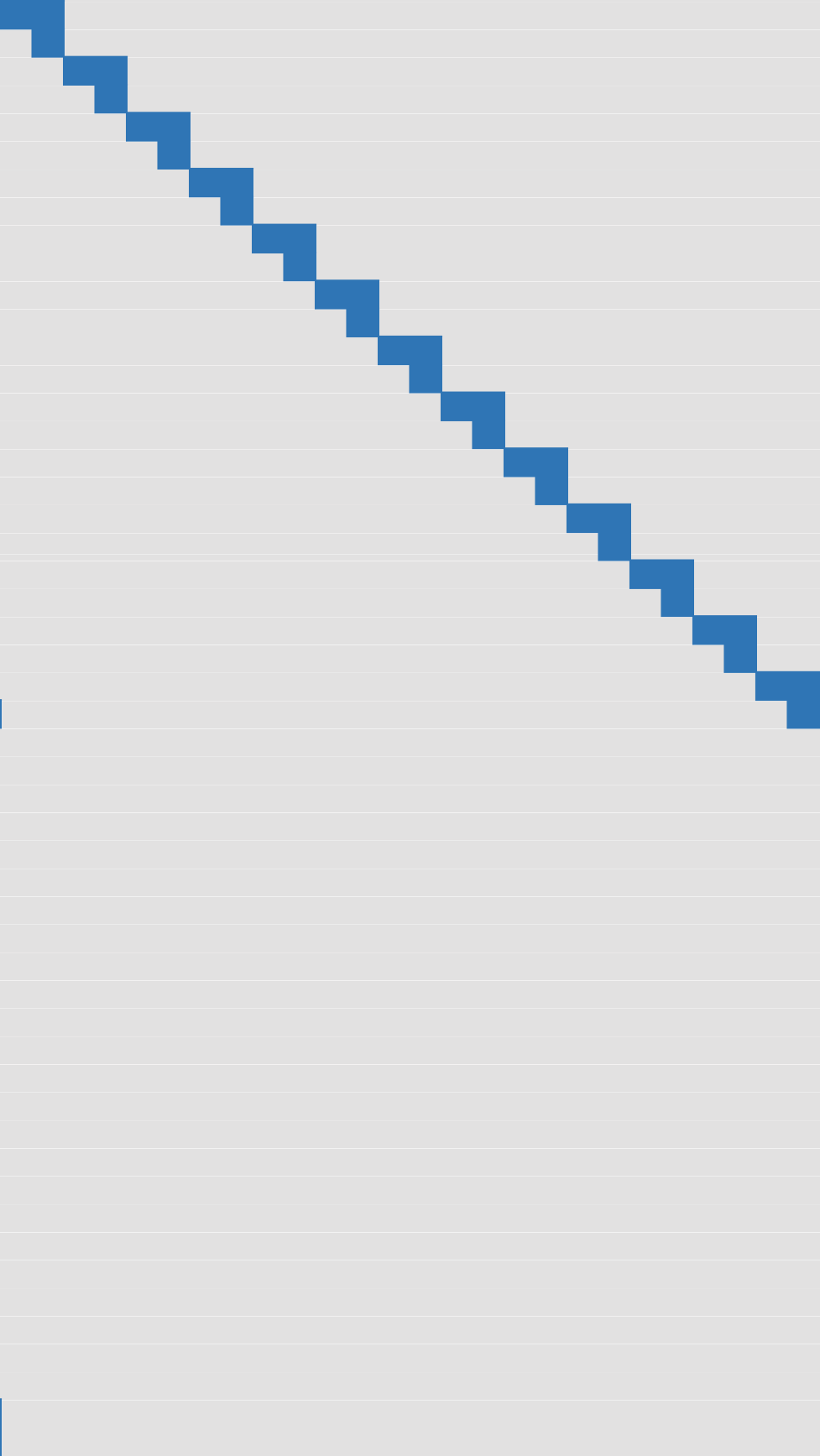}
\label{fig:blockDiagonalR}
}
\caption{ (a) An example of favorable sparsity patterns in different blocks of matrix $A$. (b-c) Typical sparsity of QR decomposition of a block diagonal matrix (e.g. block $A_1$ in $A$). In $\Qfull$, the left diagonal forms $Q$ and the right diagonal $\Qperp$, and $\Qfull = \hcat{Q}{\Qperp}$.}
\label{fig:blockDiagonalQR}
\end{figure}

If the matrix $A$ is block diagonal, i.e. $A = \blkdiag{A_1,...,A_K}$, then to find matrices $Q,R$ as above, we define $Q_k R_k = \qr{A_k}$ for all $k$. Observe that $\kron{\e j}{Q_j}$ is orthogonal to $\kron{\e k}{Q_k}$ for all $k \ne j$, because the nonzeros don't overlap, so we can simply write $Q = \blkdiag{Q_1, ..., Q_K}$, and $R = \blkdiag{R_1, ..., R_K}$.
The class \code{BlockDiagonalQR} is templated over the solver type \code{BlockSolver} of the individual blocks. 
For example
\begin{small}
\begin{verbatim}
 BlockDiagonalQR<DenseQR> slvr;
\end{verbatim}
\end{small}
For the simple block diagonal case, matrix $Q$ is very sparse and can therefore be formed explicitly as well as using a vector of \code{BlockSolver}s. The upper triangular factor $R$ exhibits strong sparsity as well and all of its elements are close to the matrix main diagonal. An example of QR factorization of a block diagonal matrix is shown Figure \ref{fig:blockDiagonalQ} and \ref{fig:blockDiagonalR}.

\subsection{Horizontal concatenation}
\label{sec:Horzcat}
Horizontal concatenation of two or more matrices with different sparsity structure is a common pattern in many applications~\cite{Agarwal2010,Cashman2013,Taylor2016,Triggs2000, Zach2014}.  
The core computational unit is the concatenation of two blocks, with 
matrix $A_1 \in \mathbb{R}^{\sz{n}{m_1}}$ and matrix $A_2 \in \mathbb{R}^{\sz{n}{m_2}}$:
\be
A = \hcat{A_1}{A_2}
\ee
We assume again that $A_1$ has a structure which makes 
\[
Q_1R_1 = \qr{A_1}
\]
efficient to compute.  Rewriting $A$ using full size $\Qfull_1$ gives
\begin{align*}
A & = \hcat{Q_1 R_1}{A_2} \\
  & = \hcat{\Qfull_1 \vcat{R_1\\0}}{\Qfull_1 \Qfull_1^\top A_2}\\
  & = \Qfull_1 \hcat{\vcat{R_1\\0}}{\Qfull_1^\top A_2}
\end{align*}
We now have the product of a unitary matrix $\Qfull_1$ and a matrix whose top $m_1$ rows are upper triangular.  We now factorize the bottom $n-m_1$ rows of $\Qfull_1^\top A_2$, which, from~\eqref{Qfull}, is ${\Qperp_1}^\top A_2$. 
\be
Q'R' = \qr{{\Qperp_1}^\top A_2}
\ee
And following from above,
\begin{align*}
A & = \Qfull_1 \hcat{\vcat{R_1\\0}}{\vcat{Q_1^\top A_2 \\ Q' R'}}\\
  & = \Qfull_1 \hcat{\vcat{R_1\\0}}{\vcat{I & 0\\0 & \Qfull'} \vcat{Q_1^\top A_2 \\ R' \\ 0}}\\
  & = \Qfull_1 \vcat{I & 0\\0 & \Qfull'} \hcat{\vcat{R_1\\0}}{\vcat{Q_1^\top A_2 \\ R' \\ 0}}
\end{align*}
An implementation may form the product, or simply store $Q_1$ and $Q'$ in $\code Q$.

\begin{figure}[t]
\centering
\subfigure[$A$]
{
\includegraphics[width=17mm,height=32mm]{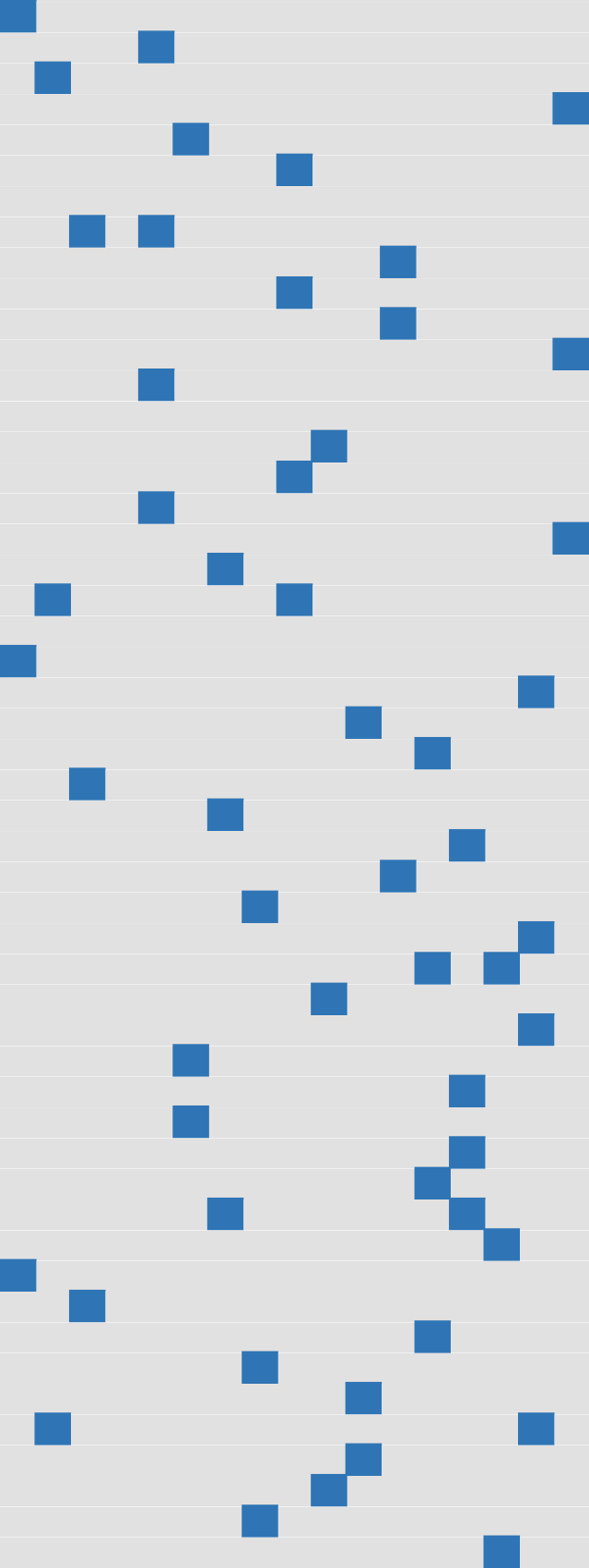}
\label{fig:rowpermBefore}
}
\subfigure[$P_rA$]
{
\includegraphics[width=19mm,height=32mm]{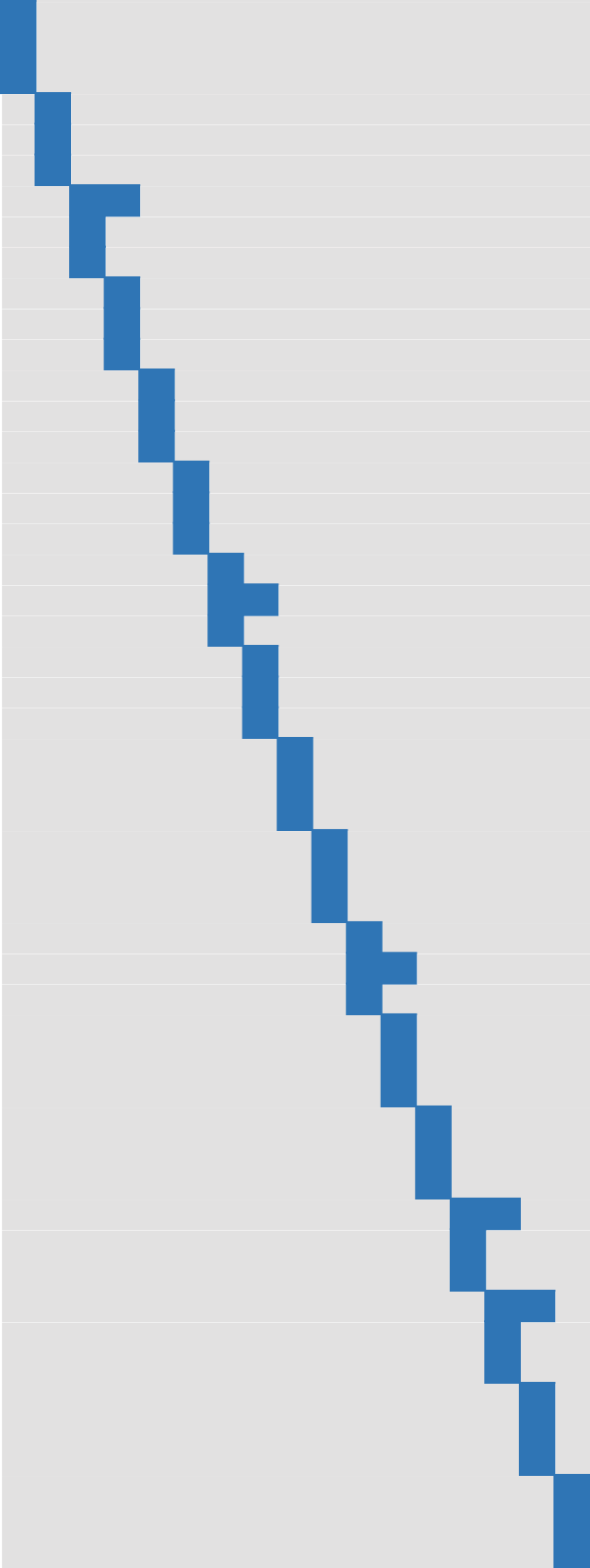}
\label{fig:rowpermAfter}
}
\subfigure[$B$]
{
\includegraphics[width=17mm,height=32mm]{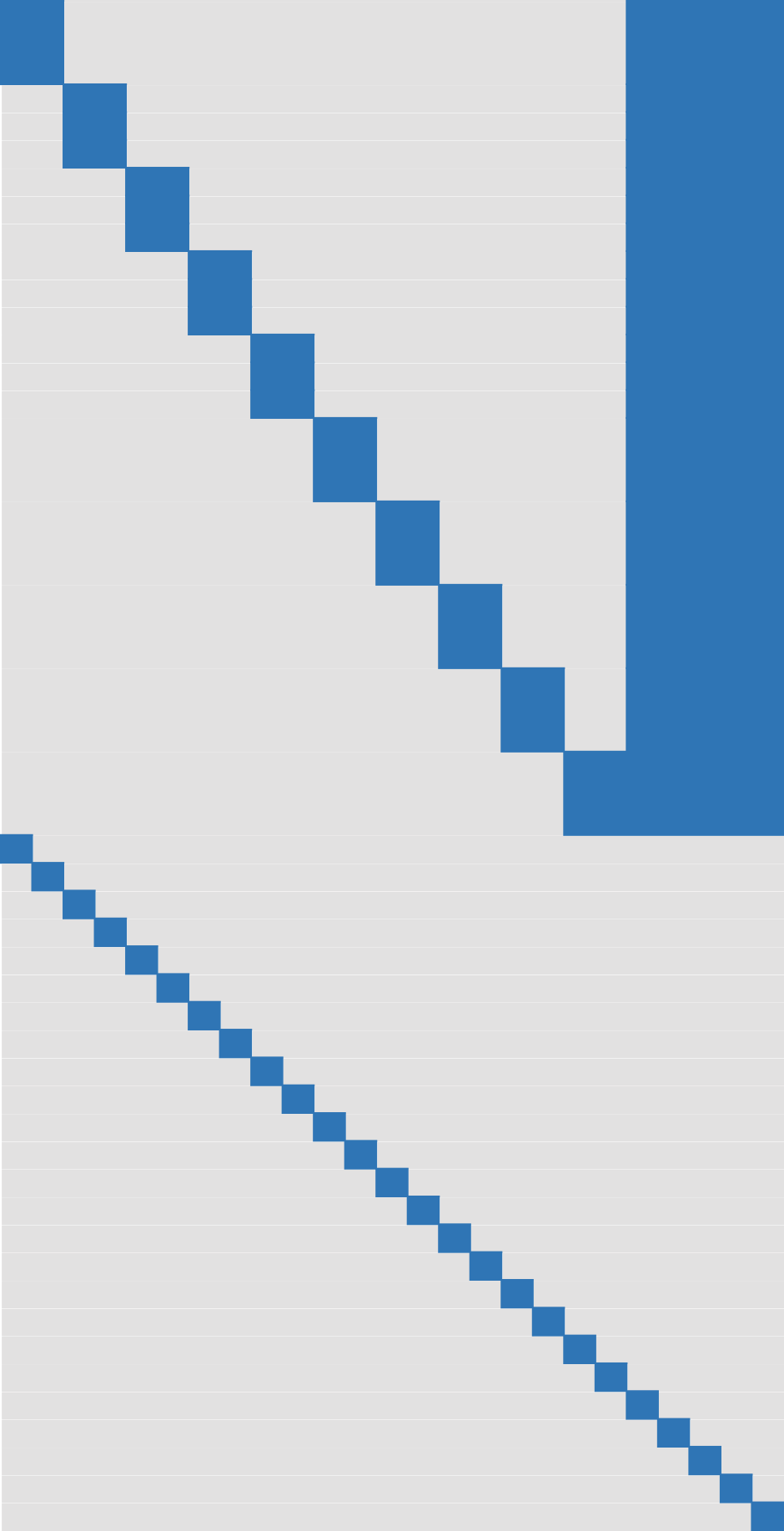}
\label{fig:levMarqJD}
}
\subfigure[$P_r B$]
{
\includegraphics[width=17mm,height=32mm]{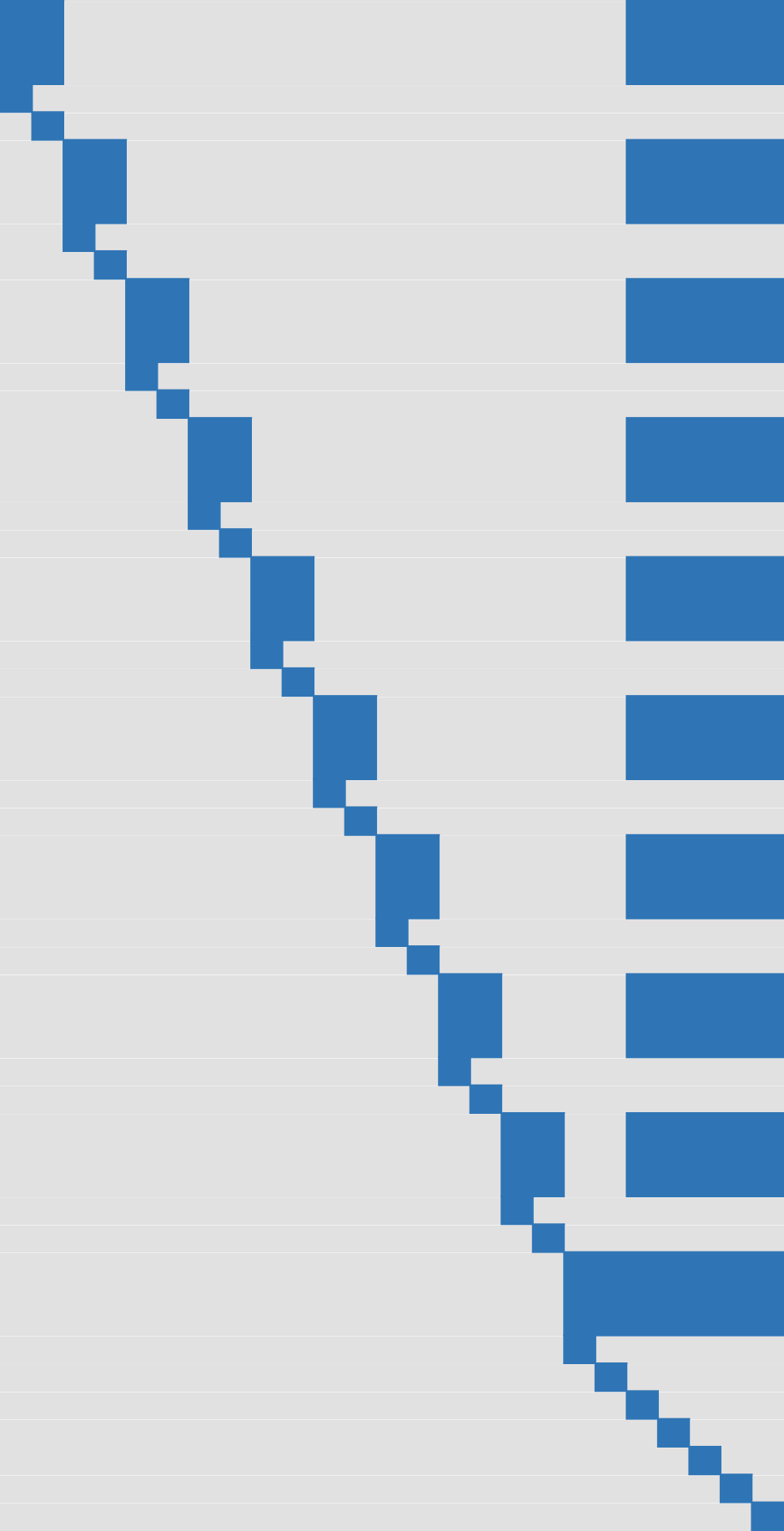}
\label{fig:levMarqJDPerm}
}
\caption{ (a-b) Row permutation $P_r$ discovering banded structure in the matrix $A$. (c-d) Row permutation $P_r$ while solving vertical concatenation of two matrices. }
\label{fig:rowPermutations}
\end{figure}

\subsection{Row and column permutations}
A common sparsity pattern is not always obvious. This can be demonstrated by generating a block diagonal/banded matrix $A \in \mathbb{R}^{\sz{n}{m}}$ and applying a random permutation $P_r \in \{0, 1\}^{\sz{n}{n}}$ to its rows (see Figure \ref{fig:rowpermBefore} and \ref{fig:rowpermAfter}).
Applying this process in reverse to a sparse input matrix $A$, we can search for a row permutation $P_r$ that would reorder the rows in order to create an `As-Banded-As-Possible' sparsity pattern~\cite{Davis2004}. The resulting matrix $A^{\prime} \in \mathbb{R}^{\sz{n}{m}}$:
\be
A^{\prime} = P_r A
\ee 
can be factorized using our efficient solvers.

Another common technique is to permute the columns of $A$ using a permutation matrix $P_c \in \{0, 1\}^{\sz{m}{m}}$, obtaining 
\be
A^{\prime} = A P_c,
\ee 
which is often used in order to reduce fill-in during the QR decomposition. 

The best practice is to combine both permutations by searching for a row-banded structure in $A$ and reducing the fill-in of the QR decomposition at the same time. The QR solver will then be performing decomposition of  
\be
A^{\prime} = P_r A P_c.
\ee

\subsection{Block banded}
\label{sec:BlockBanded}


Adding new residuals to a block-diagonal optimization problem may create overlaps between the diagonal blocks of the Jacobian and therefore break the block diagonal structure of $A$ described in Section \ref{sec:BlockDiag}. 

Given block banded $A \in \mathbb{R}^{\sz{n}{m}}$, such as $A_2$ from Figure~\ref{fig:exampleA}, let $\blk{k}{(A)} \in \mathbb{R}^{\sz{n_k}{m_k}}$ be the $k$th block of $A$. Writing $r_k \in \{0, \cdots, \textit{min}(m_k, m_{k+1})\}$ for the number of overlapping columns of blocks $\blk{k}{(A)}$ and $\blk{k+1}{(A)}$, we assume $r_k \ll m$.  

\begin{figure}[t]
\centering
\subfigure[Blocks Y]
{
\includegraphics[width=15mm,height=32mm]{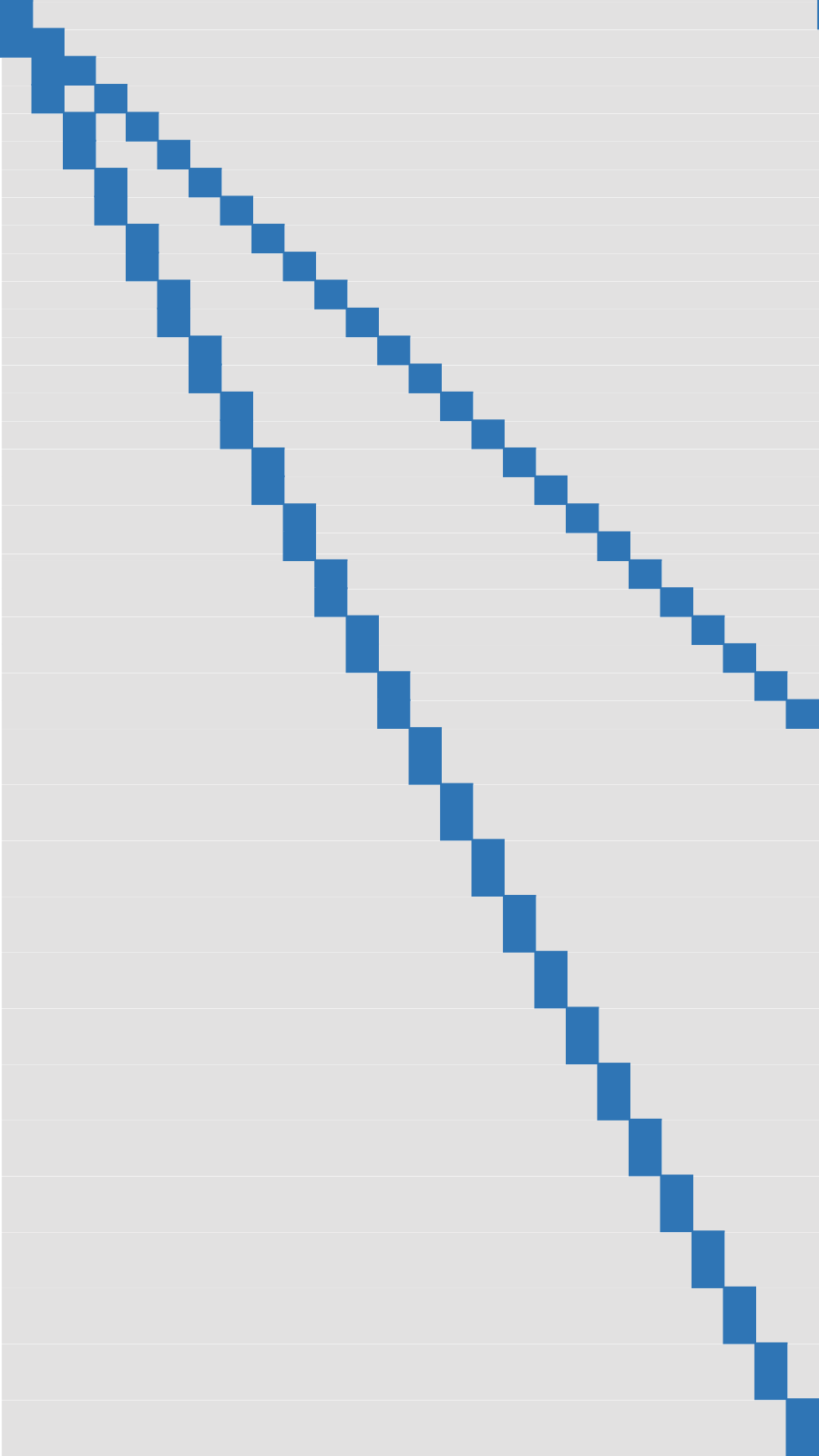}
\label{fig:blockBandedMatY}
}
\subfigure[Blocks T]
{
\includegraphics[width=15mm,height=32mm]{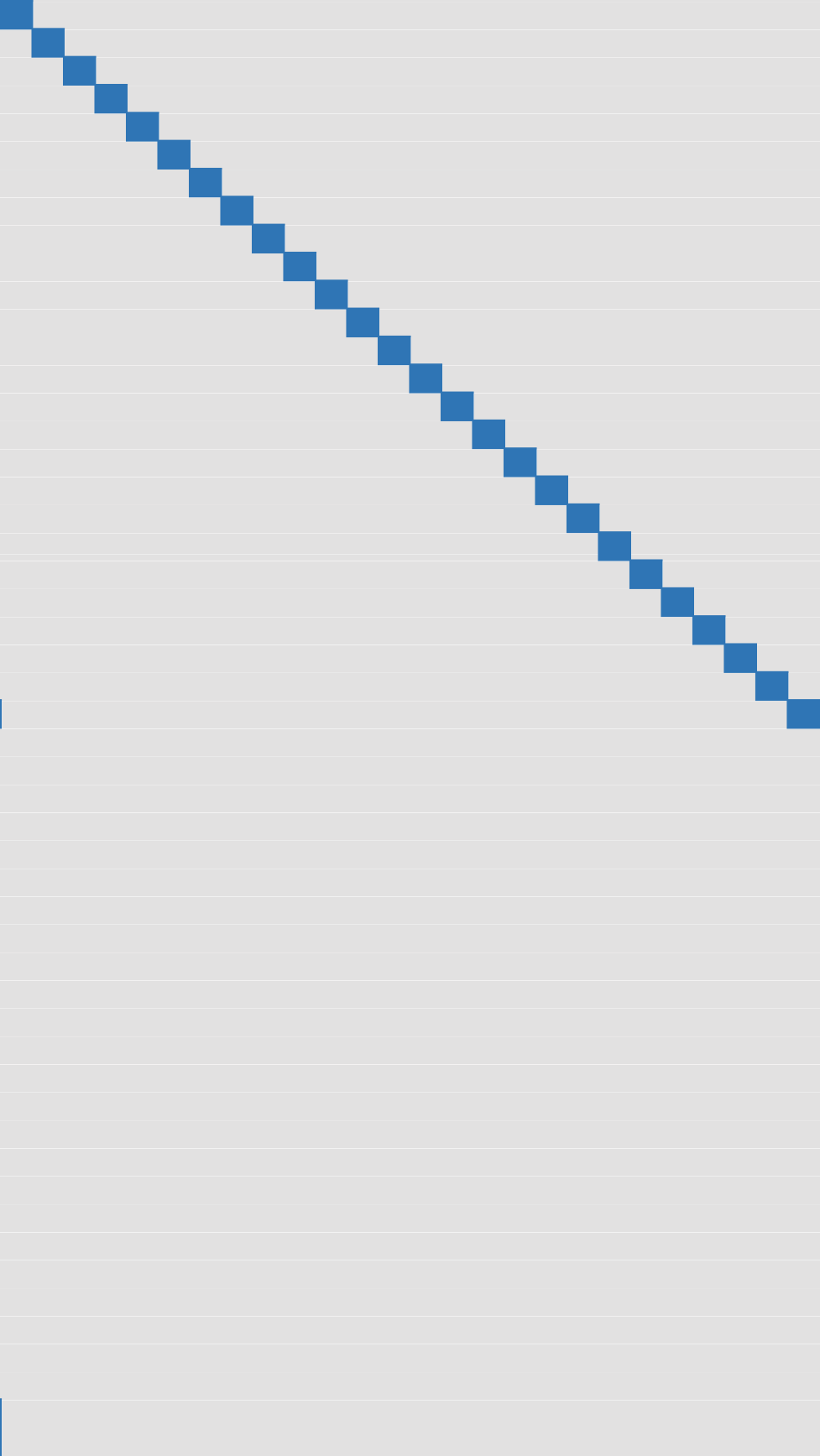}
\label{fig:blockBandedMatT}
}
\subfigure[Matrix Q]
{
\includegraphics[width=26mm,height=32mm]{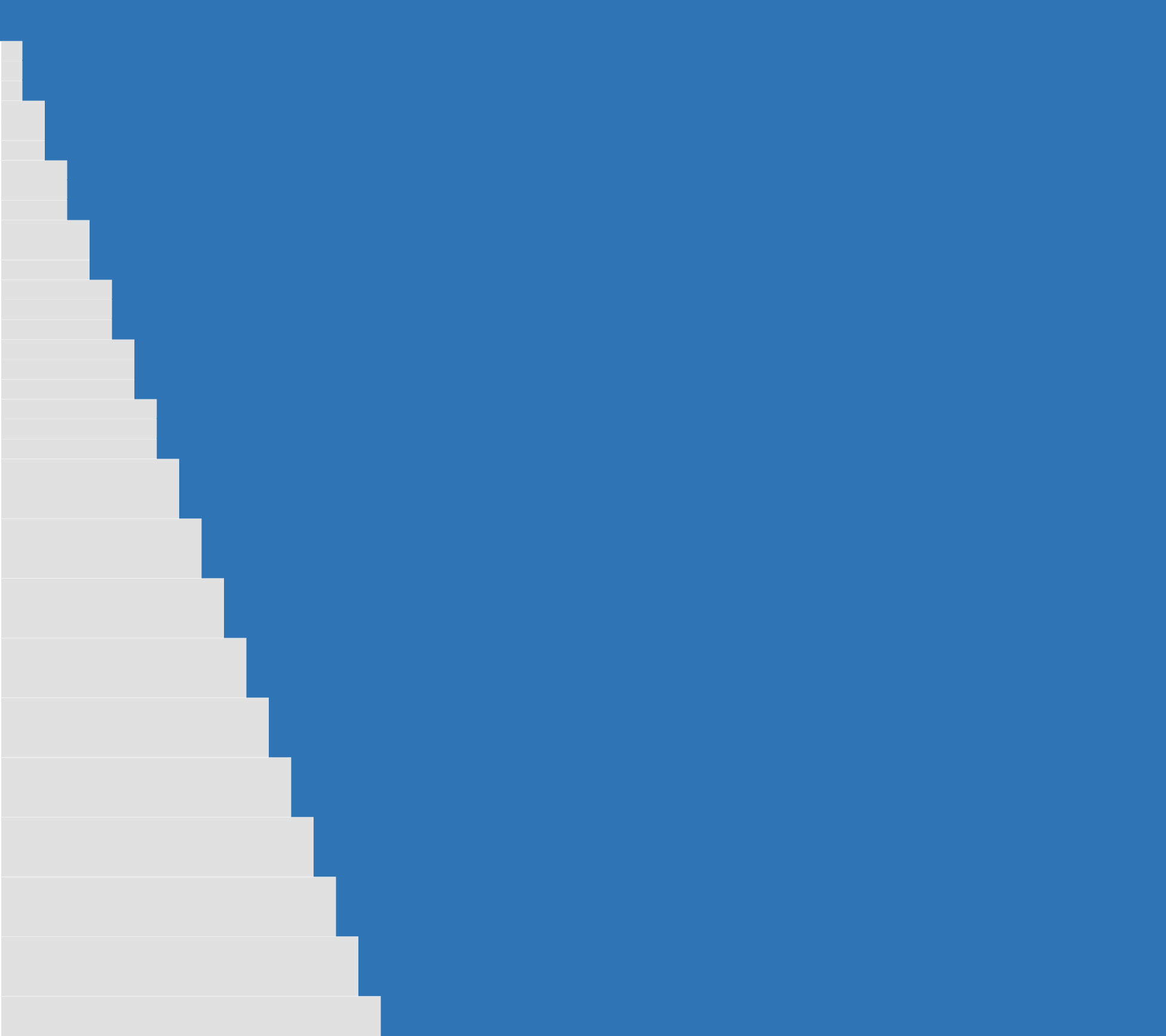}
\label{fig:blockBandedMatQ}
}
\subfigure[Matrix R]
{
\includegraphics[width=15mm,height=32mm]{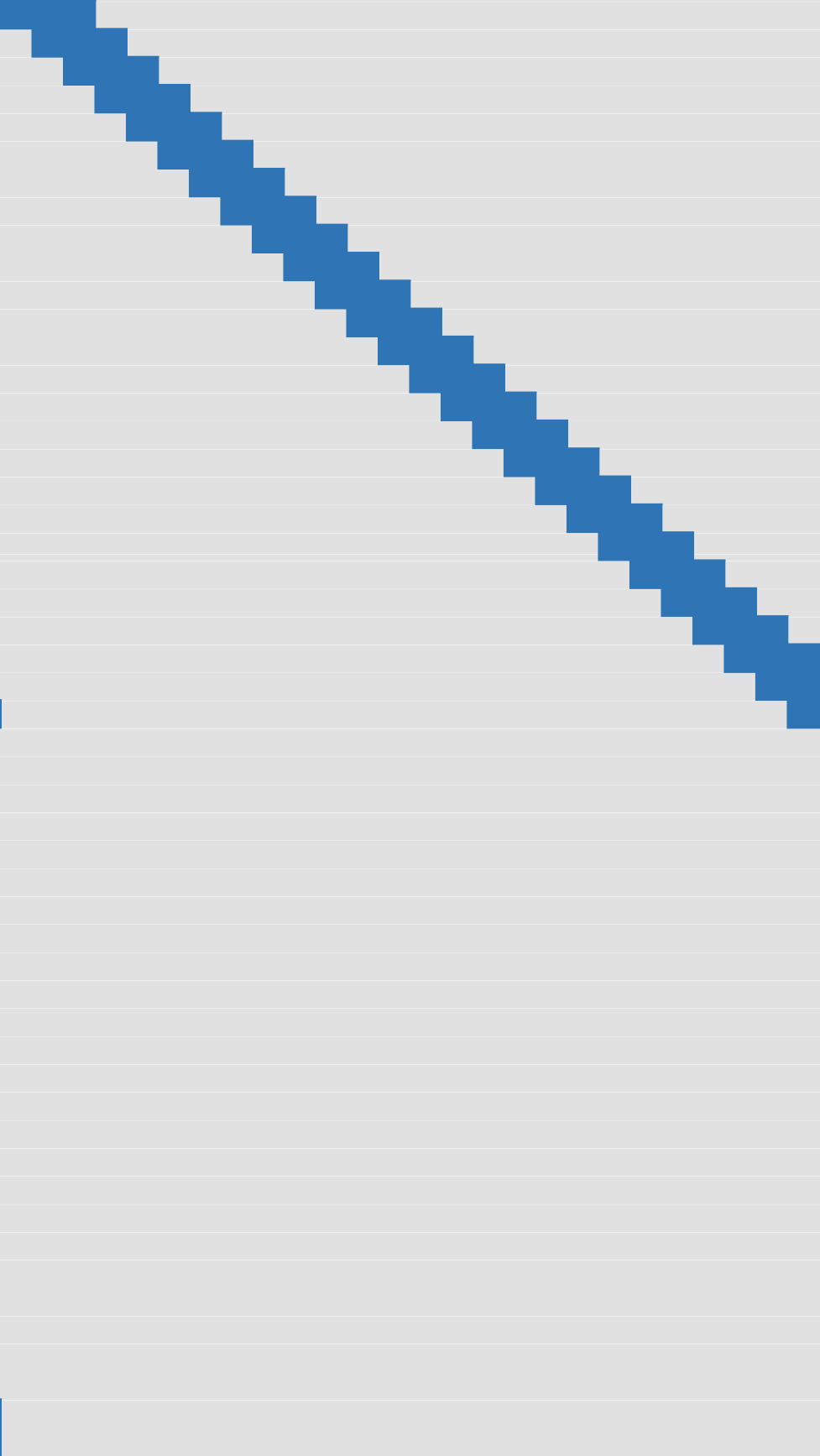}
\label{fig:blockBandedR}
}
\caption{Sparsity patterns of the QR factorization of a block banded matrix. Matrices of blocks Y and T store, column-wise using sparse representation, the small blocks Y and T that are used to factorize several columns of the input matrix at a time. }
\label{fig:blockBanded}
\end{figure}

Instead of performing QR factorization of the whole $A$, which would yield Householder vectors of length $n$, we can create block-wise Householder vectors of length $n_k + n_{k+1}$ for the $k$th block. Typically $n_k + n_{k+1} \ll n$, which results in both faster execution and lower memory consumption. The sparse Householder vectors can be stored as columns of a big sparse matrix as depicted in Figure \ref{fig:blockBandedMatY}. 

Exploiting structural zeros in such a way allows us to perform operations on small dense blocks and sequentially combine the partial results, instead of performing operations on the large sparse matrix $A$, which would suffer from both the size of $A$ and the fact that operations for sparse matrices cannot be implemented as efficiently as the dense case.

Striving for better performance, we use the compressed WY representation of Householder QR (see Section \ref{sec:BackgroundQR} for details).
Using dense blocks of size $\sz{n_j}{r}$, we factorize $r$ columns of $A$ in the $j$th step as
\be
\begin{split}
A_{j+1}^{\prime} = Q_j^\top A_j & = (I + Y_jT_jY_j^\top)^\top A_j \\
 & = A_j + Y_j (T_j^\top (Y_j^\top A_j)),
\end{split}
\ee
where $Y \in \mathbb{R}^{n_j \times r}$, $T \in \mathbb{R}^{r \times r}$ and $r < n_j$.
Factorization of $A$ is expressed as a sequence of $K$ such economy blocked transformations in form
\be
\label{eq:blockBandedYTY}
A^{\prime} = (I + Y_{k-1}T_{k-1}Y_{k-1}^\top)^\top \cdots (I + Y_0T_0Y_0^\top)^\top A,
\ee
where $k \in \{0, \dots, K-1\}$. It can be observed that matrix $Q$ expressed explicitly would be at least $60\%$ dense in this case. On the other hand, representing $Q$ as a series of the economy blocked transformations is extremely sparse. A pictorial example of the factorization is in Figure \ref{fig:blockBanded}. 

\subsection{Vertical concatenation}
\label{sec:Vertcat}
In general, vertical concatenation of two QR decompositions is much harder than the horizontal concatenation discussed in Section \ref{sec:Horzcat}. The main unit is the two-block case, having $A_1 \in \mathbb{R}^{n_1 \times m}$ and $A_2 \in \mathbb{R}^{n_2 \times m}$ with
\be
A = \vcat{A_1 \\ A_2}.
\ee
Assuming both $A_1$ and $A_2$ have a favorable sparsity pattern, we can perform their QR decompositions $Q_1 R_1 = \qr{A_1}$ and $Q_2 R_2 = \qr{A_2}$ efficiently and express the result as
\be
\label{eq:vcatDecomp}
A = \vcat{Q_1 R_1 \\ Q_2 R_2} = \vcat{Q_1 & 0 \\ 0 & Q_2} \vcat{R_1 \\ R_2} = Q^{\prime} R^{\prime}.
\ee
It is apparent from (\ref{eq:vcatDecomp}) that the matrix $R^{\prime}$ is not upper triangular and therefore we do not have a valid QR decomposition of $A$ yet.
One possibility would be to use Givens rotations or column-wise Householder reflectors to eliminate the lower block $R_2$.
However, if $A_1$ and $A_2$ have known sparsity structure, we can again exploit the known structure of $R^{\prime}$ for greater efficiency.
In particular, if $R_1$, the upper part of $R$, is output by one of the QRkit sparse solvers, then it is typically very sparse with elements accumulated in the proximity of the diagonal (see Figure \ref{fig:blockDiagonalR} and \ref{fig:blockBandedR}). We can therefore apply a row permutation matrix $P \in \mathbb{R}^{\sz{n_1 + n_2}{m}}$ to $R^{\prime}$ that will interleave $R_1$ and $R_2$ so that they create another block diagonal/banded matrix $P R^{\prime}$, as depicted in Figure \ref{fig:levMarqJD} and \ref{fig:levMarqJDPerm}. We have shown in the previous sections how to factorize such matrices efficiently. Following Equation \ref{eq:vcatDecomp}, the final QR decomposition of $A$ is expressed as
\be
\begin{split}
A = Q^{\prime} P R^{\prime} & = \vcat{Q_1 & 0 \\ 0 & Q_2} P R^{\prime} \\
  & = \vcat{Q_1 & 0 \\ 0 & Q_2} P Q_3 \vcat{R_3 \\ 0}
\end{split}
\ee

\section{Implementation}
\label{sec:implementation}
QRkit is implemented using only the Eigen C++ library~\cite{Eigen} without any calls to other external libraries.
\subsection{Solver interface}
New solver for a specific structure of the input matrix $A$ is defined as composition of our efficient solvers (see Section \ref{sec:qrkit}) as follows:
\begin{small}
\begin{verbatim}
 typedef BlockDiagonalQR<DenseQR> LeftSolver;
 typedef BlockBandedQR<DenseQR> RightSolver;
 HorzCat<LeftSolver, RightSolver> slvr;
\end{verbatim}
\end{small}
A typical use case, solution of the least squares problem $Ax = b$, is represented by the following pseudocode:
\begin{small}
\begin{verbatim}
 slvr.compute(A);
 qtb = slvr.matrixQ().transpose() * b;
 x = slvr.matrixR().solve(qtb);
\end{verbatim}
\end{small}
where \code{slvr} is an efficient solver specified using the definitions above and it therefore knows how to deal with the sparsity pattern of $A$ efficiently.

\subsection{Levenberg--Marquardt}
In order to solve non-linear optimization problems, we are using QRkit as the core building block of the Levenberg--Marquardt algorithm. The Eigen C++ library already provides an implementation of Levenberg--Marquardt using QR solvers. It is a C++ port of LAPACK Fortran code which implements the Mor{\'e} method~\cite{More1977}. 

\paragraph{Mor{\'e} Levenberg--Marquardt} performs two QR decompositions in each step. Starting from (\ref{eq:LevMarqQR}), each iteration performs QR decomposition as follows:
\be
\label{eq:LMTwoStep}
\begin{split}
 J & = QR \\
 \vcat{R \\ \lambda^{\frac{1}{2}} D} & = Q^{\prime}R^{\prime},
\end{split}
\ee
where $J \in \mathbb{R}^{n \times m}$, $R \in \mathbb{R}^{n \times m}$ is upper triangular, $D \in \mathbb{R}^{m \times m}$ is identity and $\lambda$ is the damping factor.

Unfortunately, Eigen's implementation solves the second QR decomposition as a sequence of Givens rotations, which becomes slow as number of columns $m$ increases.

However, the solver described in Section \ref{sec:Vertcat} is applicable for this vertical concatenation of two matrices with favorable sparsity patterns.
If we consider the most general case, where matrix $R$ in (\ref{eq:LMTwoStep}) has a dense upper triangle, permuting rows of $\lambda^{\frac{1}{2}}D$ into $R$ will create a skewed upper-triangular structure that can be treated with sparsity-aware blocked Householder QR, which is faster than applying Givens transformations.

We should however remind ourselves that for a Jacobian matrix $J$ with favorable sparsity pattern, the upper triangular matrix $R$ of its QR decomposition will still be very sparse, with elements concentrated close to the diagonal (see Figures \ref{fig:blockDiagonalR} and \ref{fig:blockBandedR}). Row permuting the diagonal matrix $\lambda^{\frac{1}{2}}D$ into such an $R$ will result in a block banded matrix which we know how to solve very efficiently (see Figure \ref{fig:levMarqJD} and \ref{fig:levMarqJDPerm}).

\paragraph{Backtrack Levenberg--Marquardt}
Assuming favorable sparsity pattern of the Jacobian $J$, we can reduce the required number of QR factorizations by directly row-permuting the diagonal matrix $\lambda^{\frac{1}{2}}D$ into $J$ and performing a single QR decomposition
\be
P_r \vcat{J \\ \lambda^{\frac{1}{2}} D} = QR,
\ee
where $J \in \mathbb{R}^{n \times m}$, $\lambda$ is the damping factor, $P_r$ the row permutation matrix and $D \in \mathbb{R}^{m \times m}$ is the identity matrix for the Levenberg algorithm, or $\mathrm{diag}(J^\top J)^\frac{1}{2}$ for Levenberg--Marquardt.

We further speed up the Levenberg--Marquardt iterations by adapting the approach of Lourakis \etal~\cite{Lourakis2009} described at the end of Section \ref{sec:BackgroundQR}. 
 
\section{Results}
We compare QRkit to some state-of-the-art methods on two common computer vision problems: surface fitting and bundle adjustment. We perform all our experiments using both single and double precision floating point in order to evaluate the impact of differing machine precision on the convergence. We assess accuracy by comparing the optimizers purely on their ability to minimize the objective function in question, not by comparison to any ground truth for these problems, since this is the most direct evaluation of the optimizers' success.

\subsection{Solvers}
Each experiment compares several different QR and Cholesky-based solvers described below. The solvers were used as the core part of the Backtrack Levenberg--Marquardt implementation from Section \ref{sec:implementation}, with the exception of SSBA, which is standalone and Mor{\'e} QR, which operates in two steps.
\vspace{-2mm}
\paragraph{QRkit}
Our new kit of sparse QR factorizations directly implemented as a submodule of the Eigen C++ library.
\newline
\noindent\rule{0pt}{5mm}\textbf{Eigen Sparse QR}
Current implementation of Sparse QR solver in the Eigen C++ library. It is a fill-in reducing rank-revealing QR factorization which does not leverage any sparsity structure in the input matrix. 
\newline
\noindent\rule{0pt}{5mm}\textbf{SuiteSparse QR (SPQR)}
SuiteSparse implementation of multifrontal rank-revealing QR factorization~\cite{Davis2011} using BLAS~\cite{Lawson1979} and LAPACK~\cite{Anderson1990}. It does not assume any sparsity pattern in advance.
\newline
\noindent\rule{0pt}{5mm}\textbf{Cholesky}
Variant of Cholesky factorization, in particular \\ \code{Eigen::SimplicialLDLT}, which is a sparse fill-in reducing LDLT Cholesky factorization without square root. 
\newline
\noindent\rule{0pt}{5mm}\textbf{QRkit + Cholesky}
A combination of the QRkit block diagonal solver and a Cholesky solver for decomposing block angular matrices. The fast QRkit block diagonal solver is applied on the left block diagonal part, and the generally dense right subblock is consecutively solved using Cholesky.
\newline
\noindent\rule{0pt}{5mm}\textbf{SSBA}
Complete bundle adjustment package provided by Zach~\cite{Zach2014} for large sparse bundle adjustment problems. We compare our methods to SSBA's Cholesky-based Levenberg--Marquardt bundle optimizer.
\newline
\noindent\rule{0pt}{5mm}\textbf{Mor{\'e} QR} Same as QRkit, however in each Levenberg--Marquardt iteration, the factorization is done in two steps (see Section \ref{sec:implementation}).

\begin{figure}[t]
\centering
\subfigure[Ellipse fitting]
{
\includegraphics[width=22mm,height=32mm]{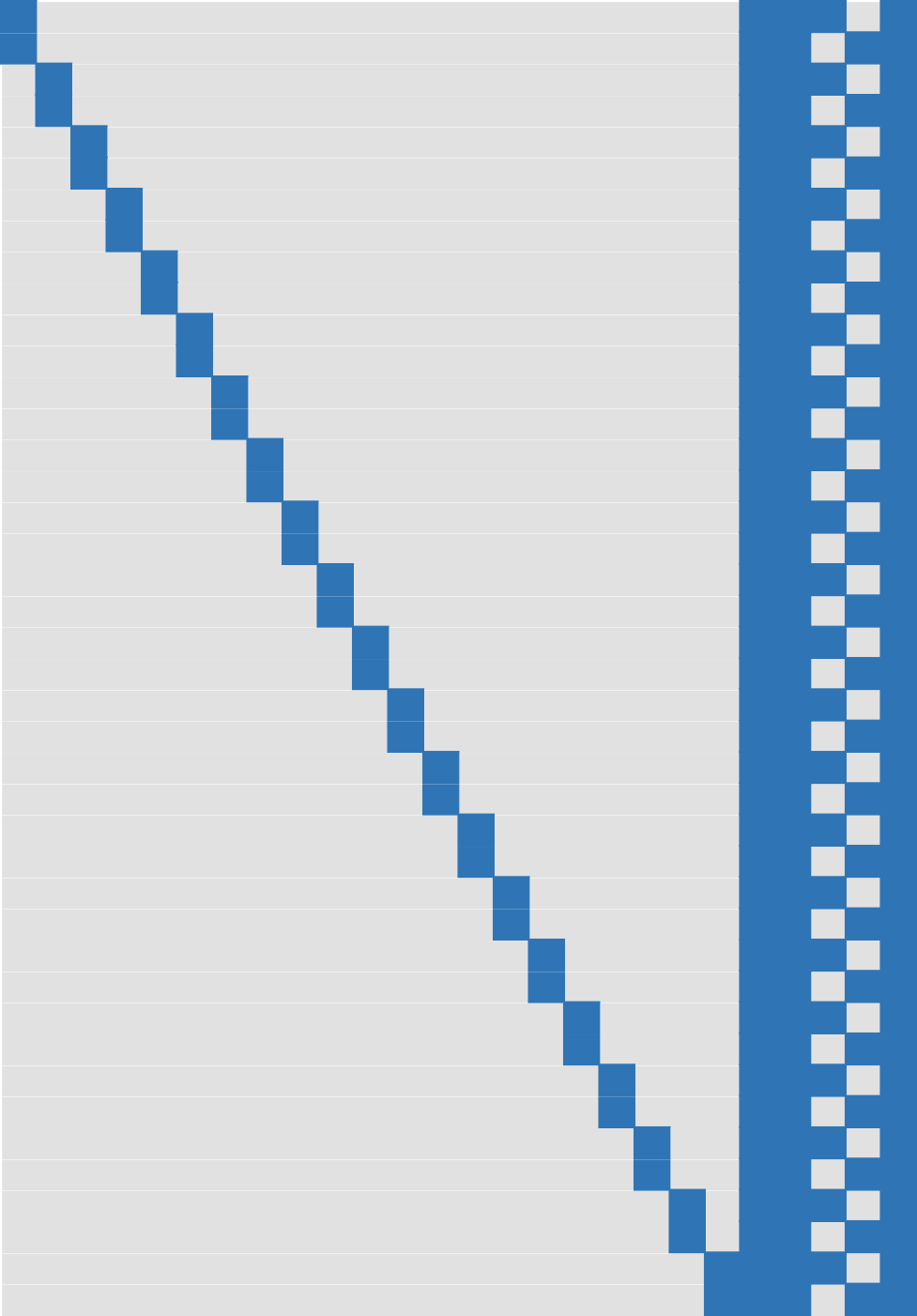}
\label{fig:ellipseFittingJac}
}
\hspace{5mm} 
\subfigure[Bundle adjustment]
{
\includegraphics[width=24mm,height=32mm]{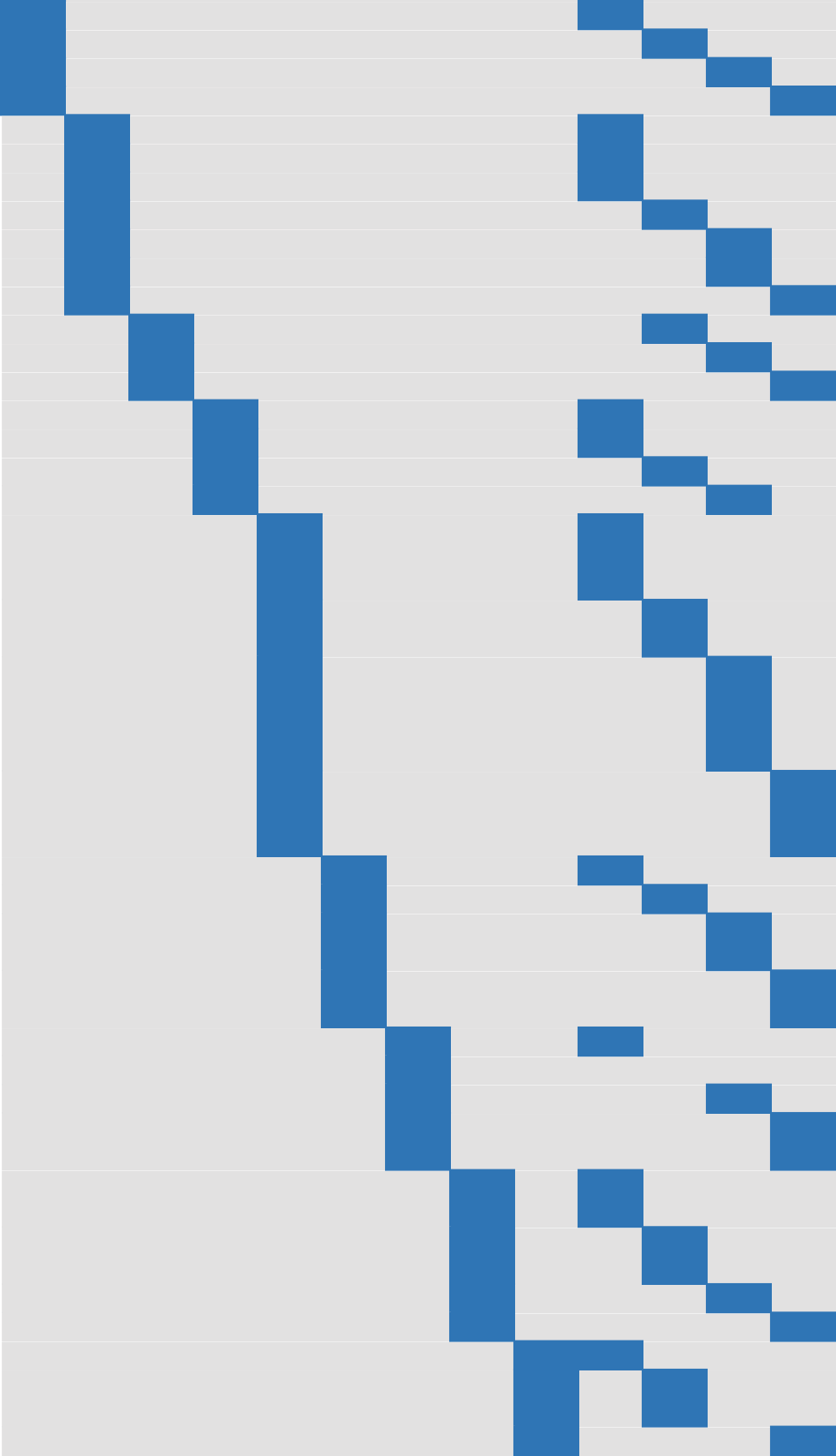}
\label{fig:bundleAdjustJac}
}
\caption{Sparsity patterns of the Jacobian matrix for our benchmark problems.}
\end{figure}

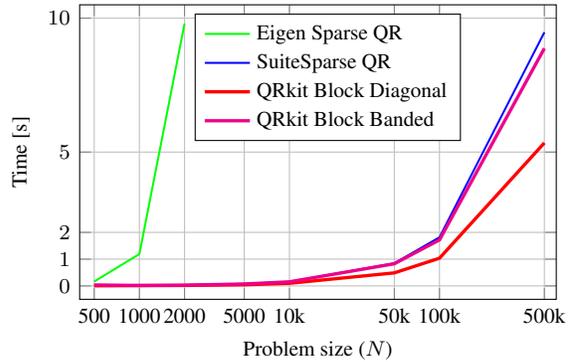
\begin{figure}[ht]
\centering
\setlength\figureheight{5.5cm} 
\setlength\figurewidth{0.95\linewidth}
\input{ellipse_fitting.tex}
\caption{Performance of different QR factorization methods on the ellipse fitting benchmark.}
\label{fig:ellipseBench}
\end{figure}

\begin{table}[ht]
\begin{small}
\begin{center}
\begin{tabular}{ccccc}
\toprule
\multirow{2}{1.5cm}{\centering\rule{0pt}{4mm} Problem size} & \multicolumn{4}{c}{Time [s]} \\
  \cmidrule{2-5}
& EigSpQR & SPQR & \textbf{QRkitBD} & \textbf{QRkitBB} \\
\midrule
500 & 0.163 & 0.016 & 0.005 & 0.037 \\
2,000 & 9.798 & 0.031 & 0.017 & 0.029 \\
10,000 & --- & 0.151 & 0.098 & 0.154 \\
100,000 & --- & 1.816 & 1.036 & 1.718 \\
500,000 & --- & 9.472 & 5.342 & 8.872 \\
\bottomrule
\end{tabular}
\end{center}
\end{small}
\caption{Timings for different methods on the ellipse fitting benchmark; EigSpQR=Eigen Sparse QR, SPQR=SuiteSparseQR, QRkitBD=QRkit Block Diagonal, QRkitBB=QRkit Block Banded. For Eigen Sparse QR, results are available only up to $N = 2000$ as the computation takes too long for bigger problems. The superiority of QRkit is clearly visible.}
\label{tab:ellipseBench}
\vspace*{-3mm}
\end{table}
\subsection{Experiments}

\paragraph{Surface fitting} is a popular technique for explaining unknown data by fitting a parametric model. In order to show the suitability of QRkit for these problems, we have implemented a simple 2D optimization that fits an ellipse to a set of 2D points. It can be considered a simplification of surface fitting tasks such as human body tracking~\cite{PerezSala2014} or hand tracking~\cite{Tagliassacchi2015,Taylor2016}, which have recently received a great deal of interest. The structure of the Jacobian is depicted in Figure~\ref{fig:ellipseFittingJac}. An efficient solution can be obtained by expressing $J$ as a horizontal concatenation, as described in Section \ref{sec:Horzcat}.

We performed the experiment at different scales. For the number of 2D points $N$, we have performed evaluations for $N$ ranging from $500$ up to $500,000$. Figure~\ref{fig:ellipseBench} shows the comparison of different solvers; the exact timings are then listed in Table~\ref{tab:ellipseBench}. We observe that for very small problems ($N \le 2000$), QRkit performance is comparable to existing methods. As we scale $N$ up to $500,000$, however, QRkit significantly outperforms state-of-the-art implementations, especially when using the block diagonal solver, but also (with a smaller margin) when a block banded solver is used.

\begin{figure}[ht!]
\centering
\subfigure
{
\setlength\figureheight{5cm} 
\setlength\figurewidth{0.95\linewidth}
\input{trafalgar_iter_double.tex}
}
\addtocounter{subfigure}{-1}
\subfigure
{
\setlength\figureheight{4cm} 
\setlength\figurewidth{0.95\linewidth}
\input{trafalgar_double.tex}
}
\caption{Trafalgar (double precision): Convergence comparison of different methods in terms of execution time and number of iterations on the Trafalgar bundle adjustment dataset.}
\label{fig:baBenchTrafalgarDouble}
\end{figure}
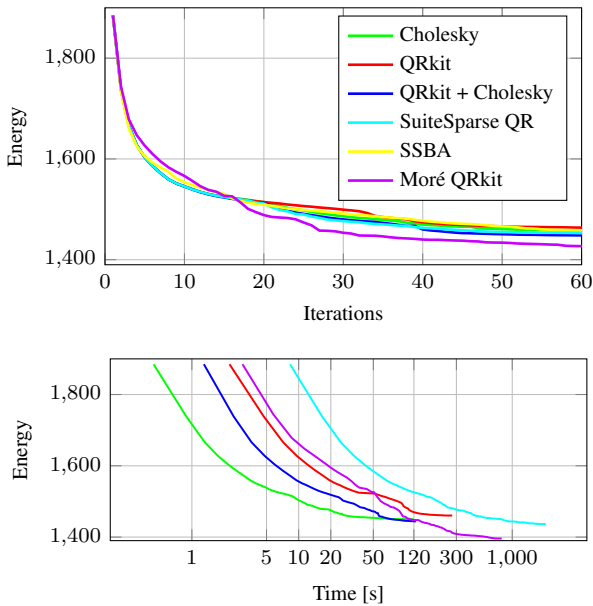

\paragraph{Bundle adjustment} is a simultaneous refinement of 3D coordinates in order to describe geometry of a scene observed by multiple cameras with unknown parameters \cite{Agarwal2010,Triggs2000,Zach2014}. The general structure of the Jacobian in a bundle adjustment problem is sketched in Figure \ref{fig:bundleAdjustJac}. Its angular structure is efficiently solvable if expressed as a horizontal concatenation (see Section \ref{sec:Horzcat}).

We perform this experiment using standard datasets from GRAIL~\cite{Agarwal2010}. In particular, we selected the small version of the and Dubrovnik and Trafalgar Square datasets (the latter we simply call `Trafalgar' below). Each camera is represented by 9 parameters, the 3D points by \textit{x,y,z-}coordinates and measurements by \textit{x,y-}coordinates.

\textbf{Trafalgar} consists of 21 cameras capturing 11,315 3D points with a total of 36,455 2D observations. This gives us Jacobian $J \in \mathbb{R}^{n \times m}$ with $n = 36455 \times 2 = 72910$ rows and $m = 21 \times 9 + 11315 \times 3 = 34134$ columns.

\textbf{Dubrovnik} has only 16 cameras, which capture 22,106 3D points producing 83,718 2D observations. In this case, the Jacobian $J \in \mathbb{R}^{n \times m}$  has $n = 83718 \times 2 = 167436$ rows and $m = 16 \times 9 + 22106 \times 3 = 66462$ columns.

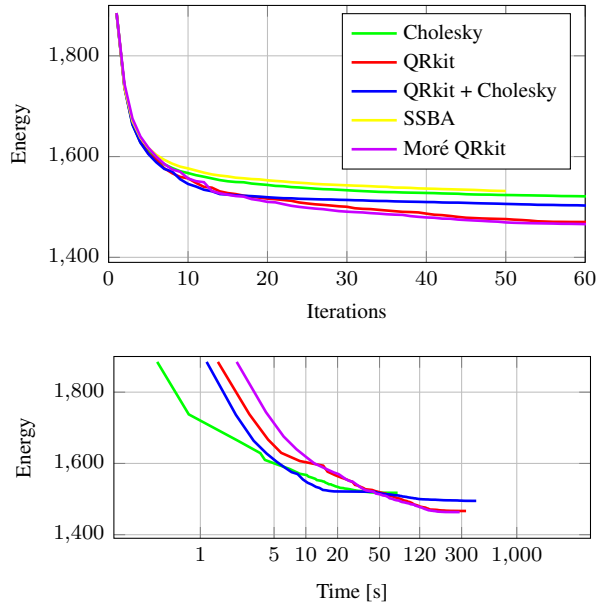
\begin{figure}[ht!]
\centering
\subfigure
{
\setlength\figureheight{5cm} 
\setlength\figurewidth{0.95\linewidth}
\input{trafalgar_iter_float.tex}
}
\addtocounter{subfigure}{-1}
\subfigure
{
\setlength\figureheight{4cm} 
\setlength\figurewidth{0.95\linewidth}
\input{trafalgar_float.tex}
}

\caption{ Trafalgar (single precision): Convergence comparison of different methods in terms of execution time and number of iterations on the Trafalgar bundle adjustment dataset.}
\label{fig:baBenchTrafalgarSingle}
\end{figure}

\begin{figure}[ht!]
\centering
\subfigure
{
\setlength\figureheight{5cm} 
\setlength\figurewidth{0.95\linewidth}
\input{dubrovnik_iter_double.tex}
}
\addtocounter{subfigure}{-1}
\subfigure
{
\setlength\figureheight{4cm} 
\setlength\figurewidth{0.95\linewidth}
\input{dubrovnik_double.tex}
}
\caption{Dubrovnik (double precision): Convergence comparison of different methods in terms of execution time and number of iterations on the Dubrovnik bundle adjustment dataset. }
\label{fig:baBenchDubrovnikDouble}
\end{figure}
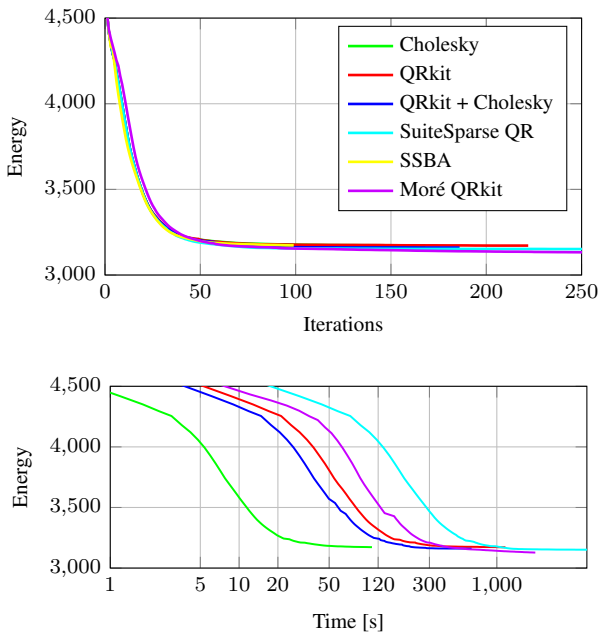

Evaluations were performed for both single (32-bit) and double (64-bit) precision floating point.  The convergence and execution times in double precision are depicted in Figures~\ref{fig:baBenchTrafalgarDouble} and~\ref{fig:baBenchDubrovnikDouble}. As expected, Cholesky factorization is the fastest here. On the other hand, it does not reach as good an optimum as the Mor{\'e} QRkit implementation, which suggests better numerical stability of QR for the cost of slower execution. It is worth mentioning that QRkit outperforms the state-of-the-art QR factorization implementation SPQR in terms of execution time and performs on-par with it in terms of convergence (see Table~\ref{tab:baBenchConvergence}).

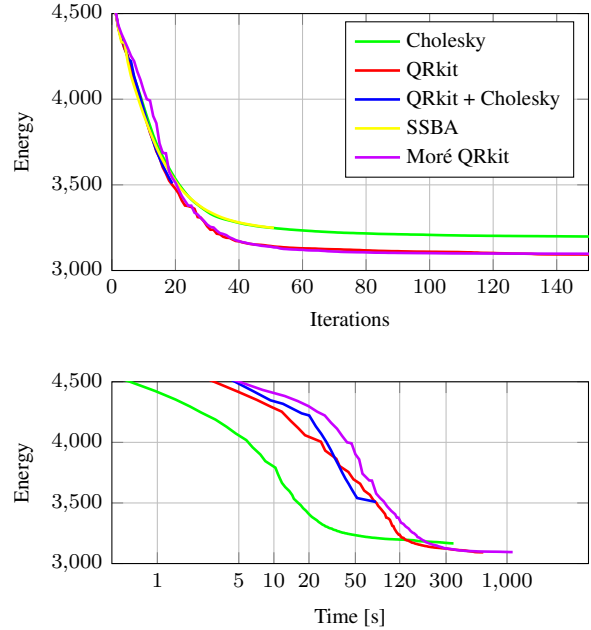
\begin{figure}[ht!]
\centering
\subfigure
{
\setlength\figureheight{5cm} 
\setlength\figurewidth{0.95\linewidth}
\input{dubrovnik_iter_float.tex}
}
\addtocounter{subfigure}{-1}
\subfigure
{
\setlength\figureheight{4cm} 
\setlength\figurewidth{0.95\linewidth}
\input{dubrovnik_float.tex}
}
\caption{Dubrovnik (single precision): Convergence comparison of different methods in terms of execution time and number of iterations on the Dubrovnik bundle adjustment dataset. }
\label{fig:baBenchDubrovnikSingle}
\end{figure}

Results for single precision in terms of convergence and execution times are shown in Figures~\ref{fig:baBenchTrafalgarSingle} and~\ref{fig:baBenchDubrovnikSingle}. It suggests that Cholesky might not be the right choice in this case as it struggles to find a good optimum. On the contrary, QRkit finds an optimum closer to the one achieved in double precision and executes approximately 50\% faster. This shows the strong advantage of using QRkit over Cholesky in single precision arithmetic.

Numerical results for both single and double precision are listed in Table~\ref{tab:baBenchTiming}. In addition, Table~\ref{tab:baBenchConvergence} displays the minimum energy for each of the methods and both precisions.



\begin{table}[!h]
\begin{small}
\begin{center}
\begin{tabular}{ccccc}
\toprule
\multirow{3}{*}{\rule{0pt}{5mm}Method} & \multicolumn{4}{c}{Minimum energy} \\
  \cmidrule{2-5}
& \multicolumn{2}{c}{Trafalgar} & \multicolumn{2}{c}{Dubrovnik}\\
  \cmidrule{2-5}
 & Double & Single & Double & Single \\
\midrule
Cholesky & 1450.34 & 1517.57 & 3172.60 & 3166.51 \\
\textbf{QRkit} & 1460.43 & 1466.59 & 3171.77 & 3092.92 \\
\textbf{QRkit + Chol} & 1444.10 & 1494.93 & 3160.13 & 3508.25 \\
\textbf{Mor{\'e} QRkit} & 1395.99 & 1463.50 & 3128.07 & 3094.73 \\
SPQR & 1436.58 & --- & 3151.63 & --- \\	
SSBA & 1454.19 & 1531.69 & 3172.89 & 3248.25 \\
\bottomrule
\end{tabular}
\end{center}
\end{small}
\caption{ Comparison of the energy minimum found by different methods. The numerical stability of QR is emphasized especially in single precision arithmetic. }
\label{tab:baBenchConvergence}
\vspace*{3mm}
\end{table}

\begin{table*}[ht!]
\begin{center}
\begin{small}
\begin{tabular}{@{\extracolsep{5pt}}ccccccccc@{}}
\toprule
\multirow{3}{*}{\rule{0pt}{5mm}Method} & \multicolumn{8}{c}{Time [s] / Number of iterations} \\
\cmidrule{2-9}
& \multicolumn{4}{c}{Trafalgar dataset} & \multicolumn{4}{c}{Dubrovnik dataset} \\
\cmidrule{2-5}\cmidrule{6-9}
 & Do@1600 & Si@1600 & Do@1460 & Si@1475 & Do@3500 & Si@3500 & Do@3175 & Si@3125 \\
\cmidrule{1-1}\cmidrule{2-5}\cmidrule{6-9}
Cholesky & 2.19 / 5 & 5.04 / 6 & 27.05 / 51 & --- & 11.37 / 19 & 15.93 / 21 & 73.71 / 112 & --- \\
\textbf{QRkit} & 12.08 / 5 & 9.98 / 5 & 222.25 / 85 & 138.99 / 51 & 80.16 / 19 & 73.22 / 18 & 551.63 / 116 & 287.83 / 67 \\
\textbf{QRkit + Chol} & 5.98 / 5 & 5.967 / 6 & 56.36 / 40 & --- & 58.45 / 19 & 75.15 / 19 & 216.05 / 67 & --- \\
\textbf{Mor{\'e} QRkit} & 18.14 / 6 & 12.063 / 6 & 94.81 / 27 & 124.70 / 43 & 123.11 / 21 & 87.47 / 20 & 404.27 / 62 & 298.80 / 55 \\
SPQR & 41.28 / 5 & --- & 451.36 / 43 & --- & 282.92 / 19 & --- & 1038.73 / 61 & --- \\
SSBA & 0.55 / 5 & 0.38 / 6 & 6.10 / 41 & --- & 6.99 / 19 & 4.17 / 21 & 29.79 / 81 & --- \\
\bottomrule
\end{tabular}
\end{small}
\end{center}
\caption{ Timing comparison for different factorization approaches based on both Cholesky and QR, showing the contrast between speed of Cholesky and stability of QR. The header is in the format Precision@Energy, where Do=Double and Si=Single. }
\label{tab:baBenchTiming}
\end{table*}


\section{Conclusions}
We presented a new suite of sparsity-aware QR factorizations for the Eigen C++ library. Our QRkit can efficiently deal with matrices that exhibit block diagonal or block banded sparsity patterns and their horizontal or vertical concatenations. QRkit is open source, fully contained in Eigen and does not have any external dependencies. It is therefore a good candidate to become part of the official Eigen release in the near future. We further adapt Eigen's Levenberg--Marquardt implementation to be able to evaluate our solvers on larger problems.

We performed experiments on a simple surface fitting problem and problems from the standard datasets in bundle adjustment. We show superior performance over the state-of-the-art sparse QR solver SPQR from SuiteSparse. Furthermore, we confirm the better numerical stability of QR over Cholesky decomposition in single precision arithmetic, which holds increasing importance for embedded computer vision applications.

For all tested problems, single precision reached lower energies with QR-based than with Cholesky-based optimizers, regardless of runtime, and for many problems this is a clear advantage of the QR-based methods.
These encouraging results motivate us to continue in the proposed direction and revisit the current state-of-the-art in sparse matrix factorization. In order to motivate the community to rethink the abundant use of Cholesky, our future work will show numerous applications of QRkit in latent variable problems.

{\small
\bibliographystyle{ieee}
\bibliography{egpaper_for_review}
}

\end{document}

%% file: ellipse_fitting.tex
%
%
%
\begin{tikzpicture}

\pgfplotsset{compat=newest} 

\definecolor{color0}{rgb}{1, 1, 1}

\tikzstyle{every node}=[font=\footnotesize]

\begin{axis}[
xlabel={Problem size ($N$)},
ylabel={Time [s]},
xmin=400, xmax=600000,
ymin=-0.5, ymax=10.5,
width=\figurewidth,
height=\figureheight,
at={(0\figurewidth,0\figureheight)},
xmajorgrids,
x grid style={lightgray},
ymajorgrids,
y grid style={lightgray},
axis line style={black},
axis background/.style={fill=color0},
xtick={500,1000,2000,5000,10000,50000,100000,500000},
xticklabels={500, 1000, 2000, 5000, 10k, 50k, 100k, 500k},
ytick={0, 1, 2, 5, 10},
xmode=log,
log ticks with fixed point,
legend style={at={(0.80,0.97)}, anchor=north east},
legend cell align={left},
legend entries={{Eigen Sparse QR},{SuiteSparse QR},{QRkit Block Diagonal},{QRkit Block Banded}}
]
%
%
%
\addplot [line width=0.75pt, green]
table {%
500 0.163
1000 1.191
2000 9.798
};
\addplot [line width=0.75pt, blue]
table {%
500 0.016
1000 0.038
2000 0.031
5000 0.077
10000 0.151
50000 0.847
100000 1.816
500000 9.472
};
\addplot [line width=1.25, red]
table {%
500 0.005
1000 0.011
2000 0.017
5000 0.042
10000 0.098
50000 0.487
100000 1.036
500000 5.342
};
\addplot [line width=1.25, magenta]
table {%
500 0.037
1000 0.017
2000 0.029
5000 0.073
10000 0.154
50000 0.827
100000 1.718
500000 8.872
};
\end{axis}

\end{tikzpicture}

%% file: trafalgar_iter_double.tex
%
%
%
\begin{tikzpicture}

\pgfplotsset{compat=newest} 

\definecolor{color0}{rgb}{1, 1, 1}

\tikzstyle{every node}=[font=\footnotesize]

\definecolor{colCholesky}{RGB}{0, 255, 0}
\definecolor{colQRkit}{RGB}{255, 0, 0}
\definecolor{colQRChol}{RGB}{0, 0, 255}
\definecolor{colSPQR}{RGB}{0, 255, 255}
\definecolor{colSSBA}{RGB}{255, 255, 0}
\definecolor{colMoreQR}{RGB}{200, 0, 255}

\begin{axis}[
xlabel={Iterations},
ylabel={Energy},
xmin=0, xmax=60,
ymin=1390, ymax=1900,
width=\figurewidth,
height=\figureheight,
at={(0\figurewidth,0\figureheight)},
xmajorgrids,
x grid style={lightgray},
ymajorgrids,
y grid style={lightgray},
axis line style={black},
axis background/.style={fill=color0},
legend style={at={(0.97,0.97)}, anchor=north east},
legend cell align={left},
legend entries={{Cholesky},{QRkit},{QRkit + Cholesky},{SuiteSparse QR},{SSBA},{Mor{\'e} QRkit}}
]
%
%
\addplot [line width=1, colCholesky]
table {%
1	1884.92000000000
2	1738.61000000000
3	1665.87000000000
4	1627.62000000000
5	1603.32000000000
6	1585.67000000000
7	1571.80000000000
8	1559.44000000000
9	1551.20000000000
10	1545.45000000000
11	1540.47000000000
12	1534.50000000000
13	1530.28000000000
14	1527.23000000000
15	1524.54000000000
16	1522.82000000000
17	1520.96000000000
18	1518.89000000000
19	1516.48000000000
20	1514.39000000000
21	1510.19000000000
22	1506.04000000000
23	1502.78000000000
24	1500.12000000000
25	1497.79000000000
26	1494.94000000000
27	1492.10000000000
28	1490.14000000000
29	1488.06000000000
30	1486.09000000000
31	1484.24000000000
32	1482.90000000000
33	1481.85000000000
34	1480.47000000000
35	1479.11000000000
36	1477.58000000000
37	1476.15000000000
38	1473.48000000000
39	1471.58000000000
40	1470.04000000000
41	1468.74000000000
42	1467.86000000000
43	1466.84000000000
44	1466.13000000000
45	1465.32000000000
46	1465.30000000000
47	1462.66000000000
48	1461.91000000000
49	1461.39000000000
50	1460.85000000000
51	1460.21000000000
52	1459.80000000000
53	1459.48000000000
54	1459.23000000000
55	1458.90000000000
56	1458.64000000000
57	1458.31000000000
58	1458.09000000000
59	1457.90000000000
60	1457.61000000000
61	1457.50000000000
62	1457.27000000000
63	1457.19000000000
64	1457.04000000000
65	1456.93000000000
66	1456.81000000000
67	1456.68000000000
68	1456.54000000000
69	1456.48000000000
70	1456.44000000000
71	1456.24000000000
72	1456.02000000000
73	1455.93000000000
74	1455.89000000000
75	1455.87000000000
76	1455.86000000000
77	1455.84000000000
78	1455.77000000000
79	1455.76000000000
80	1455.75000000000
81	1455.72000000000
82	1455.58000000000
83	1455.46000000000
84	1455.30000000000
85	1455.20000000000
86	1455.13000000000
87	1454.97000000000
88	1454.88000000000
89	1454.84000000000
90	1454.75000000000
91	1454.67000000000
92	1454.55000000000
93	1454.47000000000
94	1454.43000000000
95	1454.33000000000
96	1454.25000000000
97	1454.23000000000
98	1454.22000000000
99	1454.20000000000
100	1454.20000000000
101	1454.19000000000
102	1454.16000000000
103	1453.99000000000
104	1453.90000000000
105	1453.84000000000
106	1453.79000000000
107	1453.75000000000
108	1453.65000000000
109	1453.61000000000
110	1453.52000000000
111	1453.48000000000
112	1453.42000000000
113	1453.33000000000
114	1453.25000000000
115	1453.14000000000
116	1453.07000000000
117	1452.97000000000
118	1452.66000000000
119	1452.52000000000
120	1452.45000000000
121	1452.42000000000
122	1452.40000000000
123	1452.38000000000
124	1452.37000000000
125	1452.36000000000
126	1452.35000000000
127	1452.35000000000
128	1452.33000000000
129	1452.31000000000
130	1452.29000000000
131	1452.26000000000
132	1452.21000000000
133	1452.19000000000
134	1452.15000000000
135	1452.11000000000
136	1452.06000000000
137	1451.97000000000
138	1451.85000000000
139	1451.68000000000
140	1451.61000000000
141	1451.55000000000
142	1451.52000000000
143	1451.50000000000
144	1451.47000000000
145	1451.41000000000
146	1451.30000000000
147	1451.22000000000
148	1451.18000000000
149	1451.17000000000
150	1451.10000000000
151	1451.03000000000
152	1450.97000000000
153	1450.96000000000
154	1450.95000000000
155	1450.94000000000
156	1450.93000000000
157	1450.90000000000
158	1450.86000000000
159	1450.83000000000
160	1450.80000000000
161	1450.78000000000
162	1450.76000000000
163	1450.74000000000
164	1450.74000000000
165	1450.73000000000
166	1450.73000000000
167	1450.73000000000
168	1450.72000000000
169	1450.70000000000
170	1450.66000000000
171	1450.62000000000
172	1450.60000000000
173	1450.59000000000
174	1450.56000000000
175	1450.52000000000
176	1450.50000000000
177	1450.48000000000
178	1450.47000000000
179	1450.47000000000
180	1450.47000000000
181	1450.47000000000
182	1450.47000000000
183	1450.46000000000
184	1450.46000000000
185	1450.46000000000
186	1450.46000000000
187	1450.46000000000
188	1450.45000000000
189	1450.45000000000
190	1450.45000000000
191	1450.45000000000
192	1450.45000000000
193	1450.43000000000
194	1450.40000000000
195	1450.35000000000
196	1450.35000000000
197	1450.35000000000
198	1450.35000000000
199	1450.35000000000
200	1450.35000000000
201	1450.35000000000
202	1450.35000000000
203	1450.35000000000
204	1450.35000000000
205	1450.35000000000
206	1450.34000000000
207	1450.34000000000
208	1450.34000000000
209	1450.34000000000
210	1450.34000000000
211	1450.34000000000
212	1450.34000000000
213	1450.34000000000
214	1450.34000000000
215	1450.34000000000
216	1450.34000000000
217	1450.34000000000
218	1450.34000000000
219	1450.34000000000
220	1450.34000000000
221	1450.34000000000
};
\addplot [line width=1, colQRkit]
table {%
 1	1884.92000000000
2	1738.61000000000
3	1665.87000000000
4	1627.62000000000
5	1603.32000000000
6	1585.67000000000
7	1571.81000000000
8	1559.45000000000
9	1551.23000000000
10	1545.42000000000
11	1539.37000000000
12	1534.60000000000
13	1530.52000000000
14	1527.22000000000
15	1524.45000000000
16	1522.13000000000
17	1519.87000000000
18	1517.88000000000
19	1515.97000000000
20	1514.41000000000
21	1512.93000000000
22	1511.46000000000
23	1509.49000000000
24	1507.99000000000
25	1506.47000000000
26	1505.20000000000
27	1503.62000000000
28	1502.39000000000
29	1500.97000000000
30	1499.32000000000
31	1497.87000000000
32	1496.24000000000
33	1492.68000000000
34	1486.90000000000
35	1484.56000000000
36	1482.84000000000
37	1481.31000000000
38	1480.11000000000
39	1476.27000000000
40	1474.35000000000
41	1472.88000000000
42	1471.61000000000
43	1470.54000000000
44	1469.62000000000
45	1468.97000000000
46	1468.37000000000
47	1467.58000000000
48	1467.01000000000
49	1466.55000000000
50	1466.12000000000
51	1465.73000000000
52	1465.32000000000
53	1465.05000000000
54	1464.93000000000
55	1464.73000000000
56	1464.45000000000
57	1464.31000000000
58	1464.18000000000
59	1463.85000000000
60	1463.64000000000
61	1463.39000000000
62	1463.02000000000
63	1462.93000000000
64	1462.79000000000
65	1462.63000000000
66	1462.52000000000
67	1462.39000000000
68	1462.28000000000
69	1462.22000000000
70	1462.17000000000
71	1462.04000000000
72	1461.90000000000
73	1461.74000000000
74	1461.62000000000
75	1461.59000000000
76	1461.57000000000
77	1461.52000000000
78	1461.48000000000
79	1461.42000000000
80	1461.37000000000
81	1461.34000000000
82	1461.26000000000
83	1461.17000000000
84	1461.09000000000
85	1460.98000000000
86	1460.86000000000
87	1460.80000000000
88	1460.73000000000
89	1460.67000000000
90	1460.59000000000
91	1460.57000000000
92	1460.45000000000
93	1460.43000000000
};
\addplot [line width=1, colQRChol]
table {%
1	1884.92000000000
2	1738.61000000000
3	1665.87000000000
4	1627.62000000000
5	1603.32000000000
6	1585.67000000000
7	1571.81000000000
8	1559.45000000000
9	1551.23000000000
10	1545.42000000000
11	1539.42000000000
12	1534.10000000000
13	1529.77000000000
14	1526.27000000000
15	1523.46000000000
16	1520.80000000000
17	1518.18000000000
18	1515.95000000000
19	1513.72000000000
20	1511.97000000000
21	1507.59000000000
22	1502.39000000000
23	1497.32000000000
24	1493.54000000000
25	1490.90000000000
26	1488.47000000000
27	1486.08000000000
28	1483.62000000000
29	1481.63000000000
30	1479.88000000000
31	1478.28000000000
32	1477.01000000000
33	1475.92000000000
34	1475.03000000000
35	1472.95000000000
36	1471.22000000000
37	1469.07000000000
38	1467.20000000000
39	1467.02000000000
40	1460.11000000000
41	1458.10000000000
42	1456.86000000000
43	1455.38000000000
44	1453.84000000000
45	1452.84000000000
46	1452.53000000000
47	1451.96000000000
48	1451.41000000000
49	1450.92000000000
50	1450.62000000000
51	1450.26000000000
52	1449.75000000000
53	1449.55000000000
54	1449.45000000000
55	1449.13000000000
56	1448.84000000000
57	1448.68000000000
58	1448.56000000000
59	1448.43000000000
60	1448.27000000000
61	1448.04000000000
62	1447.86000000000
63	1447.58000000000
64	1447.12000000000
65	1446.93000000000
66	1446.66000000000
67	1446.30000000000
68	1445.99000000000
69	1445.89000000000
70	1445.73000000000
71	1445.54000000000
72	1445.42000000000
73	1445.38000000000
74	1445.33000000000
75	1445.28000000000
76	1445.14000000000
77	1444.97000000000
78	1444.82000000000
79	1444.75000000000
80	1444.68000000000
81	1444.67000000000
82	1444.63000000000
83	1444.42000000000
84	1444.23000000000
85	1444.13000000000
86	1444.11000000000
87	1444.11000000000
88	1444.10000000000
89	1444.10000000000
};
\addplot [line width=1, colSPQR]
table {%
1	1884.92000000000
2	1738.61000000000
3	1665.87000000000
4	1627.62000000000
5	1603.32000000000
6	1585.67000000000
7	1571.81000000000
8	1559.45000000000
9	1551.23000000000
10	1545.42000000000
11	1539.37000000000
12	1534.59000000000
13	1530.44000000000
14	1527.04000000000
15	1524.28000000000
16	1521.94000000000
17	1519.82000000000
18	1517.54000000000
19	1513.53000000000
20	1511.22000000000
21	1504.02000000000
22	1498.49000000000
23	1495.96000000000
24	1490.87000000000
25	1489.71000000000
26	1485.78000000000
27	1482.95000000000
28	1479.71000000000
29	1478.04000000000
30	1476.60000000000
31	1475.31000000000
32	1474.13000000000
33	1472.40000000000
34	1471.03000000000
35	1470.03000000000
36	1468.96000000000
37	1467.62000000000
38	1467.43000000000
39	1464.61000000000
40	1463.26000000000
41	1462.31000000000
42	1461.39000000000
43	1460.01000000000
44	1459.40000000000
45	1458.34000000000
46	1457.57000000000
47	1456.37000000000
48	1455.79000000000
49	1455.25000000000
50	1454.85000000000
51	1454.69000000000
52	1454.25000000000
53	1453.80000000000
54	1453.61000000000
55	1453.49000000000
56	1453.48000000000
57	1453.37000000000
58	1453.22000000000
59	1453.01000000000
60	1452.78000000000
61	1452.37000000000
62	1452
63	1451.88000000000
64	1451.71000000000
65	1451.49000000000
66	1451.36000000000
67	1451.25000000000
68	1451.13000000000
69	1451.04000000000
70	1450.89000000000
71	1450.81000000000
72	1450.62000000000
73	1450.46000000000
74	1450.26000000000
75	1450.25000000000
76	1449.73000000000
77	1449.24000000000
78	1448.83000000000
79	1448.03000000000
80	1447.60000000000
81	1446.98000000000
82	1446.46000000000
83	1445.83000000000
84	1445.22000000000
85	1445
86	1444.25000000000
87	1444.13000000000
88	1444.05000000000
89	1443.97000000000
90	1443.94000000000
91	1443.93000000000
92	1443.92000000000
93	1443.88000000000
94	1443.80000000000
95	1443.76000000000
96	1443.73000000000
97	1443.66000000000
98	1443.55000000000
99	1443.31000000000
100	1443.22000000000
101	1443.07000000000
102	1442.54000000000
103	1442.24000000000
104	1442.11000000000
105	1441.99000000000
106	1441.81000000000
107	1441.72000000000
108	1441.69000000000
109	1441.57000000000
110	1441.36000000000
111	1441.21000000000
112	1441.20000000000
113	1441.19000000000
114	1441.15000000000
115	1441.13000000000
116	1441.12000000000
117	1441.07000000000
118	1441.03000000000
119	1441
120	1440.99000000000
121	1440.97000000000
122	1440.87000000000
123	1440.77000000000
124	1440.70000000000
125	1440.60000000000
126	1440.30000000000
127	1440.16000000000
128	1440.16000000000
129	1439.93000000000
130	1439.75000000000
131	1439.44000000000
132	1439.19000000000
133	1439.02000000000
134	1438.95000000000
135	1438.92000000000
136	1438.91000000000
137	1438.91000000000
138	1438.90000000000
139	1438.90000000000
140	1438.85000000000
141	1438.81000000000
142	1438.78000000000
143	1438.77000000000
144	1438.71000000000
145	1438.70000000000
146	1438.70000000000
147	1438.69000000000
148	1438.63000000000
149	1438.62000000000
150	1438.62000000000
151	1438.57000000000
152	1438.48000000000
153	1438.36000000000
154	1438.21000000000
155	1438.18000000000
156	1438.03000000000
157	1437.70000000000
158	1437.24000000000
159	1437.21000000000
160	1437.14000000000
161	1437.11000000000
162	1437.11000000000
163	1437.10000000000
164	1437.10000000000
165	1437.10000000000
166	1437.10000000000
167	1437.06000000000
168	1436.95000000000
169	1436.86000000000
170	1436.82000000000
171	1436.71000000000
172	1436.66000000000
173	1436.65000000000
174	1436.65000000000
175	1436.65000000000
176	1436.65000000000
177	1436.65000000000
178	1436.65000000000
179	1436.65000000000
180	1436.64000000000
181	1436.64000000000
182	1436.59000000000
183	1436.58000000000
184	1436.58000000000
185	1436.58000000000
186	1436.58000000000
187	1436.58000000000
188	1436.58000000000
189	1436.58000000000
190	1436.58000000000
191	1436.58000000000
192	1436.58000000000
193	1436.58000000000
194	1436.58000000000
195	1436.58000000000
196	1436.58000000000
};
\addplot [line width=1, colSSBA]
table {%
1	1884.92280600000
2	1734.07790600000
3	1660.54016500000
4	1624.89651800000
5	1605.50586900000
6	1592.26627300000
7	1583.00016200000
8	1574.54291300000
9	1561.15270200000
10	1553.88673500000
11	1548.34720400000
12	1543.55175900000
13	1537.94512900000
14	1534.52602300000
15	1530.61781200000
16	1522.14844200000
17	1516.86848400000
18	1513.28153600000
19	1510.86623200000
20	1509.14772900000
21	1507.18934700000
22	1504.18271500000
23	1501.59653700000
24	1500.14269900000
25	1498.75389400000
26	1497.45283700000
27	1496.14583100000
28	1494.82456300000
29	1493.33720300000
30	1491.58757400000
31	1489.82975100000
32	1488.71693000000
33	1487.66106900000
34	1486.25850000000
35	1484.87744200000
36	1483.41597800000
37	1482.51171900000
38	1481.39443700000
39	1479.06908900000
40	1477.22515000000
41	1475.62140300000
42	1474.46849800000
43	1473.02806400000
44	1472.04402600000
45	1471.26495200000
46	1470.73892200000
47	1469.73563500000
48	1468.37971300000
49	1466.86907300000
50	1466.03879000000
51	1465.41386100000
52	1464.04137400000
53	1462.17842600000
54	1461.73500900000
55	1460.15103400000
56	1459.57352500000
57	1459.18837800000
58	1458.89429800000
59	1458.63022800000
60	1458.27764600000
61	1458.01781800000
62	1457.85866900000
63	1457.51823300000
64	1457.36538000000
65	1457.24448300000
66	1456.84293800000
67	1456.65203000000
68	1456.53351600000
69	1456.35962600000
70	1456.13490200000
71	1456.07037800000
72	1456.00600000000
73	1455.92433700000
74	1455.73046800000
75	1455.48771700000
76	1455.39367400000
77	1455.14336700000
78	1455.02175500000
79	1454.57533100000
80	1454.37314900000
81	1454.18904300000
};
\addplot [line width=1, colMoreQR]
table {%
1	1884.92000000000
2	1744.37000000000
3	1678.62000000000
4	1647.14000000000
5	1626.45000000000
6	1609.59000000000
7	1595.04000000000
8	1583.20000000000
9	1574.21000000000
10	1566.25000000000
11	1556.36000000000
12	1545.91000000000
13	1538.65000000000
14	1534.94000000000
15	1526.77000000000
16	1525.72000000000
17	1514.49000000000
18	1501
19	1494.46000000000
20	1488.56000000000
21	1484.87000000000
22	1483.75000000000
23	1482.75000000000
24	1480.20000000000
25	1473.88000000000
26	1468.29000000000
27	1457.77000000000
28	1457.40000000000
29	1456.27000000000
30	1453.44000000000
31	1451.92000000000
32	1448.04000000000
33	1447.78000000000
34	1447.57000000000
35	1446.26000000000
36	1444.39000000000
37	1443.01000000000
38	1442.25000000000
39	1441.37000000000
40	1440
41	1439.52000000000
42	1439.50000000000
43	1438.70000000000
44	1438.16000000000
45	1437.79000000000
46	1436.71000000000
47	1436.21000000000
48	1434.41000000000
49	1434.04000000000
50	1433.87000000000
51	1433.21000000000
52	1432.12000000000
53	1431.61000000000
54	1431.16000000000
55	1430.74000000000
56	1430.09000000000
57	1429.75000000000
58	1427.61000000000
59	1427.12000000000
60	1426.92000000000
61	1425.83000000000
62	1424.78000000000
63	1422.97000000000
64	1421.50000000000
65	1419.51000000000
66	1418.79000000000
67	1417.75000000000
68	1417.18000000000
69	1416.70000000000
70	1416.39000000000
71	1416.07000000000
72	1415.41000000000
73	1413.98000000000
74	1412.14000000000
75	1410.88000000000
76	1410.33000000000
77	1409.85000000000
78	1409.49000000000
79	1408.69000000000
80	1408.32000000000
81	1408.11000000000
82	1407.91000000000
83	1407.48000000000
84	1407.29000000000
85	1406.92000000000
86	1406.69000000000
87	1406.64000000000
88	1406.28000000000
89	1406.13000000000
90	1405.96000000000
91	1405.87000000000
92	1405.84000000000
93	1405.82000000000
94	1405.56000000000
95	1405.53000000000
96	1405.50000000000
97	1405.39000000000
98	1405.34000000000
99	1405.29000000000
100	1405.24000000000
101	1404.99000000000
102	1404.88000000000
103	1404.67000000000
104	1404.47000000000
105	1404.37000000000
106	1404.29000000000
107	1404.25000000000
108	1404.25000000000
109	1404.22000000000
110	1404.21000000000
111	1404.16000000000
112	1404.09000000000
113	1404
114	1403.94000000000
115	1403.68000000000
116	1403.68000000000
117	1403.45000000000
118	1403.22000000000
119	1403.14000000000
120	1402.93000000000
121	1402.88000000000
122	1402.80000000000
123	1402.64000000000
124	1402.56000000000
125	1402.53000000000
126	1402.47000000000
127	1402.17000000000
128	1402.08000000000
129	1402.08000000000
130	1402.03000000000
131	1402.01000000000
132	1401.95000000000
133	1401.83000000000
134	1401.83000000000
135	1401.82000000000
136	1401.47000000000
137	1401.02000000000
138	1400.91000000000
139	1400.80000000000
140	1400.80000000000
141	1400.73000000000
142	1400.64000000000
143	1400.58000000000
144	1400.49000000000
145	1400.40000000000
146	1400.32000000000
147	1399.76000000000
148	1399.51000000000
149	1399.27000000000
150	1398.95000000000
151	1398.72000000000
152	1398.47000000000
153	1398.26000000000
154	1398.14000000000
155	1397.95000000000
156	1397.90000000000
157	1397.87000000000
158	1397.66000000000
159	1397.50000000000
160	1397.43000000000
161	1397.28000000000
162	1397.01000000000
163	1396.88000000000
164	1396.82000000000
165	1396.80000000000
166	1396.58000000000
167	1396.55000000000
168	1396.43000000000
169	1396.29000000000
170	1396.29000000000
171	1396.22000000000
172	1396.09000000000
173	1396.06000000000
174	1396.04000000000
175	1396.03000000000
176	1396.01000000000
177	1396.01000000000
178	1396
179	1395.99000000000
180	1395.99000000000
181	1395.99000000000
182	1395.99000000000
183	1395.99000000000
184	1395.99000000000
185	1395.99000000000
186	1395.99000000000
187	1395.99000000000
188	1395.99000000000
189	1395.99000000000
190	1395.99000000000
191	1395.99000000000
192	1395.99000000000
193	1395.99000000000
194	1395.99000000000
};
\end{axis}

\end{tikzpicture}

%% file: trafalgar_double.tex
%
%
%
\begin{tikzpicture}

\pgfplotsset{compat=newest} 

\definecolor{color0}{rgb}{1, 1, 1}

\tikzstyle{every node}=[font=\footnotesize]

\definecolor{colCholesky}{RGB}{0, 255, 0}
\definecolor{colQRkit}{RGB}{255, 0, 0}
\definecolor{colQRChol}{RGB}{0, 0, 255}
\definecolor{colSPQR}{RGB}{0, 255, 255}
\definecolor{colSSBA}{RGB}{255, 255, 0}
\definecolor{colMoreQR}{RGB}{200, 0, 255}

\begin{axis}[
xmode=log,
log ticks with fixed point,
xtick={1,5,10,20,50,120,300,1000},
xlabel={Time [s]},
ylabel={Energy},
xmin=0, xmax=5e3,
ymin=1390, ymax=1900,
width=\figurewidth,
height=\figureheight,
at={(0\figurewidth,0\figureheight)},
xmajorgrids,
x grid style={lightgray},
ymajorgrids,
y grid style={lightgray},
axis line style={black},
axis background/.style={fill=color0},
]
%
%
\addplot [line width=0.75, colCholesky]
table {%
0.438000000000000 1884.92000000000
0.885000000000000 1738.61000000000
1.32200000000000  1665.87000000000
1.75600000000000  1627.62000000000
2.19100000000000  1603.32000000000
2.63200000000000  1585.67000000000
3.07400000000000  1571.80000000000
3.51300000000000  1559.44000000000
3.94300000000000  1551.20000000000
4.38300000000000  1545.45000000000
4.82400000000000  1540.47000000000
5.26400000000000  1534.50000000000
5.70100000000000  1530.28000000000
6.13200000000000  1527.23000000000
6.57200000000000  1524.54000000000
7.01500000000000  1522.82000000000
7.45700000000000  1520.96000000000
7.88500000000000  1518.89000000000
8.32300000000000  1516.48000000000
8.77500000000000  1514.39000000000
9.21600000000000  1510.19000000000
9.65600000000000  1506.04000000000
10.0930000000000  1502.78000000000
10.5320000000000  1500.12000000000
10.9710000000000  1497.79000000000
11.4170000000000  1494.94000000000
11.8460000000000  1492.10000000000
12.2920000000000  1490.14000000000
13.4040000000000  1488.06000000000
13.8450000000000  1486.09000000000
14.2750000000000  1484.24000000000
14.7130000000000  1482.90000000000
15.1560000000000  1481.85000000000
15.5920000000000  1480.47000000000
16.0270000000000  1479.11000000000
19.1130000000000  1477.58000000000
19.5550000000000  1476.15000000000
20.6490000000000  1473.48000000000
21.0940000000000  1471.58000000000
21.5360000000000  1470.04000000000
22.6450000000000  1468.74000000000
23.0800000000000  1467.86000000000
23.5210000000000  1466.84000000000
23.9630000000000  1466.13000000000
24.4040000000000  1465.32000000000
25.5060000000000  1465.30000000000
25.9470000000000  1462.66000000000
26.3920000000000  1461.91000000000
26.8320000000000  1461.39000000000
27.2700000000000  1460.85000000000
27.7080000000000  1460.21000000000
28.1490000000000  1459.80000000000
28.5860000000000  1459.48000000000
29.0310000000000  1459.23000000000
29.4580000000000  1458.90000000000
29.8970000000000  1458.64000000000
30.3420000000000  1458.31000000000
30.7850000000000  1458.09000000000
31.2210000000000  1457.90000000000
31.6540000000000  1457.61000000000
32.0960000000000  1457.50000000000
32.5370000000000  1457.27000000000
32.9790000000000  1457.19000000000
33.4090000000000  1457.04000000000
33.8490000000000  1456.93000000000
34.2920000000000  1456.81000000000
34.7330000000000  1456.68000000000
35.1660000000000  1456.54000000000
35.6060000000000  1456.48000000000
36.0480000000000  1456.44000000000
36.4880000000000  1456.24000000000
36.9310000000000  1456.02000000000
37.3640000000000  1455.93000000000
37.8030000000000  1455.89000000000
38.2450000000000  1455.87000000000
38.6860000000000  1455.86000000000
39.1190000000000  1455.84000000000
39.5530000000000  1455.77000000000
39.9920000000000  1455.76000000000
40.4340000000000  1455.75000000000
40.8750000000000  1455.72000000000
41.3060000000000  1455.58000000000
41.7440000000000  1455.46000000000
42.1850000000000  1455.30000000000
42.6280000000000  1455.20000000000
43.0640000000000  1455.13000000000
43.5000000000000  1454.97000000000
43.9420000000000  1454.88000000000
44.3860000000000  1454.84000000000
44.8240000000000  1454.75000000000
45.2580000000000  1454.67000000000
45.7000000000000  1454.55000000000
46.1380000000000  1454.47000000000
46.5830000000000  1454.43000000000
47.0160000000000  1454.33000000000
47.4540000000000  1454.25000000000
47.8950000000000  1454.23000000000
48.3370000000000  1454.22000000000
48.7750000000000  1454.20000000000
49.2060000000000  1454.20000000000
49.6460000000000  1454.19000000000
50.0890000000000  1454.16000000000
50.5260000000000  1453.99000000000
50.9580000000000  1453.90000000000
51.3960000000000  1453.84000000000
51.8370000000000  1453.79000000000
52.2790000000000  1453.75000000000
52.7140000000000  1453.65000000000
53.1460000000000  1453.61000000000
53.5890000000000  1453.52000000000
54.0300000000000  1453.48000000000
54.4720000000000  1453.42000000000
54.9050000000000  1453.33000000000
55.3460000000000  1453.25000000000
55.7870000000000  1453.14000000000
56.2300000000000  1453.07000000000
56.6700000000000  1452.97000000000
57.1010000000000  1452.66000000000
57.5440000000000  1452.52000000000
57.9890000000000  1452.45000000000
58.4310000000000  1452.42000000000
58.8640000000000  1452.40000000000
59.2990000000000  1452.38000000000
59.7400000000000  1452.37000000000
60.1820000000000  1452.36000000000
60.6180000000000  1452.35000000000
61.0510000000000  1452.35000000000
61.4910000000000  1452.33000000000
61.9330000000000  1452.31000000000
62.3740000000000  1452.29000000000
62.8050000000000  1452.26000000000
63.2450000000000  1452.21000000000
63.6860000000000  1452.19000000000
64.1280000000000  1452.15000000000
64.5670000000000  1452.11000000000
65  1452.06000000000
65.4410000000000  1451.97000000000
65.8840000000000  1451.85000000000
66.3280000000000  1451.68000000000
66.7600000000000  1451.61000000000
67.2000000000000  1451.55000000000
67.6380000000000  1451.52000000000
68.0800000000000  1451.50000000000
68.5210000000000  1451.47000000000
68.9560000000000  1451.41000000000
69.3960000000000  1451.30000000000
69.8390000000000  1451.22000000000
70.2850000000000  1451.18000000000
70.7310000000000  1451.17000000000
71.1690000000000  1451.10000000000
71.6080000000000  1451.03000000000
72.0500000000000  1450.97000000000
72.4830000000000  1450.96000000000
72.9150000000000  1450.95000000000
73.3530000000000  1450.94000000000
73.7910000000000  1450.93000000000
74.2330000000000  1450.90000000000
74.6660000000000  1450.86000000000
75.1010000000000  1450.83000000000
75.5470000000000  1450.80000000000
75.9890000000000  1450.78000000000
76.4210000000000  1450.76000000000
76.8570000000000  1450.74000000000
77.3020000000000  1450.74000000000
77.7430000000000  1450.73000000000
78.1850000000000  1450.73000000000
78.6130000000000  1450.73000000000
79.0510000000000  1450.72000000000
79.4900000000000  1450.70000000000
79.9330000000000  1450.66000000000
80.3680000000000  1450.62000000000
80.8050000000000  1450.60000000000
81.2480000000000  1450.59000000000
82.3570000000000  1450.56000000000
82.7860000000000  1450.52000000000
83.2210000000000  1450.50000000000
83.6640000000000  1450.48000000000
84.1010000000000  1450.47000000000
84.5400000000000  1450.47000000000
84.9700000000000  1450.47000000000
85.4150000000000  1450.47000000000
85.8550000000000  1450.47000000000
86.9660000000000  1450.46000000000
87.3970000000000  1450.46000000000
87.8390000000000  1450.46000000000
88.2840000000000  1450.46000000000
88.7260000000000  1450.46000000000
89.1580000000000  1450.45000000000
89.5970000000000  1450.45000000000
90.0380000000000  1450.45000000000
90.4800000000000  1450.45000000000
90.9199999999999  1450.45000000000
91.3490000000000  1450.43000000000
103.440000000000  1450.40000000000
103.882000000000  1450.35000000000
104.324000000000  1450.35000000000
108.748000000000  1450.35000000000
109.189000000000  1450.35000000000
109.629000000000  1450.35000000000
110.056000000000  1450.35000000000
110.493000000000  1450.35000000000
110.933000000000  1450.35000000000
111.379000000000  1450.35000000000
111.816000000000  1450.35000000000
112.274000000000  1450.35000000000
112.711000000000  1450.34000000000
113.154000000000  1450.34000000000
113.597000000000  1450.34000000000
114.024000000000  1450.34000000000
114.464000000000  1450.34000000000
114.908000000000  1450.34000000000
115.346000000000  1450.34000000000
115.785000000000  1450.34000000000
116.221000000000  1450.34000000000
116.663000000000  1450.34000000000
117.102000000000  1450.34000000000
117.566000000000  1450.34000000000
117.995000000000  1450.34000000000
118.433000000000  1450.34000000000
118.873000000000  1450.34000000000
119.316000000000  1450.34000000000
};
\addplot [line width=0.75, colQRkit]
table {%
 2.25200000000000 1884.92000000000
4.70300000000000  1738.61000000000
7.16400000000000  1665.87000000000
9.62300000000000  1627.62000000000
12.0830000000000  1603.32000000000
14.5230000000000  1585.67000000000
16.9880000000000  1571.81000000000
19.4540000000000  1559.45000000000
21.9090000000000  1551.23000000000
24.3800000000000  1545.42000000000
26.8370000000000  1539.37000000000
29.3170000000000  1534.60000000000
31.7780000000000  1530.52000000000
34.2420000000000  1527.22000000000
36.7120000000000  1524.45000000000
51.5110000000000  1522.13000000000
53.9650000000000  1519.87000000000
56.4560000000000  1517.88000000000
58.9250000000000  1515.97000000000
61.4620000000000  1514.41000000000
63.9520000000000  1512.93000000000
66.4160000000000  1511.46000000000
68.8780000000000  1509.49000000000
71.3300000000000  1507.99000000000
73.7850000000000  1506.47000000000
76.2510000000000  1505.20000000000
78.7230000000000  1503.62000000000
81.1830000000000  1502.39000000000
83.6320000000000  1500.97000000000
86.0830000000000  1499.32000000000
88.5510000000000  1497.87000000000
91.0070000000000  1496.24000000000
93.4830000000000  1492.68000000000
95.9370000000000  1486.90000000000
98.3720000000000  1484.56000000000
100.832000000000  1482.84000000000
103.296000000000  1481.31000000000
105.752000000000  1480.11000000000
108.193000000000  1476.27000000000
110.652000000000  1474.35000000000
113.103000000000  1472.88000000000
115.571000000000  1471.61000000000
118.030000000000  1470.54000000000
120.492000000000  1469.62000000000
122.945000000000  1468.97000000000
125.388000000000  1468.37000000000
127.850000000000  1467.58000000000
130.312000000000  1467.01000000000
132.767000000000  1466.55000000000
135.220000000000  1466.12000000000
137.669000000000  1465.73000000000
140.110000000000  1465.32000000000
142.562000000000  1465.05000000000
145.016000000000  1464.93000000000
147.481000000000  1464.73000000000
149.943000000000  1464.45000000000
152.399000000000  1464.31000000000
154.862000000000  1464.18000000000
157.316000000000  1463.85000000000
159.884000000000  1463.64000000000
162.464000000000  1463.39000000000
165.092000000000  1463.02000000000
167.846000000000  1462.93000000000
170.323000000000  1462.79000000000
172.772000000000  1462.63000000000
175.228000000000  1462.52000000000
177.685000000000  1462.39000000000
180.148000000000  1462.28000000000
182.585000000000  1462.22000000000
185.036000000000  1462.17000000000
187.521000000000  1462.04000000000
190.122000000000  1461.90000000000
192.694000000000  1461.74000000000
195.141000000000  1461.62000000000
197.612000000000  1461.59000000000
200.067000000000  1461.57000000000
202.517000000000  1461.52000000000
204.984000000000  1461.48000000000
207.431000000000  1461.42000000000
209.956000000000  1461.37000000000
212.420000000000  1461.34000000000
214.866000000000  1461.26000000000
217.314000000000  1461.17000000000
219.768000000000  1461.09000000000
222.247000000000  1460.98000000000
224.698000000000  1460.86000000000
227.147000000000  1460.80000000000
229.601000000000  1460.73000000000
232.061000000000  1460.67000000000
234.518000000000  1460.59000000000
236.968000000000  1460.57000000000
251.738000000000  1460.45000000000
273.857000000000  1460.43000000000
};
\addplot [line width=0.75, colQRChol]
table {%
1.29600000000000  1884.92000000000
2.45900000000000  1738.61000000000
3.62200000000000  1665.87000000000
4.81100000000000  1627.62000000000
5.97900000000000  1603.32000000000
7.16000000000000  1585.67000000000
8.33900000000000  1571.81000000000
9.51600000000000  1559.45000000000
10.6990000000000  1551.23000000000
11.8700000000000  1545.42000000000
13.0510000000000  1539.42000000000
14.2290000000000  1534.10000000000
15.4000000000000  1529.77000000000
16.5820000000000  1526.27000000000
17.7840000000000  1523.46000000000
18.9640000000000  1520.80000000000
20.1430000000000  1518.18000000000
21.3200000000000  1515.95000000000
22.5000000000000  1513.72000000000
23.6730000000000  1511.97000000000
24.8520000000000  1507.59000000000
26.0290000000000  1502.39000000000
29.4810000000000  1497.32000000000
32.9250000000000  1493.54000000000
34.1930000000000  1490.90000000000
35.3600000000000  1488.47000000000
36.5390000000000  1486.08000000000
37.7020000000000  1483.62000000000
41.1650000000000  1481.63000000000
42.3440000000000  1479.88000000000
43.5210000000000  1478.28000000000
44.7000000000000  1477.01000000000
45.8810000000000  1475.92000000000
47.0470000000000  1475.03000000000
48.2240000000000  1472.95000000000
51.6700000000000  1471.22000000000
52.8520000000000  1469.07000000000
54.0230000000000  1467.20000000000
55.1970000000000  1467.02000000000
56.3640000000000  1460.11000000000
57.5290000000000  1458.10000000000
60.9790000000000  1456.86000000000
62.1650000000000  1455.38000000000
65.6070000000000  1453.84000000000
66.7770000000000  1452.84000000000
67.9620000000000  1452.53000000000
69.1420000000000  1451.96000000000
70.3030000000000  1451.41000000000
71.4790000000000  1450.92000000000
72.6460000000000  1450.62000000000
73.8110000000000  1450.26000000000
74.9890000000000  1449.75000000000
76.1550000000000  1449.55000000000
77.3390000000000  1449.45000000000
78.5040000000000  1449.13000000000
79.6740000000000  1448.84000000000
80.8510000000000  1448.68000000000
82.1140000000000  1448.56000000000
83.3310000000000  1448.43000000000
84.5240000000000  1448.27000000000
85.6950000000000  1448.04000000000
86.8710000000000  1447.86000000000
90.3070000000000  1447.58000000000
91.4680000000001  1447.12000000000
92.6310000000000  1446.93000000000
93.8060000000000  1446.66000000000
94.9670000000000  1446.30000000000
96.1440000000001  1445.99000000000
97.3110000000001  1445.89000000000
98.4810000000001  1445.73000000000
99.6480000000001  1445.54000000000
103.087000000000  1445.42000000000
104.248000000000  1445.38000000000
105.426000000000  1445.33000000000
106.600000000000  1445.28000000000
107.766000000000  1445.14000000000
108.943000000000  1444.97000000000
110.104000000000  1444.82000000000
111.279000000000  1444.75000000000
112.448000000000  1444.68000000000
113.612000000000  1444.67000000000
117.054000000000  1444.63000000000
118.230000000000  1444.42000000000
119.397000000000  1444.23000000000
120.569000000000  1444.13000000000
121.744000000000  1444.11000000000
122.908000000000  1444.11000000000
124.084000000000  1444.10000000000
125.245000000000  1444.10000000000
};
\addplot [line width=0.75, colSPQR]
table {%
8.30000000000000  1884.92000000000
16.5550000000000  1738.61000000000
24.8050000000000  1665.87000000000
33.0450000000000  1627.62000000000
41.2840000000000  1603.32000000000
49.5140000000000  1585.67000000000
57.7480000000000  1571.81000000000
66.0100000000000  1559.45000000000
74.2710000000000  1551.23000000000
82.5300000000000  1545.42000000000
90.8790000000000  1539.37000000000
99.1140000000000  1534.59000000000
107.356000000000  1530.44000000000
115.595000000000  1527.04000000000
123.950000000000  1524.28000000000
132.322000000000  1521.94000000000
140.764000000000  1519.82000000000
149.171000000000  1517.54000000000
157.568000000000  1513.53000000000
174.328000000000  1511.22000000000
191.073000000000  1504.02000000000
207.870000000000  1498.49000000000
216.251000000000  1495.96000000000
224.627000000000  1490.87000000000
233.252000000000  1489.71000000000
241.806000000000  1485.78000000000
250.184000000000  1482.95000000000
283.948000000000  1479.71000000000
292.314000000000  1478.04000000000
300.679000000000  1476.60000000000
317.438000000000  1475.31000000000
350.904000000000  1474.13000000000
359.257000000000  1472.40000000000
367.638000000000  1471.03000000000
376.023000000000  1470.03000000000
392.777000000000  1468.96000000000
401.141000000000  1467.62000000000
409.497000000000  1467.43000000000
417.853000000000  1464.61000000000
426.225000000000  1463.26000000000
434.612000000000  1462.31000000000
442.995000000000  1461.39000000000
451.361000000000  1460.01000000000
468.099000000000  1459.40000000000
484.813000000000  1458.34000000000
511.115000000000  1457.57000000000
520.713000000000  1456.37000000000
530.371000000000  1455.79000000000
539.888000000000  1455.25000000000
548.282000000000  1454.85000000000
556.655000000000  1454.69000000000
565.404000000000  1454.25000000000
573.996000000000  1453.80000000000
582.557000000000  1453.61000000000
590.971000000000  1453.49000000000
599.340000000000  1453.48000000000
607.732000000000  1453.37000000000
616.134000000000  1453.22000000000
624.529000000000  1453.01000000000
632.910000000000  1452.78000000000
641.269000000000  1452.37000000000
674.731000000000  1452
683.114000000000  1451.88000000000
691.476000000000  1451.71000000000
699.847000000000  1451.49000000000
708.219000000000  1451.36000000000
716.593000000000  1451.25000000000
725.037000000000  1451.13000000000
733.425000000000  1451.04000000000
741.878000000000  1450.89000000000
750.247000000000  1450.81000000000
766.990000000000  1450.62000000000
775.351000000000  1450.46000000000
783.723000000000  1450.26000000000
792.085000000000  1450.25000000000
800.592000000000  1449.73000000000
808.946000000000  1449.24000000000
817.361000000000  1448.83000000000
825.727000000000  1448.03000000000
834.139000000000  1447.60000000000
842.504000000000  1446.98000000000
850.872000000000  1446.46000000000
884.357000000000  1445.83000000000
892.811000000000  1445.22000000000
901.163000000000  1445
926.202000000000  1444.25000000000
934.522000000000  1444.13000000000
942.844000000000  1444.05000000000
951.186000000000  1443.97000000000
959.521000000000  1443.94000000000
967.866000000000  1443.93000000000
976.234000000000  1443.92000000000
984.588000000000  1443.88000000000
993.029000000001  1443.80000000000
1001.37500000000  1443.76000000000
1009.73100000000  1443.73000000000
1043.35000000000  1443.66000000000
1051.70900000000  1443.55000000000
1060.07100000000  1443.31000000000
1085.15500000000  1443.22000000000
1093.50500000000  1443.07000000000
1101.87900000000  1442.54000000000
1110.22800000000  1442.24000000000
1118.59500000000  1442.11000000000
1126.95300000000  1441.99000000000
1135.30900000000  1441.81000000000
1143.68100000000  1441.72000000000
1152.05800000000  1441.69000000000
1160.50400000000  1441.57000000000
1168.86700000000  1441.36000000000
1177.22300000000  1441.21000000000
1185.58300000000  1441.20000000000
1193.92800000000  1441.19000000000
1210.63800000000  1441.15000000000
1227.35900000000  1441.13000000000
1244.08500000000  1441.12000000000
1260.78800000000  1441.07000000000
1269.15000000000  1441.03000000000
1277.52600000000  1441
1285.88400000000  1440.99000000000
1294.24700000000  1440.97000000000
1302.60100000000  1440.87000000000
1319.30100000000  1440.77000000000
1327.66000000000  1440.70000000000
1336.03600000000  1440.60000000000
1344.38100000000  1440.30000000000
1361.10800000000  1440.16000000000
1369.45600000000  1440.16000000000
1377.80900000000  1439.93000000000
1386.17300000000  1439.75000000000
1394.65300000000  1439.44000000000
1403.02800000000  1439.19000000000
1411.39400000000  1439.02000000000
1420.62700000000  1438.95000000000
1429.09800000000  1438.92000000000
1437.44300000000  1438.91000000000
1445.81100000000  1438.91000000000
1454.17900000000  1438.90000000000
1462.54500000000  1438.90000000000
1470.90300000000  1438.85000000000
1479.26900000000  1438.81000000000
1487.63300000000  1438.78000000000
1495.98000000000  1438.77000000000
1504.34600000000  1438.71000000000
1512.69600000000  1438.70000000000
1521.04700000000  1438.70000000000
1529.39900000000  1438.69000000000
1537.75300000000  1438.63000000000
1546.10100000000  1438.62000000000
1554.44200000000  1438.62000000000
1562.81100000000  1438.57000000000
1571.18400000000  1438.48000000000
1579.54000000000  1438.36000000000
1587.89100000000  1438.21000000000
1605.19200000000  1438.18000000000
1614.32300000000  1438.03000000000
1623.59600000000  1437.70000000000
1632.27700000000  1437.24000000000
1640.91500000000  1437.21000000000
1649.39500000000  1437.14000000000
1658.44500000000  1437.11000000000
1667.43300000000  1437.11000000000
1676.53900000000  1437.10000000000
1684.90700000000  1437.10000000000
1693.27300000000  1437.10000000000
1701.64500000000  1437.10000000000
1709.99300000000  1437.06000000000
1718.33500000000  1436.95000000000
1726.68700000000  1436.86000000000
1735.18800000000  1436.82000000000
1743.58000000000  1436.71000000000
1752.30400000000  1436.66000000000
1760.86200000000  1436.65000000000
1769.33100000000  1436.65000000000
1777.69200000000  1436.65000000000
1786.04700000000  1436.65000000000
1794.39100000000  1436.65000000000
1802.75000000000  1436.65000000000
1811.10500000000  1436.65000000000
1885.71800000000  1436.64000000000
1893.93800000000  1436.64000000000
1902.16300000000  1436.59000000000
1951.65800000000  1436.58000000000
1959.88100000000  1436.58000000000
1968.10200000000  1436.58000000000
1976.33300000000  1436.58000000000
1984.55600000000  1436.58000000000
1992.82500000000  1436.58000000000
2001.14900000000  1436.58000000000
2009.38100000000  1436.58000000000
2017.60700000000  1436.58000000000
2026.33900000000  1436.58000000000
2034.80700000000  1436.58000000000
2043.59800000000  1436.58000000000
2051.98000000000  1436.58000000000
2060.23100000000  1436.58000000000
};
\addplot [line width=0.75, colMoreQR]
table {%
2.97900000000000  1884.92000000000
5.78600000000000  1744.37000000000
8.58500000000000  1678.62000000000
11.4050000000000  1647.14000000000
14.3140000000000  1626.45000000000
17.1380000000000  1609.59000000000
19.9630000000000  1595.04000000000
22.8210000000000  1583.20000000000
25.6940000000000  1574.21000000000
28.5250000000000  1566.25000000000
31.3630000000000  1556.36000000000
34.1820000000000  1545.91000000000
37.0080000000000  1538.65000000000
43.3830000000000  1534.94000000000
48.0110000000000  1526.77000000000
50.8820000000000  1525.72000000000
53.9390000000000  1514.49000000000
58.5580000000000  1501
63.1850000000000  1494.46000000000
69.6030000000000  1488.56000000000
72.4490000000000  1484.87000000000
75.2710000000000  1483.75000000000
78.0960000000000  1482.75000000000
80.9300000000000  1480.20000000000
87.3280000000000  1473.88000000000
90.1630000000000  1468.29000000000
94.8050000000000  1457.77000000000
97.6340000000000  1457.40000000000
100.465000000000  1456.27000000000
103.287000000000  1453.44000000000
106.117000000000  1451.92000000000
116.132000000000  1448.04000000000
126.120000000000  1447.78000000000
132.545000000000  1447.57000000000
135.380000000000  1446.26000000000
138.206000000000  1444.39000000000
141.057000000000  1443.01000000000
143.922000000000  1442.25000000000
146.751000000000  1441.37000000000
151.384000000000  1440
154.226000000000  1439.52000000000
157.046000000000  1439.50000000000
159.881000000000  1438.70000000000
162.704000000000  1438.16000000000
165.537000000000  1437.79000000000
168.375000000000  1436.71000000000
171.300000000000  1436.21000000000
174.129000000000  1434.41000000000
180.536000000000  1434.04000000000
183.377000000000  1433.87000000000
186.223000000000  1433.21000000000
189.061000000000  1432.12000000000
191.912000000000  1431.61000000000
194.764000000000  1431.16000000000
197.644000000000  1430.74000000000
200.476000000000  1430.09000000000
203.316000000000  1429.75000000000
209.715000000000  1427.61000000000
212.544000000000  1427.12000000000
217.166000000000  1426.92000000000
220.015000000000  1425.83000000000
222.849000000000  1424.78000000000
231.030000000000  1422.97000000000
233.854000000000  1421.50000000000
238.503000000000  1419.51000000000
241.400000000000  1418.79000000000
246.023000000000  1417.75000000000
250.635000000000  1417.18000000000
257.031000000000  1416.70000000000
259.857000000000  1416.39000000000
262.679000000000  1416.07000000000
265.641000000000  1415.41000000000
270.238000000000  1413.98000000000
273.057000000000  1412.14000000000
281.241000000000  1410.88000000000
284.080000000000  1410.33000000000
288.701000000000  1409.85000000000
295.114000000000  1409.49000000000
297.952000000000  1408.69000000000
306.195000000000  1408.32000000000
309.044000000000  1408.11000000000
315.460000000000  1407.91000000000
320.116000000000  1407.48000000000
324.781000000000  1407.29000000000
327.639000000000  1406.92000000000
332.276000000000  1406.69000000000
335.126000000000  1406.64000000000
343.306000000000  1406.28000000000
346.127000000000  1406.13000000000
352.760000000000  1405.96000000000
355.602000000000  1405.87000000000
358.444000000000  1405.84000000000
361.287000000000  1405.82000000000
367.716000000000  1405.56000000000
372.348000000000  1405.53000000000
375.180000000000  1405.50000000000
383.579000000000  1405.39000000000
386.419000000000  1405.34000000000
389.267000000000  1405.29000000000
392.104000000000  1405.24000000000
400.294000000000  1404.99000000000
403.126000000000  1404.88000000000
405.950000000000  1404.67000000000
412.363000000000  1404.47000000000
415.207000000000  1404.37000000000
421.641000000000  1404.29000000000
426.296000000000  1404.25000000000
429.148000000000  1404.25000000000
437.398000000000  1404.22000000000
440.243000000000  1404.21000000000
443.082000000000  1404.16000000000
447.702000000000  1404.09000000000
450.534000000000  1404
453.371000000000  1403.94000000000
457.988000000000  1403.68000000000
460.821000000000  1403.68000000000
463.654000000000  1403.45000000000
468.308000000000  1403.22000000000
474.728000000000  1403.14000000000
481.157000000000  1402.93000000000
483.996000000000  1402.88000000000
490.484000000000  1402.80000000000
493.315000000000  1402.64000000000
499.738000000000  1402.56000000000
504.370000000000  1402.53000000000
507.211000000000  1402.47000000000
513.629000000000  1402.17000000000
516.449000000000  1402.08000000000
519.303000000001  1402.08000000000
522.135000000000  1402.03000000000
524.963000000000  1402.01000000000
527.796000000000  1401.95000000000
530.624000000000  1401.83000000000
535.260000000000  1401.83000000000
538.106000000000  1401.82000000000
540.944000000000  1401.47000000000
543.787000000000  1401.02000000000
548.413000000000  1400.91000000000
556.675000000000  1400.80000000000
559.555000000000  1400.80000000000
562.383000000000  1400.73000000000
565.219000000000  1400.64000000000
568.048000000000  1400.58000000000
572.662000000000  1400.49000000000
579.077000000000  1400.40000000000
581.898000000000  1400.32000000000
590.226000000000  1399.76000000000
593.076000000000  1399.51000000000
595.912000000000  1399.27000000000
599.050000000000  1398.95000000000
606.933000000000  1398.72000000000
611.942000000000  1398.47000000000
620.272000000000  1398.26000000000
623.167000000000  1398.14000000000
626.120000000000  1397.95000000000
634.333000000000  1397.90000000000
637.211000000000  1397.87000000000
640.121000000000  1397.66000000000
646.503000000000  1397.50000000000
649.334000000000  1397.43000000000
652.194000000000  1397.28000000000
655.028000000000  1397.01000000000
663.350000000000  1396.88000000000
666.212000000000  1396.82000000000
669.104000000000  1396.80000000000
675.592000000000  1396.58000000000
678.426000000000  1396.55000000000
684.876000000000  1396.43000000000
691.396000000000  1396.29000000000
694.234000000000  1396.29000000000
697.104000000000  1396.22000000000
703.582000000000  1396.09000000000
706.473000000000  1396.06000000000
709.359000000000  1396.04000000000
712.210000000000  1396.03000000000
716.834000000000  1396.01000000000
728.648000000000  1396.01000000000
742.196000000000  1396
745.016000000000  1395.99000000000
747.857000000000  1395.99000000000
754.234000000000  1395.99000000000
757.091000000000  1395.99000000000
761.708000000000  1395.99000000000
764.575000000000  1395.99000000000
767.397000000000  1395.99000000000
770.218000000000  1395.99000000000
773.047000000000  1395.99000000000
775.875000000000  1395.99000000000
778.707000000000  1395.99000000000
781.537000000000  1395.99000000000
784.435000000000  1395.99000000000
787.328000000000  1395.99000000000
790.151000000000  1395.99000000000
796.605000000000  1395.99000000000
};
\end{axis}

\end{tikzpicture}

%% file: trafalgar_iter_float.tex
%
%
%
\begin{tikzpicture}

\pgfplotsset{compat=newest} 

\definecolor{color0}{rgb}{1, 1, 1}

\tikzstyle{every node}=[font=\footnotesize]

\definecolor{colCholesky}{RGB}{0, 255, 0}
\definecolor{colQRkit}{RGB}{255, 0, 0}
\definecolor{colQRChol}{RGB}{0, 0, 255}
\definecolor{colSPQR}{RGB}{0, 255, 255}
\definecolor{colSSBA}{RGB}{255, 255, 0}
\definecolor{colMoreQR}{RGB}{200, 0, 255}

\begin{axis}[
xlabel={Iterations},
ylabel={Energy},
xmin=0, xmax=60,
ymin=1390, ymax=1900,
width=\figurewidth,
height=\figureheight,
at={(0\figurewidth,0\figureheight)},
xmajorgrids,
x grid style={lightgray},
ymajorgrids,
y grid style={lightgray},
axis line style={black},
axis background/.style={fill=color0},
legend style={at={(0.97,0.97)}, anchor=north east},
legend cell align={left},
legend entries={{Cholesky},{QRkit},{QRkit + Cholesky},{SSBA},{Mor{\'e} QRkit}}
]
%
%
\addplot [line width=1, colCholesky]
table {%
1	1884.92000000000
2	1737.58000000000
3	1664.07000000000
4	1628.82000000000
5	1609.81000000000
6	1595.41000000000
7	1584.77000000000
8	1577.28000000000
9	1570.96000000000
10	1567.23000000000
11	1563.06000000000
12	1559.98000000000
13	1556.65000000000
14	1553.30000000000
15	1550.96000000000
16	1549.53000000000
17	1549
18	1546.44000000000
19	1545.06000000000
20	1543.73000000000
21	1541.81000000000
22	1541.16000000000
23	1539.46000000000
24	1538.66000000000
25	1537.52000000000
26	1536.38000000000
27	1535.34000000000
28	1534.89000000000
29	1534.07000000000
30	1533.17000000000
31	1532.23000000000
32	1531.45000000000
33	1530.76000000000
34	1529.97000000000
35	1529.52000000000
36	1529.33000000000
37	1528.73000000000
38	1528.38000000000
39	1528.12000000000
40	1527.76000000000
41	1527.46000000000
42	1526.83000000000
43	1526.50000000000
44	1526.20000000000
45	1525.24000000000
46	1524.78000000000
47	1524.52000000000
48	1524.44000000000
49	1523.93000000000
50	1523.68000000000
51	1523.52000000000
52	1523.21000000000
53	1523.03000000000
54	1522.73000000000
55	1522.49000000000
56	1522.14000000000
57	1521.64000000000
58	1521.56000000000
59	1521.39000000000
60	1521.08000000000
61	1520.50000000000
62	1520.49000000000
63	1520.39000000000
64	1520.11000000000
65	1520
66	1519.81000000000
67	1519.53000000000
68	1519.26000000000
69	1519.11000000000
70	1519.03000000000
71	1518.87000000000
72	1518.68000000000
73	1518.60000000000
74	1518.58000000000
75	1518.39000000000
76	1518.33000000000
77	1518.24000000000
78	1518.19000000000
79	1518.19000000000
80	1518.15000000000
81	1518.12000000000
82	1518.05000000000
83	1518.02000000000
84	1517.83000000000
85	1517.79000000000
86	1517.74000000000
87	1517.73000000000
88	1517.71000000000
89	1517.70000000000
90	1517.70000000000
91	1517.68000000000
92	1517.68000000000
93	1517.68000000000
94	1517.68000000000
95	1517.67000000000
96	1517.66000000000
97	1517.60000000000
98	1517.58000000000
99	1517.58000000000
100	1517.58000000000
101	1517.58000000000
102	1517.58000000000
103	1517.58000000000
104	1517.58000000000
105	1517.58000000000
106	1517.58000000000
107	1517.57000000000
108	1517.57000000000
109	1517.57000000000
110	1517.57000000000
111	1517.57000000000
112	1517.57000000000
113	1517.57000000000
};
\addplot [line width=1, colQRkit]
table {%
 1	1884.92000000000
2	1737.81000000000
3	1666.85000000000
4	1629.68000000000
5	1606.70000000000
6	1593.65000000000
7	1576.90000000000
8	1566.94000000000
9	1559.61000000000
10	1556.85000000000
11	1547.83000000000
12	1539.91000000000
13	1538.29000000000
14	1530.91000000000
15	1527.57000000000
16	1525.52000000000
17	1522.41000000000
18	1520.38000000000
19	1518.37000000000
20	1516.53000000000
21	1515.37000000000
22	1513.81000000000
23	1512.22000000000
24	1508.87000000000
25	1506.96000000000
26	1505.26000000000
27	1504.13000000000
28	1502.78000000000
29	1501.27000000000
30	1500.38000000000
31	1497.44000000000
32	1496.24000000000
33	1495.73000000000
34	1494.28000000000
35	1492.95000000000
36	1491.62000000000
37	1490.33000000000
38	1489.90000000000
39	1489.12000000000
40	1486.73000000000
41	1484.61000000000
42	1483.74000000000
43	1482.12000000000
44	1480.22000000000
45	1478.64000000000
46	1477.93000000000
47	1477.20000000000
48	1476.73000000000
49	1476.33000000000
50	1475.95000000000
51	1474.84000000000
52	1474.06000000000
53	1473.02000000000
54	1471.56000000000
55	1471.04000000000
56	1470.49000000000
57	1470.44000000000
58	1470.25000000000
59	1470.19000000000
60	1470.13000000000
61	1470.05000000000
62	1469.99000000000
63	1469.86000000000
64	1469.79000000000
65	1469.71000000000
66	1469.55000000000
67	1469.52000000000
68	1469.34000000000
69	1469.09000000000
70	1468.90000000000
71	1468.67000000000
72	1468.53000000000
73	1468.43000000000
74	1468.37000000000
75	1468.33000000000
76	1468.27000000000
77	1468.18000000000
78	1468.11000000000
79	1468.09000000000
80	1468
81	1467.96000000000
82	1467.92000000000
83	1467.88000000000
84	1467.86000000000
85	1467.73000000000
86	1467.59000000000
87	1467.56000000000
88	1467.37000000000
89	1467.33000000000
90	1467.27000000000
91	1467.25000000000
92	1467.23000000000
93	1467.21000000000
94	1467.20000000000
95	1467.09000000000
96	1467.08000000000
97	1466.97000000000
98	1466.97000000000
99	1466.96000000000
100	1466.94000000000
101	1466.93000000000
102	1466.89000000000
103	1466.89000000000
104	1466.89000000000
105	1466.88000000000
106	1466.88000000000
107	1466.88000000000
108	1466.88000000000
109	1466.85000000000
110	1466.81000000000
111	1466.81000000000
112	1466.80000000000
113	1466.80000000000
114	1466.78000000000
115	1466.78000000000
116	1466.77000000000
117	1466.77000000000
118	1466.76000000000
119	1466.76000000000
120	1466.76000000000
121	1466.76000000000
122	1466.70000000000
123	1466.67000000000
124	1466.66000000000
125	1466.63000000000
126	1466.62000000000
127	1466.62000000000
128	1466.62000000000
129	1466.61000000000
130	1466.61000000000
131	1466.60000000000
132	1466.60000000000
133	1466.60000000000
134	1466.60000000000
135	1466.60000000000
136	1466.60000000000
137	1466.60000000000
138	1466.60000000000
139	1466.59000000000
140	1466.59000000000
141	1466.59000000000
142	1466.59000000000
};
\addplot [line width=1, colQRChol]
table {%
1	1884.92000000000
2	1736.18000000000
3	1663.65000000000
4	1627.78000000000
5	1605.12000000000
6	1588.68000000000
7	1575.62000000000
8	1570.17000000000
9	1555.11000000000
10	1545.70000000000
11	1540.51000000000
12	1534.18000000000
13	1531.18000000000
14	1525.57000000000
15	1524.62000000000
16	1523.10000000000
17	1522.18000000000
18	1521.49000000000
19	1520.51000000000
20	1519.27000000000
21	1518.26000000000
22	1517.34000000000
23	1517.17000000000
24	1516.23000000000
25	1515.85000000000
26	1515.36000000000
27	1514.94000000000
28	1514.84000000000
29	1514.19000000000
30	1513.88000000000
31	1513.29000000000
32	1513.07000000000
33	1512.32000000000
34	1511.92000000000
35	1511.78000000000
36	1511.08000000000
37	1510.74000000000
38	1510.42000000000
39	1510.17000000000
40	1509.77000000000
41	1509.63000000000
42	1508.89000000000
43	1508.63000000000
44	1508.59000000000
45	1508.11000000000
46	1507.57000000000
47	1507.16000000000
48	1506.85000000000
49	1506.52000000000
50	1506.02000000000
51	1505.65000000000
52	1505.20000000000
53	1504.93000000000
54	1504.25000000000
55	1504.13000000000
56	1504
57	1503.56000000000
58	1503.39000000000
59	1503.34000000000
60	1502.81000000000
61	1502.49000000000
62	1502.33000000000
63	1502.04000000000
64	1501.93000000000
65	1501.57000000000
66	1501.12000000000
67	1501.12000000000
68	1500.52000000000
69	1500.23000000000
70	1500
71	1499.91000000000
72	1499.73000000000
73	1499.70000000000
74	1499.57000000000
75	1499.40000000000
76	1499.23000000000
77	1499.04000000000
78	1498.92000000000
79	1498.86000000000
80	1498.82000000000
81	1498.70000000000
82	1498.60000000000
83	1498.56000000000
84	1498.50000000000
85	1498.47000000000
86	1498.33000000000
87	1498.27000000000
88	1498.24000000000
89	1498.21000000000
90	1498.20000000000
91	1498.18000000000
92	1498.16000000000
93	1498.07000000000
94	1497.98000000000
95	1497.74000000000
96	1497.69000000000
97	1497.67000000000
98	1497.66000000000
99	1497.61000000000
100	1497.55000000000
101	1497.41000000000
102	1497.28000000000
103	1497.13000000000
104	1496.93000000000
105	1496.86000000000
106	1496.82000000000
107	1496.74000000000
108	1496.66000000000
109	1496.65000000000
110	1496.58000000000
111	1496.55000000000
112	1496.55000000000
113	1496.54000000000
114	1496.52000000000
115	1496.52000000000
116	1496.50000000000
117	1496.48000000000
118	1496.48000000000
119	1496.47000000000
120	1496.47000000000
121	1496.46000000000
122	1496.44000000000
123	1496.35000000000
124	1496.35000000000
125	1496.34000000000
126	1496.33000000000
127	1496.31000000000
128	1496.30000000000
129	1496.26000000000
130	1496.22000000000
131	1496.19000000000
132	1496.17000000000
133	1496.16000000000
134	1496.14000000000
135	1496.11000000000
136	1496.10000000000
137	1496.05000000000
138	1496.04000000000
139	1496.04000000000
140	1495.97000000000
141	1495.90000000000
142	1495.78000000000
143	1495.78000000000
144	1495.69000000000
145	1495.68000000000
146	1495.64000000000
147	1495.57000000000
148	1495.54000000000
149	1495.51000000000
150	1495.47000000000
151	1495.39000000000
152	1495.38000000000
153	1495.38000000000
154	1495.35000000000
155	1495.33000000000
156	1495.28000000000
157	1495.26000000000
158	1495.18000000000
159	1495.17000000000
160	1495.14000000000
161	1495.09000000000
162	1495.07000000000
163	1495.06000000000
164	1495.04000000000
165	1495.03000000000
166	1495.02000000000
167	1495.02000000000
168	1495.01000000000
169	1495.01000000000
170	1495
171	1495
172	1495
173	1495
174	1494.94000000000
175	1494.93000000000
176	1494.93000000000
177	1494.93000000000
178	1494.93000000000
179	1494.93000000000
180	1494.93000000000
181	1494.93000000000
182	1494.93000000000
};
\addplot [line width=1, colSSBA]
table {%
1	1884.92163100000
2	1736.59716800000
3	1670.61999500000
4	1637.33715800000
5	1617.18603500000
6	1603.29138200000
7	1593.64697300000
8	1585.97204600000
9	1580.04113800000
10	1576.49939000000
11	1572.35510300000
12	1568.10339400000
13	1564.63647500000
14	1562.33728000000
15	1560.68237300000
16	1558.10937500000
17	1556.92468300000
18	1555.41467300000
19	1554.74694800000
20	1552.92565900000
21	1551.64709500000
22	1550.62463400000
23	1549.07714800000
24	1548.08581500000
25	1547.13464400000
26	1545.77941900000
27	1544.91613800000
28	1544.32128900000
29	1543.58154300000
30	1542.94531300000
31	1542.30029300000
32	1541.84558100000
33	1541.15502900000
34	1540.29589800000
35	1539.97290000000
36	1538.48584000000
37	1538.35974100000
38	1537.89709500000
39	1536.88305700000
40	1536.45532200000
41	1536.10253900000
42	1535.48498500000
43	1535.12280300000
44	1534.64941400000
45	1534.27258300000
46	1533.83142100000
47	1533.12927200000
48	1532.37255900000
49	1531.94934100000
50	1531.68884300000
};
\addplot [line width=1, colMoreQR]
table {%
1	1884.92000000000
2	1743.80000000000
3	1675.47000000000
4	1639.90000000000
5	1617.67000000000
6	1600.02000000000
7	1585.61000000000
8	1574.41000000000
9	1569.87000000000
10	1556.59000000000
11	1551.20000000000
12	1549.12000000000
13	1531.61000000000
14	1527.10000000000
15	1524.21000000000
16	1521.65000000000
17	1520.53000000000
18	1515.57000000000
19	1513.27000000000
20	1510.08000000000
21	1509.11000000000
22	1505.57000000000
23	1502.32000000000
24	1500.14000000000
25	1498.36000000000
26	1496.76000000000
27	1495.59000000000
28	1493.60000000000
29	1491.70000000000
30	1490.63000000000
31	1489.79000000000
32	1489.11000000000
33	1487.75000000000
34	1486.83000000000
35	1485.56000000000
36	1484.73000000000
37	1484.27000000000
38	1481.84000000000
39	1480.40000000000
40	1479.19000000000
41	1478.55000000000
42	1476.83000000000
43	1476.64000000000
44	1474.92000000000
45	1474.01000000000
46	1473.14000000000
47	1472.38000000000
48	1471.41000000000
49	1470.79000000000
50	1469.30000000000
51	1468.62000000000
52	1468.18000000000
53	1467.79000000000
54	1467.55000000000
55	1467.31000000000
56	1467.04000000000
57	1466.59000000000
58	1466.52000000000
59	1466.25000000000
60	1465.98000000000
61	1465.66000000000
62	1465.23000000000
63	1464.92000000000
64	1464.89000000000
65	1464.50000000000
66	1464.41000000000
67	1464.26000000000
68	1464.10000000000
69	1464.08000000000
70	1463.92000000000
71	1463.89000000000
72	1463.82000000000
73	1463.65000000000
74	1463.63000000000
75	1463.63000000000
76	1463.62000000000
77	1463.60000000000
78	1463.58000000000
79	1463.58000000000
80	1463.57000000000
81	1463.56000000000
82	1463.54000000000
83	1463.54000000000
84	1463.53000000000
85	1463.53000000000
86	1463.52000000000
87	1463.51000000000
88	1463.51000000000
89	1463.51000000000
90	1463.51000000000
91	1463.51000000000
92	1463.51000000000
93	1463.50000000000
94	1463.50000000000
95	1463.50000000000
96	1463.50000000000
97	1463.50000000000
98	1463.50000000000
99	1463.50000000000
100	1463.50000000000
101	1463.50000000000
102	1463.50000000000
};
\end{axis}

\end{tikzpicture}

%% file: trafalgar_float.tex
%
%
%
\begin{tikzpicture}

\pgfplotsset{compat=newest} 

\definecolor{color0}{rgb}{1, 1, 1}

\tikzstyle{every node}=[font=\footnotesize]

\definecolor{colCholesky}{RGB}{0, 255, 0}
\definecolor{colQRkit}{RGB}{255, 0, 0}
\definecolor{colQRChol}{RGB}{0, 0, 255}
\definecolor{colSPQR}{RGB}{0, 255, 255}
\definecolor{colSSBA}{RGB}{255, 255, 0}
\definecolor{colMoreQR}{RGB}{200, 0, 255}

\begin{axis}[
xmode=log,
log ticks with fixed point,
xtick={1,5,10,20,50,120,300,1000},
xlabel={Time [s]},
ylabel={Energy},
xmin=0, xmax=5e3,
ymin=1390, ymax=1900,
width=\figurewidth,
height=\figureheight,
at={(0\figurewidth,0\figureheight)},
xmajorgrids,
x grid style={lightgray},
ymajorgrids,
y grid style={lightgray},
axis line style={black},
axis background/.style={fill=color0},
]
%
%
\addplot [line width=1, colCholesky]
table {%
0.390000000000000	1884.92000000000
0.777000000000000	1737.58000000000
2.25600000000000	1664.07000000000
3.70600000000000	1628.82000000000
4.09400000000000	1609.81000000000
5.54100000000000	1595.41000000000
6.98500000000000	1584.77000000000
7.90600000000000	1577.28000000000
8.77200000000000	1570.96000000000
10.2180000000000	1567.23000000000
10.6070000000000	1563.06000000000
11.4680000000000	1559.98000000000
12.3190000000000	1556.65000000000
12.7080000000000	1553.30000000000
14.1440000000000	1550.96000000000
14.5380000000000	1549.53000000000
14.9260000000000	1549
15.3170000000000	1546.44000000000
15.7050000000000	1545.06000000000
16.0920000000000	1543.73000000000
16.4850000000000	1541.81000000000
17.3450000000000	1541.16000000000
17.7340000000000	1539.46000000000
18.1230000000000	1538.66000000000
18.5140000000000	1537.52000000000
18.9020000000000	1536.38000000000
19.2940000000000	1535.34000000000
19.6850000000000	1534.89000000000
20.0740000000000	1534.07000000000
20.4620000000000	1533.17000000000
20.8590000000000	1532.23000000000
21.2470000000000	1531.45000000000
22.7020000000000	1530.76000000000
23.0910000000000	1529.97000000000
23.4790000000000	1529.52000000000
23.8680000000000	1529.33000000000
24.7300000000000	1528.73000000000
25.1190000000000	1528.38000000000
25.5080000000000	1528.12000000000
25.9000000000000	1527.76000000000
26.2890000000000	1527.46000000000
26.6820000000000	1526.83000000000
27.0720000000000	1526.50000000000
27.4680000000000	1526.20000000000
28.4270000000000	1525.24000000000
28.8190000000000	1524.78000000000
29.2090000000000	1524.52000000000
29.5980000000000	1524.44000000000
29.9920000000000	1523.93000000000
30.3810000000000	1523.68000000000
30.7720000000000	1523.52000000000
31.1620000000000	1523.21000000000
31.5530000000000	1523.03000000000
32.5060000000000	1522.73000000000
32.8960000000000	1522.49000000000
33.2860000000000	1522.14000000000
33.6780000000000	1521.64000000000
34.6300000000000	1521.56000000000
35.0230000000000	1521.39000000000
35.4110000000000	1521.08000000000
36.3620000000000	1520.50000000000
36.7520000000000	1520.49000000000
37.1450000000000	1520.39000000000
38.0950000000000	1520.11000000000
38.4820000000000	1520
38.8690000000000	1519.81000000000
39.7320000000000	1519.53000000000
40.1220000000000	1519.26000000000
40.5090000000000	1519.11000000000
40.9020000000000	1519.03000000000
41.2930000000000	1518.87000000000
42.2460000000000	1518.68000000000
42.6330000000000	1518.60000000000
43.0220000000000	1518.58000000000
43.9690000000000	1518.39000000000
45.6550000000000	1518.33000000000
46.0440000000000	1518.24000000000
47.7240000000000	1518.19000000000
48.1110000000000	1518.19000000000
48.5020000000000	1518.15000000000
48.8910000000000	1518.12000000000
49.2800000000000	1518.05000000000
49.6710000000000	1518.02000000000
50.6270000000000	1517.83000000000
51.0180000000000	1517.79000000000
51.4090000000000	1517.74000000000
52.3570000000000	1517.73000000000
54.0340000000000	1517.71000000000
54.4240000000000	1517.70000000000
56.1090000000000	1517.70000000000
60.9930000000000	1517.68000000000
61.3820000000000	1517.68000000000
61.7690000000000	1517.68000000000
62.1580000000000	1517.68000000000
62.5490000000000	1517.67000000000
62.9420000000000	1517.66000000000
63.3330000000000	1517.60000000000
64.2820000000000	1517.58000000000
66.8620000000000	1517.58000000000
68.5320000000000	1517.58000000000
68.9220000000000	1517.58000000000
69.3130000000001	1517.58000000000
69.7030000000001	1517.58000000000
70.0960000000001	1517.58000000000
70.4860000000001	1517.58000000000
70.8760000000001	1517.58000000000
71.2650000000000	1517.57000000000
71.6590000000001	1517.57000000000
72.0490000000001	1517.57000000000
72.4430000000001	1517.57000000000
72.8340000000001	1517.57000000000
73.2240000000001	1517.57000000000
74.1720000000001	1517.57000000000
};
\addplot [line width=1, colQRkit]
table {%
  1.46700000000000	1884.92000000000
2.91800000000000	1737.81000000000
4.38800000000000	1666.85000000000
5.84400000000000	1629.68000000000
8.77500000000000	1606.70000000000
14.6270000000000	1593.65000000000
16.0840000000000	1576.90000000000
19.0160000000000	1566.94000000000
21.9390000000000	1559.61000000000
23.4010000000000	1556.85000000000
27.7870000000000	1547.83000000000
29.2480000000000	1539.91000000000
30.7070000000000	1538.29000000000
35.1500000000000	1530.91000000000
36.6800000000000	1527.57000000000
42.5260000000000	1525.52000000000
45.4380000000000	1522.41000000000
48.3610000000000	1520.38000000000
51.2830000000000	1518.37000000000
54.1940000000000	1516.53000000000
57.1130000000000	1515.37000000000
58.5770000000000	1513.81000000000
60.0380000000000	1512.22000000000
62.9640000000000	1508.87000000000
67.3550000000000	1506.96000000000
68.8140000000000	1505.26000000000
71.7250000000000	1504.13000000000
73.1840000000000	1502.78000000000
77.5680000000000	1501.27000000000
79.0280000000000	1500.38000000000
83.4220000000000	1497.44000000000
84.8900000000000	1496.24000000000
86.3520000000000	1495.73000000000
90.8050000000000	1494.28000000000
92.2620000000000	1492.95000000000
96.6340000000000	1491.62000000000
98.0860000000000	1490.33000000000
99.5640000000000	1489.90000000000
101.026000000000	1489.12000000000
103.945000000000	1486.73000000000
106.869000000000	1484.61000000000
109.792000000000	1483.74000000000
114.180000000000	1482.12000000000
118.562000000000	1480.22000000000
121.492000000000	1478.64000000000
125.863000000000	1477.93000000000
128.804000000000	1477.20000000000
131.714000000000	1476.73000000000
134.626000000000	1476.33000000000
136.081000000000	1475.95000000000
138.998000000000	1474.84000000000
140.451000000000	1474.06000000000
144.820000000000	1473.02000000000
147.737000000000	1471.56000000000
149.200000000000	1471.04000000000
152.118000000000	1470.49000000000
158.005000000000	1470.44000000000
162.370000000000	1470.25000000000
163.831000000000	1470.19000000000
165.287000000000	1470.13000000000
168.196000000000	1470.05000000000
169.653000000000	1469.99000000000
171.109000000000	1469.86000000000
172.559000000000	1469.79000000000
174.019000000000	1469.71000000000
175.476000000000	1469.55000000000
176.939000000000	1469.52000000000
178.405000000000	1469.34000000000
179.865000000000	1469.09000000000
181.322000000000	1468.90000000000
184.236000000000	1468.67000000000
185.696000000000	1468.53000000000
187.165000000000	1468.43000000000
188.619000000000	1468.37000000000
190.075000000000	1468.33000000000
191.538000000000	1468.27000000000
192.996000000000	1468.18000000000
194.445000000000	1468.11000000000
195.904000000000	1468.09000000000
198.826000000000	1468
201.747000000000	1467.96000000000
203.202000000000	1467.92000000000
204.657000000000	1467.88000000000
206.120000000000	1467.86000000000
207.575000000000	1467.73000000000
209.035000000000	1467.59000000000
210.493000000000	1467.56000000000
213.415000000000	1467.37000000000
214.866000000000	1467.33000000000
216.325000000000	1467.27000000000
219.248000000000	1467.25000000000
220.700000000000	1467.23000000000
222.157000000000	1467.21000000000
223.613000000000	1467.20000000000
226.523000000000	1467.09000000000
228.058000000000	1467.08000000000
231.011000000000	1466.97000000000
232.470000000000	1466.97000000000
233.924000000000	1466.96000000000
235.374000000000	1466.94000000000
236.826000000000	1466.93000000000
241.204000000000	1466.89000000000
242.661000000000	1466.89000000000
244.124000000000	1466.89000000000
245.745000000000	1466.88000000000
247.212000000000	1466.88000000000
250.132000000000	1466.88000000000
251.587000000000	1466.88000000000
253.052000000000	1466.85000000000
254.515000000000	1466.81000000000
255.970000000000	1466.81000000000
257.428000000000	1466.80000000000
258.877000000000	1466.80000000000
261.796000000000	1466.78000000000
263.250000000000	1466.78000000000
264.705000000000	1466.77000000000
266.177000000000	1466.77000000000
269.094000000000	1466.76000000000
270.550000000000	1466.76000000000
271.997000000000	1466.76000000000
273.456000000000	1466.76000000000
276.380000000000	1466.70000000000
277.836000000000	1466.67000000000
279.292000000000	1466.66000000000
282.210000000000	1466.63000000000
283.667000000000	1466.62000000000
285.121000000000	1466.62000000000
288.049000000000	1466.62000000000
290.966000000000	1466.61000000000
293.877000000000	1466.61000000000
295.334000000000	1466.60000000000
296.807000000000	1466.60000000000
298.261000000000	1466.60000000000
299.718000000000	1466.60000000000
302.640000000000	1466.60000000000
304.095000000000	1466.60000000000
307.013000000000	1466.60000000000
308.465000000000	1466.60000000000
312.898000000000	1466.59000000000
317.272000000000	1466.59000000000
318.778000000000	1466.59000000000
330.533000000000	1466.59000000000
};
\addplot [line width=1, colQRChol]
table {%
1.14700000000000	1884.92000000000
2.17400000000000	1736.18000000000
3.20100000000000	1663.65000000000
4.24500000000000	1627.78000000000
5.26700000000000	1605.12000000000
6.30600000000000	1588.68000000000
7.32300000000000	1575.62000000000
8.35900000000000	1570.17000000000
9.37900000000000	1555.11000000000
10.4150000000000	1545.70000000000
11.4340000000000	1540.51000000000
12.4650000000000	1534.18000000000
13.4790000000000	1531.18000000000
14.5140000000000	1525.57000000000
15.5290000000000	1524.62000000000
16.5570000000000	1523.10000000000
17.5730000000000	1522.18000000000
18.6120000000000	1521.49000000000
44.8270000000000	1520.51000000000
45.8570000000000	1519.27000000000
46.8760000000000	1518.26000000000
47.9050000000000	1517.34000000000
48.9210000000000	1517.17000000000
49.9490000000000	1516.23000000000
50.9630000000000	1515.85000000000
51.9870000000000	1515.36000000000
53.0040000000000	1514.94000000000
54.0350000000000	1514.84000000000
55.0540000000000	1514.19000000000
56.0830000000000	1513.88000000000
59.0180000000000	1513.29000000000
61.9600000000000	1513.07000000000
64.9230000000000	1512.32000000000
65.9380000000000	1511.92000000000
66.9640000000000	1511.78000000000
67.9760000000000	1511.08000000000
70.9130000000000	1510.74000000000
73.8640000000000	1510.42000000000
74.8890000000000	1510.17000000000
75.8960000000000	1509.77000000000
76.9240000000000	1509.63000000000
77.9350000000000	1508.89000000000
78.9650000000000	1508.63000000000
79.9800000000000	1508.59000000000
81.0110000000000	1508.11000000000
82.0230000000000	1507.57000000000
83.0540000000000	1507.16000000000
85.9900000000000	1506.85000000000
87.0040000000000	1506.52000000000
89.9620000000000	1506.02000000000
90.9840000000000	1505.65000000000
91.9960000000000	1505.20000000000
93.0240000000000	1504.93000000000
95.9660000000000	1504.25000000000
96.9770000000000	1504.13000000000
98.0010000000000	1504
99.0090000000000	1503.56000000000
100.037000000000	1503.39000000000
101.046000000000	1503.34000000000
104.006000000000	1502.81000000000
105.037000000000	1502.49000000000
106.046000000000	1502.33000000000
107.075000000000	1502.04000000000
110.021000000000	1501.93000000000
111.036000000000	1501.57000000000
112.062000000000	1501.12000000000
113.081000000000	1501.12000000000
116.032000000000	1500.52000000000
117.068000000000	1500.23000000000
120.016000000000	1500
121.032000000000	1499.91000000000
122.060000000000	1499.73000000000
123.072000000000	1499.70000000000
126.026000000000	1499.57000000000
127.049000000000	1499.40000000000
132.821000000000	1499.23000000000
135.763000000000	1499.04000000000
141.562000000000	1498.92000000000
144.506000000000	1498.86000000000
145.514000000000	1498.82000000000
146.544000000000	1498.70000000000
147.556000000000	1498.60000000000
148.580000000000	1498.56000000000
149.608000000000	1498.50000000000
150.632000000000	1498.47000000000
153.572000000000	1498.33000000000
156.520000000000	1498.27000000000
159.476000000000	1498.24000000000
165.262000000000	1498.21000000000
166.286000000000	1498.20000000000
169.226000000000	1498.18000000000
170.243000000000	1498.16000000000
171.262000000000	1498.07000000000
172.280000000000	1497.98000000000
178.123000000000	1497.74000000000
179.145000000000	1497.69000000000
180.168000000000	1497.67000000000
181.185000000000	1497.66000000000
182.246000000000	1497.61000000000
183.263000000000	1497.55000000000
184.329000000000	1497.41000000000
190.104000000000	1497.28000000000
193.060000000000	1497.13000000000
198.863000000000	1496.93000000000
199.879000000000	1496.86000000000
200.899000000000	1496.82000000000
206.673000000000	1496.74000000000
209.622000000000	1496.66000000000
210.640000000000	1496.65000000000
216.421000000000	1496.58000000000
219.353000000000	1496.55000000000
220.361000000000	1496.55000000000
223.307000000000	1496.54000000000
224.333000000000	1496.52000000000
230.105000000000	1496.52000000000
231.131000000000	1496.50000000000
232.146000000000	1496.48000000000
233.174000000000	1496.48000000000
236.115000000000	1496.47000000000
237.122000000000	1496.47000000000
238.146000000000	1496.46000000000
239.161000000000	1496.44000000000
240.188000000000	1496.35000000000
241.200000000000	1496.35000000000
244.160000000000	1496.34000000000
245.184000000000	1496.33000000000
248.115000000000	1496.31000000000
249.125000000000	1496.30000000000
254.894000000000	1496.26000000000
264.375000000000	1496.22000000000
265.388000000000	1496.19000000000
266.413000000000	1496.17000000000
269.349000000000	1496.16000000000
270.370000000000	1496.14000000000
273.321000000000	1496.11000000000
274.343000000000	1496.10000000000
275.357000000000	1496.05000000000
276.383000000000	1496.04000000000
277.404000000000	1496.04000000000
278.428000000000	1495.97000000000
279.439000000000	1495.90000000000
282.388000000000	1495.78000000000
283.410000000000	1495.78000000000
284.419000000000	1495.69000000000
285.437000000000	1495.68000000000
291.199000000000	1495.64000000000
300.744000000000	1495.57000000000
301.766000000000	1495.54000000000
302.791000000000	1495.51000000000
303.815000000000	1495.47000000000
304.831000000000	1495.39000000000
305.848000000000	1495.38000000000
306.870000000000	1495.38000000000
307.897000000000	1495.35000000000
308.917000000000	1495.33000000000
309.946000000000	1495.28000000000
310.964000000000	1495.26000000000
313.906000000000	1495.18000000000
316.864000000000	1495.17000000000
317.888000000000	1495.14000000000
323.648000000000	1495.09000000000
324.667000000000	1495.07000000000
325.683000000000	1495.06000000000
328.637000000000	1495.04000000000
338.155000000000	1495.03000000000
339.170000000000	1495.02000000000
340.194000000000	1495.02000000000
343.123000000000	1495.01000000000
352.726000000000	1495.01000000000
358.761000000000	1495
368.289000000000	1495
369.369000000000	1495
370.395000000000	1495
371.421000000000	1494.94000000000
374.353000000000	1494.93000000000
383.825000000000	1494.93000000000
403.592000000000	1494.93000000000
404.623000000000	1494.93000000000
405.635000000000	1494.93000000000
406.655000000000	1494.93000000000
409.662000000000	1494.93000000000
412.582000000000	1494.93000000000
};
\addplot [line width=1, colMoreQR]
table {%
2.21100000000000	1884.92000000000
4.19000000000000	1743.80000000000
6.16400000000000	1675.47000000000
8.14600000000000	1639.90000000000
10.1070000000000	1617.67000000000
12.0630000000000	1600.02000000000
15.3090000000000	1585.61000000000
18.5380000000000	1574.41000000000
20.4930000000000	1569.87000000000
23.7380000000000	1556.59000000000
25.7010000000000	1551.20000000000
27.6690000000000	1549.12000000000
33.4980000000000	1531.61000000000
35.4550000000000	1527.10000000000
37.4140000000000	1524.21000000000
42.0890000000000	1521.65000000000
44.0930000000000	1520.53000000000
46.0620000000000	1515.57000000000
50.5930000000000	1513.27000000000
53.8440000000000	1510.08000000000
57.0850000000000	1509.11000000000
60.3210000000000	1505.57000000000
62.2880000000000	1502.32000000000
66.8190000000000	1500.14000000000
68.7820000000000	1498.36000000000
72.0270000000000	1496.76000000000
76.5420000000000	1495.59000000000
78.5030000000000	1493.60000000000
81.7410000000000	1491.70000000000
86.2530000000000	1490.63000000000
88.2150000000000	1489.79000000000
90.1690000000000	1489.11000000000
94.6840000000000	1487.75000000000
96.6520000000000	1486.83000000000
98.6160000000000	1485.56000000000
103.151000000000	1484.73000000000
106.506000000000	1484.27000000000
108.494000000000	1481.84000000000
113.027000000000	1480.40000000000
116.249000000000	1479.19000000000
118.207000000000	1478.55000000000
121.453000000000	1476.83000000000
124.701000000000	1476.64000000000
127.948000000000	1474.92000000000
129.909000000000	1474.01000000000
131.860000000000	1473.14000000000
133.824000000000	1472.38000000000
137.062000000000	1471.41000000000
139.016000000000	1470.79000000000
143.530000000000	1469.30000000000
145.489000000000	1468.62000000000
149.999000000000	1468.18000000000
151.938000000000	1467.79000000000
153.882000000000	1467.55000000000
155.834000000000	1467.31000000000
160.335000000000	1467.04000000000
162.292000000000	1466.59000000000
165.525000000000	1466.52000000000
167.478000000000	1466.25000000000
169.449000000000	1465.98000000000
172.687000000000	1465.66000000000
174.636000000000	1465.23000000000
176.587000000000	1464.92000000000
178.532000000000	1464.89000000000
184.328000000000	1464.50000000000
190.109000000000	1464.41000000000
192.060000000000	1464.26000000000
196.576000000000	1464.10000000000
198.520000000000	1464.08000000000
201.736000000000	1463.92000000000
203.688000000000	1463.89000000000
205.636000000000	1463.82000000000
210.147000000000	1463.65000000000
215.913000000000	1463.63000000000
220.428000000000	1463.63000000000
222.383000000000	1463.62000000000
224.330000000000	1463.60000000000
226.294000000000	1463.58000000000
228.241000000000	1463.58000000000
231.477000000000	1463.57000000000
234.697000000000	1463.56000000000
236.654000000000	1463.54000000000
239.888000000000	1463.54000000000
241.846000000000	1463.53000000000
243.799000000000	1463.53000000000
245.751000000000	1463.52000000000
247.801000000000	1463.51000000000
252.297000000000	1463.51000000000
255.527000000000	1463.51000000000
258.760000000000	1463.51000000000
261.990000000000	1463.51000000000
263.941000000000	1463.51000000000
267.176000000000	1463.50000000000
269.140000000000	1463.50000000000
271.091000000000	1463.50000000000
273.041000000000	1463.50000000000
274.991000000000	1463.50000000000
278.206000000000	1463.50000000000
280.156000000000	1463.50000000000
282.108000000000	1463.50000000000
284.067000000000	1463.50000000000
287.324000000000	1463.50000000000
};
\end{axis}

\end{tikzpicture}

%% file: dubrovnik_iter_double.tex
%
%
%
\begin{tikzpicture}

\pgfplotsset{compat=newest} 

\definecolor{color0}{rgb}{1, 1, 1}

\tikzstyle{every node}=[font=\footnotesize]

\definecolor{colCholesky}{RGB}{0, 255, 0}
\definecolor{colQRkit}{RGB}{255, 0, 0}
\definecolor{colQRChol}{RGB}{0, 0, 255}
\definecolor{colSPQR}{RGB}{0, 255, 255}
\definecolor{colSSBA}{RGB}{255, 255, 0}
\definecolor{colMoreQR}{RGB}{200, 0, 255}

\begin{axis}[
xlabel={Iterations},
ylabel={Energy},
xmin=0, xmax=250,
ymin=3000, ymax=4500,
width=\figurewidth,
height=\figureheight,
at={(0\figurewidth,0\figureheight)},
xmajorgrids,
x grid style={lightgray},
ymajorgrids,
y grid style={lightgray},
axis line style={black},
axis background/.style={fill=color0},
legend style={at={(0.97,0.97)}, anchor=north east},
legend cell align={left},
legend entries={{Cholesky},{QRkit},{QRkit + Cholesky},{SuiteSparse QR},{SSBA},{Mor{\'e} QRkit}}
]
%
%
\addplot [line width=1, colCholesky]
table {%
1	4529.42000000000
2	4419.14000000000
3	4350.72000000000
4	4294.49000000000
5	4254.80000000000
6	4179.39000000000
7	4121.36000000000
8	4058.63000000000
9	3994.87000000000
10	3927.32000000000
11	3862.32000000000
12	3799.52000000000
13	3740.75000000000
14	3693.61000000000
15	3651.06000000000
16	3610.38000000000
17	3571.86000000000
18	3537.47000000000
19	3505.20000000000
20	3474.06000000000
21	3446.91000000000
22	3421.52000000000
23	3398.87000000000
24	3379.25000000000
25	3361.90000000000
26	3346.97000000000
27	3333.44000000000
28	3321.12000000000
29	3309.64000000000
30	3299.40000000000
31	3289.53000000000
32	3279.38000000000
33	3269.81000000000
34	3262.68000000000
35	3256.23000000000
36	3249.39000000000
37	3243.61000000000
38	3237.78000000000
39	3232.99000000000
40	3228.74000000000
41	3224.97000000000
42	3222.10000000000
43	3219.61000000000
44	3217.11000000000
45	3214.50000000000
46	3212.62000000000
47	3210.85000000000
48	3209.39000000000
49	3207.62000000000
50	3205.97000000000
51	3204.50000000000
52	3203.44000000000
53	3202.31000000000
54	3200.89000000000
55	3199.52000000000
56	3198.33000000000
57	3197.28000000000
58	3196.05000000000
59	3195.16000000000
60	3194.29000000000
61	3193.47000000000
62	3192.59000000000
63	3191.67000000000
64	3190.90000000000
65	3189.94000000000
66	3189.13000000000
67	3188.28000000000
68	3187.53000000000
69	3186.80000000000
70	3186.04000000000
71	3185.58000000000
72	3185.05000000000
73	3184.54000000000
74	3184.10000000000
75	3183.64000000000
76	3183.50000000000
77	3183.25000000000
78	3182.95000000000
79	3182.76000000000
80	3182.61000000000
81	3182.26000000000
82	3181.93000000000
83	3181.51000000000
84	3181.05000000000
85	3180.88000000000
86	3180.54000000000
87	3180.22000000000
88	3179.96000000000
89	3179.71000000000
90	3179.46000000000
91	3179.15000000000
92	3178.89000000000
93	3178.70000000000
94	3178.57000000000
95	3178.40000000000
96	3178.25000000000
97	3177.98000000000
98	3177.66000000000
99	3177.49000000000
100	3177.35000000000
101	3177.13000000000
102	3177.03000000000
103	3176.88000000000
104	3176.79000000000
105	3176.75000000000
106	3176.68000000000
107	3176.49000000000
108	3176.35000000000
109	3176.32000000000
110	3176.26000000000
111	3176.12000000000
112	3175.88000000000
113	3175.68000000000
114	3175.49000000000
115	3175.41000000000
116	3175.34000000000
117	3175.12000000000
118	3174.94000000000
119	3174.88000000000
120	3174.85000000000
121	3174.83000000000
122	3174.79000000000
123	3174.70000000000
124	3174.64000000000
125	3174.63000000000
126	3174.56000000000
127	3174.53000000000
128	3174.47000000000
129	3174.40000000000
130	3174.39000000000
131	3174.38000000000
132	3174.30000000000
133	3174.26000000000
134	3174.18000000000
135	3174.13000000000
136	3174.02000000000
137	3174.01000000000
138	3173.99000000000
139	3173.96000000000
140	3173.89000000000
141	3173.87000000000
142	3173.84000000000
143	3173.80000000000
144	3173.78000000000
145	3173.71000000000
146	3173.68000000000
147	3173.64000000000
148	3173.63000000000
149	3173.62000000000
150	3173.56000000000
151	3173.49000000000
152	3173.47000000000
153	3173.45000000000
154	3173.36000000000
155	3173.22000000000
156	3173.09000000000
157	3173
158	3172.97000000000
159	3172.96000000000
160	3172.91000000000
161	3172.75000000000
162	3172.71000000000
163	3172.68000000000
164	3172.65000000000
165	3172.64000000000
166	3172.62000000000
167	3172.60000000000
168	3172.60000000000
};
\addplot [line width=1, colQRkit]
table {%
1	4529.42000000000
2	4419.14000000000
3	4350.72000000000
4	4294.49000000000
5	4254.80000000000
6	4179.39000000000
7	4121.36000000000
8	4058.62000000000
9	3994.86000000000
10	3927.26000000000
11	3862.23000000000
12	3799.50000000000
13	3740.69000000000
14	3693.51000000000
15	3650.87000000000
16	3610.23000000000
17	3571.28000000000
18	3537.09000000000
19	3504.30000000000
20	3473.41000000000
21	3446.64000000000
22	3421.61000000000
23	3399.25000000000
24	3380.32000000000
25	3364.30000000000
26	3349.90000000000
27	3337.11000000000
28	3324.50000000000
29	3312.46000000000
30	3301.78000000000
31	3291.37000000000
32	3280.83000000000
33	3271.19000000000
34	3263.29000000000
35	3256.07000000000
36	3248.70000000000
37	3242.90000000000
38	3237.67000000000
39	3233.59000000000
40	3229.70000000000
41	3225.58000000000
42	3222.45000000000
43	3220.11000000000
44	3217.87000000000
45	3216.04000000000
46	3214.49000000000
47	3212.37000000000
48	3210.95000000000
49	3210.73000000000
50	3206.17000000000
51	3204.10000000000
52	3202.22000000000
53	3200.86000000000
54	3199.44000000000
55	3198.01000000000
56	3196.77000000000
57	3195.70000000000
58	3194.67000000000
59	3193.84000000000
60	3192.65000000000
61	3191.49000000000
62	3190.31000000000
63	3189.12000000000
64	3188.36000000000
65	3187.74000000000
66	3187.01000000000
67	3186.28000000000
68	3185.69000000000
69	3184.66000000000
70	3184.04000000000
71	3183.55000000000
72	3183.06000000000
73	3182.71000000000
74	3182.31000000000
75	3181.97000000000
76	3181.52000000000
77	3181.13000000000
78	3180.71000000000
79	3180.11000000000
80	3179.62000000000
81	3179.28000000000
82	3179.09000000000
83	3178.94000000000
84	3178.77000000000
85	3178.63000000000
86	3178.54000000000
87	3178.48000000000
88	3178.41000000000
89	3178.27000000000
90	3178.07000000000
91	3178
92	3177.94000000000
93	3177.84000000000
94	3177.73000000000
95	3177.60000000000
96	3177.55000000000
97	3177.46000000000
98	3177.42000000000
99	3177.36000000000
100	3177.26000000000
101	3177.12000000000
102	3177.03000000000
103	3176.76000000000
104	3176.68000000000
105	3176.62000000000
106	3176.54000000000
107	3176.48000000000
108	3176.42000000000
109	3176.37000000000
110	3176.33000000000
111	3176.31000000000
112	3176.29000000000
113	3176.26000000000
114	3176.24000000000
115	3176.12000000000
116	3175.91000000000
117	3175.82000000000
118	3175.80000000000
119	3175.78000000000
120	3175.75000000000
121	3175.72000000000
122	3175.70000000000
123	3175.63000000000
124	3175.57000000000
125	3175.56000000000
126	3175.56000000000
127	3175.52000000000
128	3175.51000000000
129	3175.50000000000
130	3175.50000000000
131	3175.49000000000
132	3175.40000000000
133	3175.32000000000
134	3175.28000000000
135	3175.23000000000
136	3175.22000000000
137	3175.22000000000
138	3175.22000000000
139	3175.21000000000
140	3175.21000000000
141	3175.18000000000
142	3175.16000000000
143	3175.07000000000
144	3175.05000000000
145	3174.89000000000
146	3174.78000000000
147	3174.38000000000
148	3174.31000000000
149	3174.29000000000
150	3174.27000000000
151	3174.20000000000
152	3174.05000000000
153	3173.99000000000
154	3173.81000000000
155	3173.72000000000
156	3173.68000000000
157	3173.64000000000
158	3173.55000000000
159	3173.53000000000
160	3173.52000000000
161	3173.51000000000
162	3173.51000000000
163	3173.49000000000
164	3173.46000000000
165	3173.32000000000
166	3173.27000000000
167	3173.27000000000
168	3173.27000000000
169	3173.27000000000
170	3173.26000000000
171	3173.26000000000
172	3173.26000000000
173	3173.25000000000
174	3173.24000000000
175	3173.11000000000
176	3173.06000000000
177	3173.06000000000
178	3173.01000000000
179	3172.93000000000
180	3172.89000000000
181	3172.87000000000
182	3172.87000000000
183	3172.87000000000
184	3172.87000000000
185	3172.87000000000
186	3172.87000000000
187	3172.87000000000
188	3172.84000000000
189	3172.77000000000
190	3172.74000000000
191	3172.59000000000
192	3172.46000000000
193	3172.38000000000
194	3172.23000000000
195	3172.22000000000
196	3172.21000000000
197	3172.13000000000
198	3172.02000000000
199	3171.91000000000
200	3171.90000000000
201	3171.89000000000
202	3171.89000000000
203	3171.83000000000
204	3171.83000000000
205	3171.83000000000
206	3171.80000000000
207	3171.78000000000
208	3171.78000000000
209	3171.78000000000
210	3171.78000000000
211	3171.78000000000
212	3171.78000000000
213	3171.78000000000
214	3171.78000000000
215	3171.78000000000
216	3171.78000000000
217	3171.78000000000
218	3171.77000000000
219	3171.77000000000
220	3171.77000000000
221	3171.77000000000
222	3171.77000000000
};
\addplot [line width=1, colQRChol]
table {%
1	4529.42000000000
2	4419.14000000000
3	4350.72000000000
4	4294.49000000000
5	4254.80000000000
6	4179.39000000000
7	4121.36000000000
8	4058.62000000000
9	3994.86000000000
10	3927.27000000000
11	3862.23000000000
12	3799.48000000000
13	3740.67000000000
14	3693.51000000000
15	3650.95000000000
16	3610.57000000000
17	3571.86000000000
18	3536.89000000000
19	3504.05000000000
20	3473.23000000000
21	3446.79000000000
22	3421.37000000000
23	3398.23000000000
24	3378.25000000000
25	3360.93000000000
26	3346.37000000000
27	3333.03000000000
28	3319.88000000000
29	3308.45000000000
30	3298.05000000000
31	3287.91000000000
32	3278.38000000000
33	3269.93000000000
34	3262.60000000000
35	3256.11000000000
36	3249.12000000000
37	3242.85000000000
38	3241.16000000000
39	3230.93000000000
40	3225.79000000000
41	3224.49000000000
42	3217.81000000000
43	3214.82000000000
44	3212.57000000000
45	3210.79000000000
46	3209.32000000000
47	3206.81000000000
48	3204.66000000000
49	3203.08000000000
50	3201.67000000000
51	3199.88000000000
52	3198.55000000000
53	3197.46000000000
54	3196.42000000000
55	3195.07000000000
56	3193.96000000000
57	3193.68000000000
58	3188.29000000000
59	3186.36000000000
60	3184.87000000000
61	3183.56000000000
62	3182.04000000000
63	3181.09000000000
64	3180.17000000000
65	3179.26000000000
66	3176.19000000000
67	3174.88000000000
68	3173.98000000000
69	3173.27000000000
70	3172.56000000000
71	3171.86000000000
72	3171.06000000000
73	3170.31000000000
74	3170.04000000000
75	3169.70000000000
76	3169.38000000000
77	3168.79000000000
78	3168.56000000000
79	3168.23000000000
80	3167.97000000000
81	3167.81000000000
82	3167.50000000000
83	3167.20000000000
84	3167.03000000000
85	3165.58000000000
86	3164.64000000000
87	3164.22000000000
88	3164.05000000000
89	3163.93000000000
90	3163.76000000000
91	3163.68000000000
92	3163.61000000000
93	3163.54000000000
94	3163.31000000000
95	3163.11000000000
96	3162.99000000000
97	3162.82000000000
98	3162.62000000000
99	3162.46000000000
100	3162.38000000000
101	3162.32000000000
102	3162.29000000000
103	3162.22000000000
104	3162.01000000000
105	3161.88000000000
106	3161.82000000000
107	3161.80000000000
108	3161.79000000000
109	3161.78000000000
110	3161.78000000000
111	3161.74000000000
112	3161.71000000000
113	3161.70000000000
114	3161.69000000000
115	3161.68000000000
116	3161.66000000000
117	3161.64000000000
118	3161.61000000000
119	3161.57000000000
120	3161.55000000000
121	3161.51000000000
122	3161.41000000000
123	3161.39000000000
124	3161.35000000000
125	3161.33000000000
126	3161.32000000000
127	3161.29000000000
128	3161.23000000000
129	3161.20000000000
130	3161.16000000000
131	3161.16000000000
132	3161.16000000000
133	3161.16000000000
134	3161.16000000000
135	3161.16000000000
136	3161.16000000000
137	3161.16000000000
138	3161.15000000000
139	3161.15000000000
140	3161.15000000000
141	3161.15000000000
142	3161.15000000000
143	3161.15000000000
144	3161.15000000000
145	3161.14000000000
146	3161.14000000000
147	3161.14000000000
148	3161.13000000000
149	3161.13000000000
150	3161.13000000000
151	3161.11000000000
152	3161.07000000000
153	3161.06000000000
154	3161
155	3160.92000000000
156	3160.91000000000
157	3160.88000000000
158	3160.82000000000
159	3160.76000000000
160	3160.71000000000
161	3160.71000000000
162	3160.71000000000
163	3160.70000000000
164	3160.50000000000
165	3160.38000000000
166	3160.38000000000
167	3160.37000000000
168	3160.37000000000
169	3160.36000000000
170	3160.36000000000
171	3160.36000000000
172	3160.36000000000
173	3160.36000000000
174	3160.35000000000
175	3160.35000000000
176	3160.18000000000
177	3160.13000000000
178	3160.13000000000
179	3160.13000000000
180	3160.13000000000
181	3160.13000000000
182	3160.13000000000
183	3160.13000000000
184	3160.13000000000
185	3160.13000000000
186	3160.13000000000
};
\addplot [line width=1, colSPQR]
table {%
1	4529.42000000000
2	4419.14000000000
3	4350.72000000000
4	4294.49000000000
5	4254.80000000000
6	4179.39000000000
7	4121.36000000000
8	4058.62000000000
9	3994.86000000000
10	3927.26000000000
11	3862.23000000000
12	3799.50000000000
13	3740.69000000000
14	3693.51000000000
15	3650.87000000000
16	3610.02000000000
17	3571.16000000000
18	3536.64000000000
19	3503.15000000000
20	3471.66000000000
21	3444.36000000000
22	3416.68000000000
23	3393.72000000000
24	3374.23000000000
25	3357.09000000000
26	3341.80000000000
27	3327.45000000000
28	3314.18000000000
29	3302.89000000000
30	3291.72000000000
31	3282
32	3273.57000000000
33	3261.27000000000
34	3253.69000000000
35	3246.36000000000
36	3241.40000000000
37	3232.45000000000
38	3226.97000000000
39	3222.78000000000
40	3218.08000000000
41	3214.90000000000
42	3209.97000000000
43	3207.24000000000
44	3203.31000000000
45	3200.48000000000
46	3197.45000000000
47	3194.97000000000
48	3192.91000000000
49	3191.03000000000
50	3189.35000000000
51	3187.88000000000
52	3186.59000000000
53	3185.31000000000
54	3183.84000000000
55	3182.16000000000
56	3180.94000000000
57	3180
58	3179.07000000000
59	3177.49000000000
60	3176.13000000000
61	3174.88000000000
62	3174.06000000000
63	3173.09000000000
64	3171.90000000000
65	3170.66000000000
66	3169.85000000000
67	3169.07000000000
68	3167.73000000000
69	3166.94000000000
70	3166.11000000000
71	3164.89000000000
72	3163.86000000000
73	3163.05000000000
74	3162.26000000000
75	3161.47000000000
76	3160.85000000000
77	3160.24000000000
78	3159.74000000000
79	3159.32000000000
80	3158.93000000000
81	3158.74000000000
82	3158.56000000000
83	3158.48000000000
84	3158.32000000000
85	3158.13000000000
86	3157.91000000000
87	3157.46000000000
88	3157.04000000000
89	3156.83000000000
90	3156.73000000000
91	3156.64000000000
92	3156.47000000000
93	3156.22000000000
94	3156.12000000000
95	3156.02000000000
96	3155.95000000000
97	3155.89000000000
98	3155.81000000000
99	3155.71000000000
100	3155.62000000000
101	3155.54000000000
102	3155.47000000000
103	3155.43000000000
104	3155.40000000000
105	3155.35000000000
106	3155.26000000000
107	3155.19000000000
108	3155.14000000000
109	3155.07000000000
110	3155
111	3154.96000000000
112	3154.93000000000
113	3154.90000000000
114	3154.88000000000
115	3154.86000000000
116	3154.84000000000
117	3154.80000000000
118	3154.76000000000
119	3154.73000000000
120	3154.72000000000
121	3154.71000000000
122	3154.68000000000
123	3154.47000000000
124	3154.31000000000
125	3154.26000000000
126	3154.16000000000
127	3154.13000000000
128	3154.11000000000
129	3154.01000000000
130	3153.97000000000
131	3153.85000000000
132	3153.83000000000
133	3153.82000000000
134	3153.81000000000
135	3153.81000000000
136	3153.81000000000
137	3153.80000000000
138	3153.78000000000
139	3153.69000000000
140	3153.65000000000
141	3153.65000000000
142	3153.61000000000
143	3153.56000000000
144	3153.54000000000
145	3153.30000000000
146	3153.24000000000
147	3153.24000000000
148	3153.23000000000
149	3153.21000000000
150	3153.15000000000
151	3153.04000000000
152	3153
153	3152.99000000000
154	3152.94000000000
155	3152.86000000000
156	3152.81000000000
157	3152.80000000000
158	3152.79000000000
159	3152.66000000000
160	3152.61000000000
161	3152.60000000000
162	3152.50000000000
163	3152.49000000000
164	3152.47000000000
165	3152.46000000000
166	3152.46000000000
167	3152.45000000000
168	3152.45000000000
169	3152.45000000000
170	3152.45000000000
171	3152.40000000000
172	3152.39000000000
173	3152.39000000000
174	3152.38000000000
175	3152.38000000000
176	3152.38000000000
177	3152.38000000000
178	3152.38000000000
179	3152.38000000000
180	3152.38000000000
181	3152.38000000000
182	3152.38000000000
183	3152.38000000000
184	3152.38000000000
185	3152.38000000000
186	3152.38000000000
187	3152.38000000000
188	3152.36000000000
189	3152.32000000000
190	3152.32000000000
191	3152.32000000000
192	3152.32000000000
193	3152.32000000000
194	3152.32000000000
195	3152.32000000000
196	3152.32000000000
197	3152.32000000000
198	3152.32000000000
199	3152.32000000000
200	3152.32000000000
201	3152.32000000000
202	3152.30000000000
203	3152.28000000000
204	3152.28000000000
205	3152.28000000000
206	3152.27000000000
207	3152.27000000000
208	3152.27000000000
209	3152.25000000000
210	3152.22000000000
211	3152.22000000000
212	3152.21000000000
213	3152.21000000000
214	3152.21000000000
215	3152.21000000000
216	3152.21000000000
217	3152.21000000000
218	3152.20000000000
219	3152.18000000000
220	3152.17000000000
221	3152.17000000000
222	3152.17000000000
223	3152.17000000000
224	3152.14000000000
225	3152.13000000000
226	3152.12000000000
227	3152.12000000000
228	3152.12000000000
229	3152.12000000000
230	3152.09000000000
231	3152.04000000000
232	3151.86000000000
233	3151.75000000000
234	3151.69000000000
235	3151.67000000000
236	3151.67000000000
237	3151.67000000000
238	3151.67000000000
239	3151.67000000000
240	3151.67000000000
241	3151.67000000000
242	3151.67000000000
243	3151.67000000000
244	3151.67000000000
245	3151.67000000000
246	3151.67000000000
247	3151.67000000000
248	3151.66000000000
249	3151.66000000000
250	3151.66000000000
251	3151.66000000000
252	3151.66000000000
253	3151.66000000000
254	3151.66000000000
255	3151.66000000000
256	3151.66000000000
257	3151.66000000000
258	3151.66000000000
259	3151.66000000000
260	3151.66000000000
261	3151.66000000000
262	3151.66000000000
263	3151.66000000000
264	3151.66000000000
265	3151.66000000000
266	3151.66000000000
267	3151.66000000000
268	3151.66000000000
269	3151.66000000000
270	3151.66000000000
271	3151.66000000000
272	3151.66000000000
273	3151.65000000000
274	3151.65000000000
275	3151.65000000000
276	3151.65000000000
277	3151.65000000000
278	3151.65000000000
279	3151.65000000000
280	3151.65000000000
281	3151.65000000000
282	3151.65000000000
283	3151.65000000000
284	3151.65000000000
285	3151.65000000000
286	3151.65000000000
287	3151.65000000000
288	3151.65000000000
289	3151.65000000000
290	3151.65000000000
291	3151.65000000000
292	3151.65000000000
293	3151.65000000000
294	3151.65000000000
295	3151.65000000000
296	3151.65000000000
297	3151.65000000000
298	3151.65000000000
299	3151.65000000000
300	3151.65000000000
301	3151.65000000000
302	3151.65000000000
303	3151.65000000000
304	3151.65000000000
305	3151.65000000000
306	3151.65000000000
307	3151.65000000000
308	3151.64000000000
309	3151.64000000000
310	3151.64000000000
311	3151.64000000000
312	3151.64000000000
313	3151.63000000000
314	3151.63000000000
315	3151.63000000000
316	3151.63000000000
317	3151.63000000000
318	3151.63000000000
319	3151.63000000000
320	3151.63000000000
321	3151.63000000000
322	3151.63000000000
323	3151.63000000000
324	3151.63000000000
325	3151.63000000000
326	3151.63000000000
327	3151.63000000000
328	3151.63000000000
329	3151.63000000000
330	3151.63000000000
331	3151.63000000000
332	3151.63000000000
};
\addplot [line width=1, colSSBA]
table {%
1	4529.42332500000
2	4423.38846000000
3	4331.87566100000
4	4324.49816200000
5	4226.41055500000
6	4132.25494400000
7	4053.44162700000
8	3982.60084700000
9	3915.75531100000
10	3853.29087000000
11	3793.92521500000
12	3742.98906400000
13	3699.40961600000
14	3658.55036800000
15	3623.04629600000
16	3589.26066900000
17	3556.41591300000
18	3523.40045900000
19	3493.04007200000
20	3464.10340400000
21	3436.12536400000
22	3411.70414700000
23	3390.93160000000
24	3371.96944000000
25	3354.51395200000
26	3338.34248600000
27	3323.14201800000
28	3310.22309000000
29	3297.52635700000
30	3285.81583600000
31	3275.78121200000
32	3266.83710900000
33	3258.62035100000
34	3250.89201400000
35	3244.93917100000
36	3239.33396300000
37	3234.15780500000
38	3230.03384300000
39	3226.35468600000
40	3222.99182700000
41	3219.95099200000
42	3217.11228200000
43	3214.06991800000
44	3211.06283100000
45	3208.70817900000
46	3206.50243400000
47	3203.75364000000
48	3202.00860200000
49	3199.78940400000
50	3197.34227800000
51	3195.63550400000
52	3194.05119500000
53	3192.84394100000
54	3191.81983000000
55	3190.46506400000
56	3189.33326600000
57	3188.56874600000
58	3188.00524000000
59	3187.24094100000
60	3186.57504600000
61	3185.90343700000
62	3184.90090500000
63	3184.14154200000
64	3183.29258500000
65	3182.84067300000
66	3182.47916400000
67	3182.12828500000
68	3181.68602200000
69	3181.18995000000
70	3180.64460500000
71	3180.12980900000
72	3179.74318400000
73	3179.36760400000
74	3179.14165000000
75	3179.00064300000
76	3178.79821400000
77	3178.59690000000
78	3178.40811800000
79	3178.09317000000
80	3177.73643600000
81	3177.51769400000
82	3177.25448400000
83	3176.96257700000
84	3176.50857300000
85	3176.26994900000
86	3176.00173800000
87	3175.81672500000
88	3175.65014100000
89	3175.58231600000
90	3175.10447800000
91	3174.96888500000
92	3174.68214200000
93	3174.09720100000
94	3173.77139000000
95	3173.50584500000
96	3173.37630100000
97	3173.17190900000
98	3173.06112800000
99	3172.89879400000
};
\addplot [line width=1, colMoreQR]
table {%
1	4529.42000000000
2	4433.03000000000
3	4382.58000000000
4	4340.66000000000
5	4298.34000000000
6	4252.94000000000
7	4221.29000000000
8	4160.80000000000
9	4113.30000000000
10	4052.55000000000
11	3990.26000000000
12	3926
13	3864.72000000000
14	3803.23000000000
15	3744.17000000000
16	3690.07000000000
17	3645.85000000000
18	3607.30000000000
19	3572.88000000000
20	3541.28000000000
21	3510.83000000000
22	3481.01000000000
23	3453.35000000000
24	3427.78000000000
25	3405.82000000000
26	3387.44000000000
27	3368.07000000000
28	3352.25000000000
29	3338.84000000000
30	3326.57000000000
31	3313.87000000000
32	3301.82000000000
33	3290.21000000000
34	3280.30000000000
35	3271.85000000000
36	3264.25000000000
37	3256.83000000000
38	3249.79000000000
39	3242.90000000000
40	3236.44000000000
41	3231.46000000000
42	3226.64000000000
43	3222.53000000000
44	3218.72000000000
45	3214.98000000000
46	3211.60000000000
47	3208.63000000000
48	3205.40000000000
49	3202.01000000000
50	3199.20000000000
51	3195.99000000000
52	3193.94000000000
53	3191.49000000000
54	3189.34000000000
55	3187.08000000000
56	3185.18000000000
57	3184.63000000000
58	3180.87000000000
59	3179.19000000000
60	3177.81000000000
61	3176.34000000000
62	3175.12000000000
63	3173.96000000000
64	3173.01000000000
65	3171.69000000000
66	3170.81000000000
67	3170.29000000000
68	3169.65000000000
69	3168.85000000000
70	3168.11000000000
71	3167.58000000000
72	3167.11000000000
73	3166.62000000000
74	3166.20000000000
75	3165.77000000000
76	3165.12000000000
77	3164.83000000000
78	3164.43000000000
79	3163.54000000000
80	3163.03000000000
81	3162.68000000000
82	3162.50000000000
83	3162.40000000000
84	3162.13000000000
85	3161.86000000000
86	3161.54000000000
87	3161.18000000000
88	3160.15000000000
89	3159.73000000000
90	3159.53000000000
91	3157.81000000000
92	3157.01000000000
93	3156.50000000000
94	3156.12000000000
95	3155.76000000000
96	3155.51000000000
97	3155.24000000000
98	3154.96000000000
99	3154.75000000000
100	3154.45000000000
101	3154.24000000000
102	3154.15000000000
103	3154.05000000000
104	3153.87000000000
105	3153.69000000000
106	3153.55000000000
107	3153.36000000000
108	3153.23000000000
109	3153.10000000000
110	3152.91000000000
111	3152.66000000000
112	3152.33000000000
113	3152.23000000000
114	3152.17000000000
115	3151.99000000000
116	3151.84000000000
117	3151.69000000000
118	3151.52000000000
119	3151.40000000000
120	3151.22000000000
121	3150.75000000000
122	3150.45000000000
123	3150.30000000000
124	3149.95000000000
125	3149.54000000000
126	3149.28000000000
127	3149.14000000000
128	3148.94000000000
129	3148.86000000000
130	3148.64000000000
131	3148.38000000000
132	3148.02000000000
133	3147.71000000000
134	3147.51000000000
135	3147.40000000000
136	3147.31000000000
137	3147.22000000000
138	3147.10000000000
139	3147.06000000000
140	3146.83000000000
141	3146.56000000000
142	3146.28000000000
143	3146.14000000000
144	3145.89000000000
145	3145.79000000000
146	3145.74000000000
147	3145.62000000000
148	3145.50000000000
149	3145.36000000000
150	3145.08000000000
151	3144.93000000000
152	3144.78000000000
153	3144.58000000000
154	3144.38000000000
155	3144.10000000000
156	3143.90000000000
157	3143.74000000000
158	3143.55000000000
159	3143.35000000000
160	3143.05000000000
161	3142.70000000000
162	3142.47000000000
163	3142.24000000000
164	3142.05000000000
165	3141.95000000000
166	3141.87000000000
167	3141.70000000000
168	3141.60000000000
169	3141.48000000000
170	3141.34000000000
171	3141.13000000000
172	3141
173	3140.84000000000
174	3140.62000000000
175	3140.37000000000
176	3140.27000000000
177	3140.13000000000
178	3140.05000000000
179	3139.93000000000
180	3139.75000000000
181	3139.54000000000
182	3139.40000000000
183	3139.25000000000
184	3139.14000000000
185	3138.92000000000
186	3138.83000000000
187	3138.80000000000
188	3138.78000000000
189	3138.77000000000
190	3138.73000000000
191	3138.61000000000
192	3138.49000000000
193	3138.31000000000
194	3138.23000000000
195	3138.17000000000
196	3137.94000000000
197	3137.83000000000
198	3137.75000000000
199	3137.69000000000
200	3137.61000000000
201	3137.52000000000
202	3137.30000000000
203	3137.10000000000
204	3136.95000000000
205	3136.79000000000
206	3136.74000000000
207	3136.73000000000
208	3136.67000000000
209	3136.39000000000
210	3136.25000000000
211	3136.18000000000
212	3136.17000000000
213	3136.16000000000
214	3136.14000000000
215	3136.13000000000
216	3136.10000000000
217	3136.05000000000
218	3135.90000000000
219	3135.74000000000
220	3135.71000000000
221	3135.64000000000
222	3135.58000000000
223	3135.47000000000
224	3135.29000000000
225	3134.74000000000
226	3134.63000000000
227	3134.59000000000
228	3134.53000000000
229	3134.38000000000
230	3134.29000000000
231	3134.26000000000
232	3134.20000000000
233	3134.14000000000
234	3133.97000000000
235	3133.91000000000
236	3133.79000000000
237	3133.66000000000
238	3133.62000000000
239	3133.61000000000
240	3133.56000000000
241	3133.52000000000
242	3133.45000000000
243	3133.28000000000
244	3133.19000000000
245	3133.14000000000
246	3133.06000000000
247	3133
248	3132.96000000000
249	3132.83000000000
250	3132.64000000000
251	3132.60000000000
252	3132.58000000000
253	3132.56000000000
254	3132.53000000000
255	3132.50000000000
256	3132.39000000000
257	3132.34000000000
258	3132.23000000000
259	3131.98000000000
260	3131.69000000000
261	3131.50000000000
262	3131.32000000000
263	3131.29000000000
264	3131.26000000000
265	3131.24000000000
266	3131.17000000000
267	3131.05000000000
268	3131.02000000000
269	3130.91000000000
270	3130.85000000000
271	3130.83000000000
272	3130.81000000000
273	3130.78000000000
274	3130.78000000000
275	3130.78000000000
276	3130.78000000000
277	3130.70000000000
278	3130.62000000000
279	3130.46000000000
280	3130.42000000000
281	3130.40000000000
282	3130.36000000000
283	3130.33000000000
284	3130.32000000000
285	3130.26000000000
286	3130.16000000000
287	3129.95000000000
288	3129.89000000000
289	3129.88000000000
290	3129.82000000000
291	3129.64000000000
292	3129.45000000000
293	3129.36000000000
294	3129.34000000000
295	3129.34000000000
296	3129.30000000000
297	3129.24000000000
298	3129.21000000000
299	3129.21000000000
300	3129.20000000000
301	3129.13000000000
302	3129.13000000000
303	3129.11000000000
304	3129.04000000000
305	3129.01000000000
306	3129
307	3128.96000000000
308	3128.85000000000
309	3128.81000000000
310	3128.78000000000
311	3128.74000000000
312	3128.71000000000
313	3128.70000000000
314	3128.59000000000
315	3128.50000000000
316	3128.33000000000
317	3128.26000000000
318	3128.20000000000
319	3128.10000000000
320	3128.07000000000
321	3128.07000000000
322	3128.07000000000
323	3128.07000000000
};
\end{axis}

\end{tikzpicture}

%% file: dubrovnik_double.tex
%
%
%
\begin{tikzpicture}

\pgfplotsset{compat=newest} 

\definecolor{color0}{rgb}{1, 1, 1}

\tikzstyle{every node}=[font=\footnotesize]

\definecolor{colCholesky}{RGB}{0, 255, 0}
\definecolor{colQRkit}{RGB}{255, 0, 0}
\definecolor{colQRChol}{RGB}{0, 0, 255}
\definecolor{colSPQR}{RGB}{0, 255, 255}
\definecolor{colSSBA}{RGB}{255, 255, 0}
\definecolor{colMoreQR}{RGB}{200, 0, 255}

\begin{axis}[
xmode=log,
log ticks with fixed point,
xtick={1,5,10,20,50,120,300,1000},
xlabel={Time [s]},
ylabel={Energy},
xmin=1, xmax=5e3,
ymin=3000, ymax=4500,
width=\figurewidth,
height=\figureheight,
at={(0\figurewidth,0\figureheight)},
xmajorgrids,
x grid style={lightgray},
ymajorgrids,
y grid style={lightgray},
axis line style={black},
axis background/.style={fill=color0},
]
%
%
\addplot [line width=0.75, colCholesky] 
table {%
0.602000000000000	4529.42000000000
1.18900000000000	4419.14000000000
1.78900000000000	4350.72000000000
2.39200000000000	4294.49000000000
2.98700000000000	4254.80000000000
3.59100000000000	4179.39000000000
4.19000000000000	4121.36000000000
4.79000000000000	4058.63000000000
5.38600000000000	3994.87000000000
5.98400000000000	3927.32000000000
6.57900000000000	3862.32000000000
7.17700000000000	3799.52000000000
7.77400000000000	3740.75000000000
8.37900000000000	3693.61000000000
8.97400000000000	3651.06000000000
9.57600000000000	3610.38000000000
10.1680000000000	3571.86000000000
10.7660000000000	3537.47000000000
11.3680000000000	3505.20000000000
11.9670000000000	3474.06000000000
12.5640000000000	3446.91000000000
13.1660000000000	3421.52000000000
13.7640000000000	3398.87000000000
14.3950000000000	3379.25000000000
14.9980000000000	3361.90000000000
15.5930000000000	3346.97000000000
16.1930000000000	3333.44000000000
16.7890000000000	3321.12000000000
17.4000000000000	3309.64000000000
17.9940000000000	3299.40000000000
18.5890000000000	3289.53000000000
19.1860000000000	3279.38000000000
19.7840000000000	3269.81000000000
20.3870000000000	3262.68000000000
20.9870000000000	3256.23000000000
21.5850000000000	3249.39000000000
22.1800000000000	3243.61000000000
24.5750000000000	3237.78000000000
25.1760000000000	3232.99000000000
25.7810000000000	3228.74000000000
26.3780000000000	3224.97000000000
26.9800000000000	3222.10000000000
27.5790000000000	3219.61000000000
28.1800000000000	3217.11000000000
28.7750000000000	3214.50000000000
29.3760000000000	3212.62000000000
29.9730000000000	3210.85000000000
32.3730000000000	3209.39000000000
32.9760000000000	3207.62000000000
33.5700000000000	3205.97000000000
34.1800000000000	3204.50000000000
34.7770000000000	3203.44000000000
35.3700000000000	3202.31000000000
35.9760000000000	3200.89000000000
36.5740000000000	3199.52000000000
37.1760000000000	3198.33000000000
37.7750000000000	3197.28000000000
38.3700000000000	3196.05000000000
38.9750000000000	3195.16000000000
39.5700000000000	3194.29000000000
40.1680000000000	3193.47000000000
40.7680000000000	3192.59000000000
41.3640000000000	3191.67000000000
41.9670000000000	3190.90000000000
42.5600000000000	3189.94000000000
43.1720000000000	3189.13000000000
43.7730000000000	3188.28000000000
44.3710000000000	3187.53000000000
44.9840000000000	3186.80000000000
45.5790000000000	3186.04000000000
46.1790000000000	3185.58000000000
46.7760000000000	3185.05000000000
47.3750000000000	3184.54000000000
47.9790000000000	3184.10000000000
48.5770000000000	3183.64000000000
49.1750000000000	3183.50000000000
49.7690000000000	3183.25000000000
50.3680000000000	3182.95000000000
50.9660000000000	3182.76000000000
51.5610000000000	3182.61000000000
52.1570000000000	3182.26000000000
52.7550000000000	3181.93000000000
53.3510000000000	3181.51000000000
53.9510000000000	3181.05000000000
54.5430000000000	3180.88000000000
55.1390000000000	3180.54000000000
55.7450000000000	3180.22000000000
56.3440000000000	3179.96000000000
56.9470000000000	3179.71000000000
57.5420000000000	3179.46000000000
58.1400000000000	3179.15000000000
61.7430000000000	3178.89000000000
62.3380000000000	3178.70000000000
62.9340000000000	3178.57000000000
63.5350000000000	3178.40000000000
64.1310000000000	3178.25000000000
64.7310000000000	3177.98000000000
65.3300000000000	3177.66000000000
65.9270000000000	3177.49000000000
66.5250000000000	3177.35000000000
67.1240000000000	3177.13000000000
67.7250000000000	3177.03000000000
68.3250000000000	3176.88000000000
68.9229999999999	3176.79000000000
69.5250000000000	3176.75000000000
70.1240000000000	3176.68000000000
70.7220000000000	3176.49000000000
71.3190000000000	3176.35000000000
71.9139999999999	3176.32000000000
72.5109999999999	3176.26000000000
73.1089999999999	3176.12000000000
73.7119999999999	3175.88000000000
74.3099999999999	3175.68000000000
74.9109999999999	3175.49000000000
75.5179999999999	3175.41000000000
76.1199999999999	3175.34000000000
76.7219999999999	3175.12000000000
77.3209999999999	3174.94000000000
77.9189999999999	3174.88000000000
78.5219999999999	3174.85000000000
79.1179999999999	3174.83000000000
79.7159999999999	3174.79000000000
80.3139999999999	3174.70000000000
80.9149999999999	3174.64000000000
81.5149999999999	3174.63000000000
82.1109999999999	3174.56000000000
82.7149999999999	3174.53000000000
83.3079999999999	3174.47000000000
83.9069999999999	3174.40000000000
84.5099999999999	3174.39000000000
85.1079999999999	3174.38000000000
85.7069999999999	3174.30000000000
86.3029999999999	3174.26000000000
86.9009999999999	3174.18000000000
87.5009999999999	3174.13000000000
88.1029999999999	3174.02000000000
88.7089999999999	3174.01000000000
89.3119999999999	3173.99000000000
89.9199999999999	3173.96000000000
90.5189999999999	3173.89000000000
91.1159999999999	3173.87000000000
91.7269999999999	3173.84000000000
92.3259999999999	3173.80000000000
92.9179999999999	3173.78000000000
93.5179999999999	3173.71000000000
94.1169999999999	3173.68000000000
94.7179999999999	3173.64000000000
95.3219999999999	3173.63000000000
95.9219999999999	3173.62000000000
96.5179999999999	3173.56000000000
97.1139999999999	3173.49000000000
97.7119999999999	3173.47000000000
98.3099999999999	3173.45000000000
98.9089999999999	3173.36000000000
99.5129999999999	3173.22000000000
100.133000000000	3173.09000000000
100.739000000000	3173
101.337000000000	3172.97000000000
101.933000000000	3172.96000000000
102.532000000000	3172.91000000000
103.167000000000	3172.75000000000
103.829000000000	3172.71000000000
104.432000000000	3172.68000000000
105.025000000000	3172.65000000000
105.625000000000	3172.64000000000
106.223000000000	3172.62000000000
106.821000000000	3172.60000000000
107.422000000000	3172.60000000000
};
\addplot [line width=0.75, colQRkit]
table {%
4.46400000000000	4529.42000000000
8.62300000000000	4419.14000000000
12.7850000000000	4350.72000000000
17.0370000000000	4294.49000000000
21.1860000000000	4254.80000000000
25.3370000000000	4179.39000000000
29.4680000000000	4121.36000000000
33.6510000000000	4058.62000000000
37.8080000000000	3994.86000000000
41.9480000000000	3927.26000000000
46.3120000000000	3862.23000000000
50.5000000000000	3799.50000000000
54.6480000000000	3740.69000000000
58.8120000000000	3693.51000000000
63.4160000000000	3650.87000000000
67.6960000000000	3610.23000000000
71.8450000000000	3571.28000000000
76.0150000000000	3537.09000000000
80.1550000000000	3504.30000000000
84.3530000000000	3473.41000000000
88.6500000000000	3446.64000000000
92.7910000000000	3421.61000000000
96.9230000000000	3399.25000000000
101.050000000000	3380.32000000000
105.192000000000	3364.30000000000
109.371000000000	3349.90000000000
113.500000000000	3337.11000000000
117.654000000000	3324.50000000000
121.799000000000	3312.46000000000
125.963000000000	3301.78000000000
130.268000000000	3291.37000000000
134.428000000000	3280.83000000000
138.661000000000	3271.19000000000
142.934000000000	3263.29000000000
147.105000000000	3256.07000000000
151.432000000000	3248.70000000000
155.601000000000	3242.90000000000
159.977000000000	3237.67000000000
176.842000000000	3233.59000000000
180.973000000000	3229.70000000000
185.144000000000	3225.58000000000
189.277000000000	3222.45000000000
193.420000000000	3220.11000000000
197.587000000000	3217.87000000000
201.960000000000	3216.04000000000
206.676000000000	3214.49000000000
211.532000000000	3212.37000000000
228.548000000000	3210.95000000000
232.683000000000	3210.73000000000
236.813000000000	3206.17000000000
241.329000000000	3204.10000000000
245.794000000000	3202.22000000000
249.941000000000	3200.86000000000
254.086000000000	3199.44000000000
258.248000000000	3198.01000000000
262.448000000000	3196.77000000000
266.748000000000	3195.70000000000
270.958000000000	3194.67000000000
275.242000000000	3193.84000000000
279.458000000000	3192.65000000000
283.767000000000	3191.49000000000
287.977000000000	3190.31000000000
292.271000000000	3189.12000000000
296.441000000000	3188.36000000000
300.748000000000	3187.74000000000
304.986000000000	3187.01000000000
309.296000000000	3186.28000000000
313.450000000000	3185.69000000000
317.773000000000	3184.66000000000
339.063000000000	3184.04000000000
343.208000000000	3183.55000000000
347.384000000000	3183.06000000000
351.802000000000	3182.71000000000
356.028000000000	3182.31000000000
360.172000000000	3181.97000000000
364.371000000000	3181.52000000000
368.665000000000	3181.13000000000
372.897000000000	3180.71000000000
377.298000000000	3180.11000000000
381.526000000000	3179.62000000000
385.798000000000	3179.28000000000
390.278000000000	3179.09000000000
394.452000000000	3178.94000000000
398.614000000000	3178.77000000000
419.853000000000	3178.63000000000
424.132000000000	3178.54000000000
428.344000000000	3178.48000000000
432.626000000000	3178.41000000000
436.822000000000	3178.27000000000
440.968000000000	3178.07000000000
445.134000000000	3178
449.294000000000	3177.94000000000
453.504000000000	3177.84000000000
457.788000000000	3177.73000000000
461.975000000000	3177.60000000000
466.283000000000	3177.55000000000
470.461000000000	3177.46000000000
474.787000000000	3177.42000000000
479.020000000000	3177.36000000000
483.337000000000	3177.26000000000
487.537000000000	3177.12000000000
492.023000000000	3177.03000000000
496.514000000000	3176.76000000000
500.618000000000	3176.68000000000
504.846000000000	3176.62000000000
509.088000000000	3176.54000000000
513.305000000000	3176.48000000000
517.604000000000	3176.42000000000
521.843000000000	3176.37000000000
526.126000000000	3176.33000000000
530.389000000000	3176.31000000000
534.643000000000	3176.29000000000
538.813000000000	3176.26000000000
543.091000000000	3176.24000000000
547.336000000000	3176.12000000000
551.627000000000	3175.91000000000
555.826000000000	3175.82000000000
560.149000000000	3175.80000000000
564.335000000000	3175.78000000000
568.635000000000	3175.75000000000
572.842000000000	3175.72000000000
577.174000000000	3175.70000000000
581.363000000000	3175.63000000000
585.670000000000	3175.57000000000
589.820000000000	3175.56000000000
594.163000000000	3175.56000000000
598.302000000000	3175.52000000000
602.659000000000	3175.51000000000
636.593000000000	3175.50000000000
640.744000000000	3175.50000000000
644.920000000000	3175.49000000000
649.468000000000	3175.40000000000
653.642000000000	3175.32000000000
658.055000000000	3175.28000000000
662.234000000000	3175.23000000000
679.338000000000	3175.22000000000
683.503000000000	3175.22000000000
687.650000000000	3175.22000000000
691.865000000000	3175.21000000000
696.149000000000	3175.21000000000
700.366000000000	3175.18000000000
704.640000000000	3175.16000000000
708.822000000000	3175.07000000000
713.136000000000	3175.05000000000
717.298000000000	3174.89000000000
721.621000000000	3174.78000000000
725.781000000000	3174.38000000000
730.132000000000	3174.31000000000
734.333000000000	3174.29000000000
738.647000000000	3174.27000000000
742.804000000000	3174.20000000000
747.148000000000	3174.05000000000
751.500000000000	3173.99000000000
755.657000000000	3173.81000000000
760.479000000000	3173.72000000000
764.646000000000	3173.68000000000
768.890000000000	3173.64000000000
773.194000000000	3173.55000000000
777.419000000000	3173.53000000000
781.960000000000	3173.52000000000
786.189000000000	3173.51000000000
790.329000000000	3173.51000000000
794.475000000000	3173.49000000000
798.684000000000	3173.46000000000
802.975000000000	3173.32000000000
807.154000000000	3173.27000000000
815.593000000000	3173.27000000000
828.469000000000	3173.27000000000
845.364000000000	3173.27000000000
849.679000000000	3173.26000000000
875.164000000000	3173.26000000000
879.362000000000	3173.26000000000
883.530000000000	3173.25000000000
913.358000000000	3173.24000000000
917.614000000000	3173.11000000000
921.866000000000	3173.06000000000
926.114000000000	3173.06000000000
930.383000000000	3173.01000000000
934.626000000000	3172.93000000000
938.913000000000	3172.89000000000
951.586000000000	3172.87000000000
955.967000000000	3172.87000000000
960.122000000000	3172.87000000000
964.488000000000	3172.87000000000
968.633000000000	3172.87000000000
973.059000000000	3172.87000000000
977.189000000000	3172.87000000000
981.561000000000	3172.84000000000
985.684000000000	3172.77000000000
990.040000000000	3172.74000000000
994.164000000000	3172.59000000000
998.542000000000	3172.46000000000
1002.72100000000	3172.38000000000
1007.10000000000	3172.23000000000
1011.27000000000	3172.22000000000
1015.69600000000	3172.21000000000
1019.85900000000	3172.13000000000
1024.58600000000	3172.02000000000
1029.17800000000	3171.91000000000
1033.50900000000	3171.90000000000
1041.94100000000	3171.89000000000
1046.27000000000	3171.89000000000
1076.28000000000	3171.83000000000
1080.43900000000	3171.83000000000
1084.69800000000	3171.83000000000
1088.94100000000	3171.80000000000
1101.72600000000	3171.78000000000
1106.08600000000	3171.78000000000
1110.23900000000	3171.78000000000
1114.60500000000	3171.78000000000
1118.89300000000	3171.78000000000
1123.04800000000	3171.78000000000
1127.41700000000	3171.78000000000
1131.72100000000	3171.78000000000
1135.89200000000	3171.78000000000
1140.28600000000	3171.78000000000
1144.45000000000	3171.78000000000
1148.83500000000	3171.77000000000
1152.99600000000	3171.77000000000
1157.37500000000	3171.77000000000
1161.53900000000	3171.77000000000
1165.90800000000	3171.77000000000
};
\addplot [line width=0.75, colQRChol]
table {%
3.19200000000000	4529.42000000000
6.09800000000000	4419.14000000000
8.98700000000000	4350.72000000000
11.8910000000000	4294.49000000000
14.8030000000000	4254.80000000000
17.7190000000000	4179.39000000000
20.6290000000000	4121.36000000000
23.5470000000000	4058.62000000000
26.4550000000000	3994.86000000000
29.3640000000000	3927.27000000000
32.2720000000000	3862.23000000000
35.1780000000000	3799.48000000000
38.0790000000000	3740.67000000000
40.9940000000000	3693.51000000000
43.8970000000000	3650.95000000000
46.8210000000000	3610.57000000000
49.7240000000000	3571.86000000000
55.5400000000000	3536.89000000000
58.4460000000000	3504.05000000000
61.3480000000000	3473.23000000000
67.1810000000000	3446.79000000000
70.0930000000000	3421.37000000000
73	3398.23000000000
75.9140000000000	3378.25000000000
78.8290000000000	3360.93000000000
81.7300000000000	3346.37000000000
84.6280000000000	3333.03000000000
87.5260000000000	3319.88000000000
90.4300000000000	3308.45000000000
93.3400000000000	3298.05000000000
96.2750000000000	3287.91000000000
99.1810000000000	3278.38000000000
102.078000000000	3269.93000000000
104.991000000000	3262.60000000000
107.894000000000	3256.11000000000
110.812000000000	3249.12000000000
119.527000000000	3242.85000000000
122.503000000000	3241.16000000000
125.399000000000	3230.93000000000
131.203000000000	3225.79000000000
134.107000000000	3224.49000000000
137.016000000000	3217.81000000000
139.928000000000	3214.82000000000
142.834000000000	3212.57000000000
145.741000000000	3210.79000000000
148.806000000000	3209.32000000000
152.014000000000	3206.81000000000
154.909000000000	3204.66000000000
157.811000000000	3203.08000000000
160.714000000000	3201.67000000000
163.621000000000	3199.88000000000
166.533000000000	3198.55000000000
169.439000000000	3197.46000000000
172.348000000000	3196.42000000000
175.257000000000	3195.07000000000
181.085000000000	3193.96000000000
184.006000000000	3193.68000000000
186.896000000000	3188.29000000000
189.802000000000	3186.36000000000
192.701000000000	3184.87000000000
195.595000000000	3183.56000000000
198.511000000000	3182.04000000000
201.503000000000	3181.09000000000
207.331000000000	3180.17000000000
210.234000000000	3179.26000000000
213.138000000000	3176.19000000000
216.046000000000	3174.88000000000
218.946000000000	3173.98000000000
221.847000000000	3173.27000000000
224.911000000000	3172.56000000000
228.355000000000	3171.86000000000
231.589000000000	3171.06000000000
238.325000000000	3170.31000000000
241.716000000000	3170.04000000000
245.073000000000	3169.70000000000
248.443000000000	3169.38000000000
252.067000000000	3168.79000000000
255.847000000000	3168.56000000000
259.499000000000	3168.23000000000
262.991000000000	3167.97000000000
266.814000000000	3167.81000000000
270.465000000000	3167.50000000000
274.006000000000	3167.20000000000
279.951000000000	3167.03000000000
282.855000000000	3165.58000000000
285.800000000000	3164.64000000000
288.929000000000	3164.22000000000
291.836000000000	3164.05000000000
294.737000000000	3163.93000000000
297.655000000000	3163.76000000000
300.686000000000	3163.68000000000
303.593000000000	3163.61000000000
315.227000000000	3163.54000000000
318.134000000000	3163.31000000000
321.027000000000	3163.11000000000
323.932000000000	3162.99000000000
326.835000000000	3162.82000000000
329.736000000000	3162.62000000000
332.645000000000	3162.46000000000
335.540000000000	3162.38000000000
338.470000000000	3162.32000000000
341.521000000000	3162.29000000000
344.421000000000	3162.22000000000
347.315000000000	3162.01000000000
350.255000000000	3161.88000000000
353.156000000000	3161.82000000000
356.049000000000	3161.80000000000
358.952000000000	3161.79000000000
361.860000000000	3161.78000000000
364.756000000000	3161.78000000000
367.656000000000	3161.74000000000
370.562000000000	3161.71000000000
373.459000000000	3161.70000000000
376.383000000000	3161.69000000000
379.305000000000	3161.68000000000
382.211000000000	3161.66000000000
385.126000000000	3161.64000000000
388.029000000000	3161.61000000000
390.928000000000	3161.57000000000
393.867000000000	3161.55000000000
396.911000000000	3161.51000000000
399.825000000000	3161.41000000000
402.751000000000	3161.39000000000
405.657000000000	3161.35000000000
408.572000000000	3161.33000000000
411.467000000000	3161.32000000000
414.381000000000	3161.29000000000
417.284000000000	3161.23000000000
420.203000000000	3161.20000000000
426.163000000000	3161.16000000000
449.846000000000	3161.16000000000
452.743000000000	3161.16000000000
455.652000000000	3161.16000000000
458.553000000000	3161.16000000000
461.458000000000	3161.16000000000
464.364000000000	3161.16000000000
473.470000000000	3161.16000000000
476.371000000000	3161.15000000000
479.283000000000	3161.15000000000
482.193000000000	3161.15000000000
485.099000000000	3161.15000000000
488.013000000000	3161.15000000000
490.917000000000	3161.15000000000
493.821000000000	3161.15000000000
496.726000000000	3161.14000000000
499.640000000000	3161.14000000000
502.543000000000	3161.14000000000
505.449000000000	3161.13000000000
508.345000000000	3161.13000000000
511.253000000000	3161.13000000000
514.194000000000	3161.11000000000
517.371000000000	3161.07000000000
520.269000000000	3161.06000000000
523.196000000000	3161
526.235000000000	3160.92000000000
529.139000000000	3160.91000000000
532.044000000000	3160.88000000000
534.955000000000	3160.82000000000
537.869000000000	3160.76000000000
540.778000000000	3160.71000000000
552.445000000000	3160.71000000000
555.345000000000	3160.71000000000
558.245000000000	3160.70000000000
561.179000000000	3160.50000000000
564.220000000000	3160.38000000000
567.132000000000	3160.38000000000
572.946000000000	3160.37000000000
575.863000000000	3160.37000000000
578.760000000000	3160.36000000000
584.578000000000	3160.36000000000
587.477000000000	3160.36000000000
590.381000000000	3160.36000000000
593.305000000000	3160.36000000000
596.371000000000	3160.35000000000
599.269000000000	3160.35000000000
602.170000000000	3160.18000000000
610.888000000000	3160.13000000000
613.796000000000	3160.13000000000
616.703000000000	3160.13000000000
619.597000000000	3160.13000000000
622.515000000000	3160.13000000000
625.409000000000	3160.13000000000
628.351000000000	3160.13000000000
631.426000000000	3160.13000000000
634.327000000000	3160.13000000000
637.236000000000	3160.13000000000
};
\addplot [line width=0.75, colSPQR]
table {%
14.6120000000000	4529.42000000000
29.2030000000000	4419.14000000000
43.7640000000000	4350.72000000000
58.3450000000000	4294.49000000000
72.9200000000000	4254.80000000000
87.5150000000000	4179.39000000000
102.439000000000	4121.36000000000
117.039000000000	4058.62000000000
131.641000000000	3994.86000000000
146.947000000000	3927.26000000000
161.558000000000	3862.23000000000
176.426000000000	3799.50000000000
191.296000000000	3740.69000000000
206.161000000000	3693.51000000000
221.210000000000	3650.87000000000
236.297000000000	3610.02000000000
251.946000000000	3571.16000000000
267.296000000000	3536.64000000000
282.921000000000	3503.15000000000
298.008000000000	3471.66000000000
313.421000000000	3444.36000000000
328.772000000000	3416.68000000000
344.126000000000	3393.72000000000
359.242000000000	3374.23000000000
374.598000000000	3357.09000000000
390.523000000000	3341.80000000000
405.668000000000	3327.45000000000
421.298000000000	3314.18000000000
436.663000000000	3302.89000000000
452.001000000000	3291.72000000000
467.347000000000	3282
482.714000000000	3273.57000000000
498.357000000000	3261.27000000000
513.433000000000	3253.69000000000
528.519000000000	3246.36000000000
543.902000000000	3241.40000000000
559.507000000000	3232.45000000000
574.843000000000	3226.97000000000
622.119000000000	3222.78000000000
652.735000000000	3218.08000000000
668.363000000000	3214.90000000000
684.310000000000	3209.97000000000
699.877000000000	3207.24000000000
714.970000000000	3203.31000000000
730.567000000000	3200.48000000000
745.905000000000	3197.45000000000
760.986000000000	3194.97000000000
776.331000000000	3192.91000000000
791.671000000000	3191.03000000000
807.056000000000	3189.35000000000
822.412000000000	3187.88000000000
837.511000000000	3186.59000000000
852.857000000000	3185.31000000000
868.186000000000	3183.84000000000
883.787000000000	3182.16000000000
899.162000000000	3180.94000000000
975.408000000000	3180
990.688000000000	3179.07000000000
1006.62300000000	3177.49000000000
1022.64600000000	3176.13000000000
1038.73100000000	3174.88000000000
1054.57200000000	3174.06000000000
1070.63800000000	3173.09000000000
1086.57100000000	3171.90000000000
1102.42100000000	3170.66000000000
1118.28200000000	3169.85000000000
1134.14600000000	3169.07000000000
1150.23000000000	3167.73000000000
1166.05800000000	3166.94000000000
1181.81700000000	3166.11000000000
1197.69200000000	3164.89000000000
1213.64000000000	3163.86000000000
1229.59800000000	3163.05000000000
1245.46700000000	3162.26000000000
1261.54000000000	3161.47000000000
1277.51900000000	3160.85000000000
1293.66800000000	3160.24000000000
1309.51700000000	3159.74000000000
1325.37300000000	3159.32000000000
1341.21700000000	3158.93000000000
1357.03600000000	3158.74000000000
1372.85200000000	3158.56000000000
1388.68500000000	3158.48000000000
1404.50400000000	3158.32000000000
1420.33900000000	3158.13000000000
1436.16800000000	3157.91000000000
1452.02200000000	3157.46000000000
1467.82800000000	3157.04000000000
1483.67500000000	3156.83000000000
1499.52900000000	3156.73000000000
1578.79000000000	3156.64000000000
1594.59100000000	3156.47000000000
1610.24300000000	3156.22000000000
1625.94000000000	3156.12000000000
1642.00000000000	3156.02000000000
1657.97900000000	3155.95000000000
1673.84100000000	3155.89000000000
1689.69100000000	3155.81000000000
1705.54400000000	3155.71000000000
1721.39300000000	3155.62000000000
1737.48400000000	3155.54000000000
1753.39300000000	3155.47000000000
1769.23200000000	3155.43000000000
1785.07500000000	3155.40000000000
1800.96100000000	3155.35000000000
1817.04000000000	3155.26000000000
1832.89300000000	3155.19000000000
1849.14500000000	3155.14000000000
1864.76800000000	3155.07000000000
1880.80800000000	3155
1896.73300000000	3154.96000000000
1913.37800000000	3154.93000000000
1929.44000000000	3154.90000000000
1945.26200000000	3154.88000000000
1961.11400000000	3154.86000000000
1976.94300000000	3154.84000000000
1992.84800000000	3154.80000000000
2008.78600000000	3154.76000000000
2024.62400000000	3154.73000000000
2040.47200000000	3154.72000000000
2056.05400000000	3154.71000000000
2072.03300000000	3154.68000000000
2087.85800000000	3154.47000000000
2103.71900000000	3154.31000000000
2119.57100000000	3154.26000000000
2135.42000000000	3154.16000000000
2151.44500000000	3154.13000000000
2167.26900000000	3154.11000000000
2183.16800000000	3154.01000000000
2199.00800000000	3153.97000000000
2214.95000000000	3153.85000000000
2230.76800000000	3153.83000000000
2325.91600000000	3153.82000000000
2435.26300000000	3153.81000000000
2450.64900000000	3153.81000000000
2466.03400000000	3153.81000000000
2481.80200000000	3153.80000000000
2497.43800000000	3153.78000000000
2513.04800000000	3153.69000000000
2528.83000000000	3153.65000000000
2544.49900000000	3153.65000000000
2560.10300000000	3153.61000000000
2575.73900000000	3153.56000000000
2591.37800000000	3153.54000000000
2606.99700000000	3153.30000000000
2622.89800000000	3153.24000000000
2638.53600000000	3153.24000000000
2658.21800000000	3153.23000000000
2674.21400000000	3153.21000000000
2690.75700000000	3153.15000000000
2706.75900000000	3153.04000000000
2722.90200000000	3153
2738.98400000000	3152.99000000000
2754.55600000000	3152.94000000000
2770.77400000000	3152.86000000000
2786.85900000000	3152.81000000000
2803.14700000000	3152.80000000000
2819.22100000000	3152.79000000000
2835.41300000000	3152.66000000000
2851.68700000000	3152.61000000000
2868.27000000000	3152.60000000000
2884.48200000000	3152.50000000000
2901.00700000000	3152.49000000000
2917.45600000000	3152.47000000000
2933.46800000000	3152.46000000000
2949.89400000000	3152.46000000000
2966.35200000000	3152.45000000000
2983.21200000000	3152.45000000000
3005.78900000000	3152.45000000000
3026.98000000000	3152.45000000000
3044.08200000000	3152.40000000000
3060.76700000000	3152.39000000000
3077.28200000000	3152.39000000000
3093.75300000000	3152.38000000000
3110.30700000000	3152.38000000000
3127.06700000000	3152.38000000000
3143.47900000000	3152.38000000000
3159.99600000000	3152.38000000000
3176.28700000000	3152.38000000000
3192.92100000000	3152.38000000000
3209.23600000000	3152.38000000000
3226.02300000000	3152.38000000000
3242.35800000000	3152.38000000000
3258.90400000000	3152.38000000000
3275.22800000000	3152.38000000000
3292.07800000000	3152.38000000000
3308.89000000000	3152.38000000000
3325.17600000000	3152.36000000000
3341.71600000000	3152.32000000000
3358.42400000000	3152.32000000000
3374.89800000000	3152.32000000000
3391.52100000000	3152.32000000000
3408.01800000000	3152.32000000000
3424.51600000000	3152.32000000000
3441.00400000000	3152.32000000000
3457.83900000000	3152.32000000000
3474.45500000000	3152.32000000000
3490.99000000000	3152.32000000000
3507.51900000000	3152.32000000000
3524.14100000000	3152.32000000000
3540.76800000000	3152.32000000000
3557.29500000000	3152.30000000000
3573.77600000000	3152.28000000000
3590.33600000000	3152.28000000000
3606.96700000000	3152.28000000000
3623.51800000000	3152.27000000000
3640.03900000000	3152.27000000000
3657.63500000000	3152.27000000000
3674.14800000000	3152.25000000000
3690.67800000000	3152.22000000000
3707.35300000000	3152.22000000000
3723.85600000000	3152.21000000000
3740.36500000000	3152.21000000000
3756.86100000000	3152.21000000000
3773.38800000000	3152.21000000000
3790.05600000000	3152.21000000000
3806.60700000000	3152.21000000000
3822.82900000000	3152.20000000000
3839.28600000000	3152.18000000000
3855.70800000000	3152.17000000000
3872.22300000000	3152.17000000000
3888.47700000000	3152.17000000000
3904.97400000000	3152.17000000000
3921.52600000000	3152.14000000000
3938.01200000000	3152.13000000000
3954.69900000000	3152.12000000000
3971.35900000000	3152.12000000000
3988.03800000000	3152.12000000000
4004.54700000000	3152.12000000000
4021.09700000000	3152.09000000000
4037.38000000000	3152.04000000000
4053.98200000000	3151.86000000000
4070.61700000000	3151.75000000000
4087.13500000000	3151.69000000000
4103.67000000000	3151.67000000000
4119.97000000000	3151.67000000000
4136.45000000000	3151.67000000000
4152.96800000000	3151.67000000000
4169.64500000000	3151.67000000000
4186.26100000000	3151.67000000000
4203.46700000000	3151.67000000000
4220.83900000000	3151.67000000000
4238.30100000000	3151.67000000000
4255.54800000000	3151.67000000000
4272.03900000000	3151.67000000000
4288.65000000000	3151.67000000000
4304.89600000000	3151.67000000000
4321.55700000000	3151.66000000000
4338.02500000000	3151.66000000000
4354.50200000000	3151.66000000000
4370.75700000000	3151.66000000000
4387.30800000000	3151.66000000000
4403.79200000000	3151.66000000000
4420.33100000000	3151.66000000000
4437.01200000000	3151.66000000000
4453.36000000000	3151.66000000000
4469.68500000000	3151.66000000000
4486.16000000000	3151.66000000000
4502.43700000000	3151.66000000000
4519.35800000000	3151.66000000000
4535.92400000000	3151.66000000000
4552.41000000000	3151.66000000000
4568.91100000000	3151.66000000000
4585.40100000000	3151.66000000000
4601.91700000000	3151.66000000000
4618.36100000000	3151.66000000000
4635.50000000000	3151.66000000000
4652.61600000000	3151.66000000000
4669.07000000000	3151.66000000000
4685.54200000000	3151.66000000000
4702.02400000000	3151.66000000000
4718.53700000000	3151.66000000000
4735.11300000000	3151.65000000000
4751.61100000000	3151.65000000000
4768.14700000000	3151.65000000000
4784.40800000000	3151.65000000000
4800.89800000000	3151.65000000000
4817.41800000000	3151.65000000000
4834.03200000000	3151.65000000000
4850.66800000000	3151.65000000000
4867.16500000000	3151.65000000000
4883.64900000000	3151.65000000000
4900.16700000000	3151.65000000000
4916.78300000000	3151.65000000000
4933.26800000000	3151.65000000000
4949.81700000000	3151.65000000000
4966.20300000000	3151.65000000000
4982.70000000000	3151.65000000000
4999.37200000000	3151.65000000000
5015.85200000000	3151.65000000000
5032.37700000000	3151.65000000000
5048.86600000000	3151.65000000000
5065.36600000000	3151.65000000000
5081.59300000000	3151.65000000000
5098.11500000000	3151.65000000000
5114.77800000000	3151.65000000000
5131.25000000000	3151.65000000000
5147.70600000000	3151.65000000000
5164.24400000000	3151.65000000000
5180.74400000000	3151.65000000000
5197.35000000000	3151.65000000000
5213.79900000000	3151.65000000000
5230.27800000000	3151.65000000000
5246.78100000000	3151.65000000000
5263.25400000000	3151.65000000000
5280.87800000000	3151.65000000000
5298.29500000000	3151.65000000000
5314.78600000000	3151.64000000000
5331.27300000000	3151.64000000000
5347.56200000000	3151.64000000000
5364.23400000000	3151.64000000000
5380.73800000000	3151.64000000000
5397.53700000000	3151.63000000000
5414.07500000000	3151.63000000000
5430.76100000000	3151.63000000000
5447.31600000000	3151.63000000000
5463.89800000000	3151.63000000000
5480.35100000000	3151.63000000000
5496.84000000000	3151.63000000000
5513.56800000000	3151.63000000000
5530.31900000000	3151.63000000000
5546.76400000000	3151.63000000000
5563.25200000000	3151.63000000000
5579.81400000000	3151.63000000000
5596.38900000000	3151.63000000000
5612.85800000000	3151.63000000000
5629.31400000000	3151.63000000000
5645.82400000000	3151.63000000000
5662.31600000000	3151.63000000000
5678.99800000000	3151.63000000000
5695.47600000000	3151.63000000000
5711.99200000000	3151.63000000000
};
\addplot [line width=0.75, colMoreQR]
table {%
6.36800000000000	4529.42000000000
12.1550000000000	4433.03000000000
17.7620000000000	4382.58000000000
23.4640000000000	4340.66000000000
29.2470000000000	4298.34000000000
35.0700000000000	4252.94000000000
40.7840000000000	4221.29000000000
46.4700000000000	4160.80000000000
52.1220000000000	4113.30000000000
57.8680000000000	4052.55000000000
63.6670000000000	3990.26000000000
69.5320000000000	3926
75.3940000000000	3864.72000000000
81.1840000000000	3803.23000000000
86.9390000000000	3744.17000000000
93.5190000000000	3690.07000000000
99.8260000000000	3645.85000000000
105.877000000000	3607.30000000000
111.652000000000	3572.88000000000
117.334000000000	3541.28000000000
123.106000000000	3510.83000000000
128.912000000000	3481.01000000000
134.713000000000	3453.35000000000
159.705000000000	3427.78000000000
165.331000000000	3405.82000000000
171.051000000000	3387.44000000000
176.896000000000	3368.07000000000
182.705000000000	3352.25000000000
188.543000000000	3338.84000000000
194.342000000000	3326.57000000000
200.159000000000	3313.87000000000
205.970000000000	3301.82000000000
211.747000000000	3290.21000000000
217.423000000000	3280.30000000000
223.054000000000	3271.85000000000
228.755000000000	3264.25000000000
234.545000000000	3256.83000000000
240.472000000000	3249.79000000000
246.449000000000	3242.90000000000
252.269000000000	3236.44000000000
258.106000000000	3231.46000000000
263.960000000000	3226.64000000000
269.825000000000	3222.53000000000
275.621000000000	3218.72000000000
281.293000000000	3214.98000000000
286.951000000000	3211.60000000000
305.330000000000	3208.63000000000
311.337000000000	3205.40000000000
318.809000000000	3202.01000000000
325.600000000000	3199.20000000000
332.683000000000	3195.99000000000
339.817000000000	3193.94000000000
347.011000000000	3191.49000000000
353.930000000000	3189.34000000000
361.121000000000	3187.08000000000
368.373000000000	3185.18000000000
374.716000000000	3184.63000000000
380.479000000000	3180.87000000000
386.422000000000	3179.19000000000
392.591000000000	3177.81000000000
398.469000000000	3176.34000000000
404.267000000000	3175.12000000000
410.101000000000	3173.96000000000
415.963000000000	3173.01000000000
421.834000000000	3171.69000000000
427.656000000000	3170.81000000000
433.505000000000	3170.29000000000
439.330000000000	3169.65000000000
445.167000000000	3168.85000000000
450.990000000000	3168.11000000000
456.699000000000	3167.58000000000
462.507000000000	3167.11000000000
468.473000000000	3166.62000000000
474.161000000000	3166.20000000000
480.052000000000	3165.77000000000
485.920000000000	3165.12000000000
491.964000000000	3164.83000000000
497.763000000000	3164.43000000000
503.687000000000	3163.54000000000
509.534000000000	3163.03000000000
515.393000000000	3162.68000000000
521.236000000000	3162.50000000000
527.070000000000	3162.40000000000
532.914000000000	3162.13000000000
538.762000000000	3161.86000000000
544.642000000000	3161.54000000000
550.438000000000	3161.18000000000
556.090000000000	3160.15000000000
561.885000000000	3159.73000000000
567.760000000000	3159.53000000000
573.583000000000	3157.81000000000
579.375000000000	3157.01000000000
585.174000000000	3156.50000000000
590.954000000000	3156.12000000000
596.787000000000	3155.76000000000
602.622000000000	3155.51000000000
608.514000000000	3155.24000000000
614.361000000000	3154.96000000000
620.075000000000	3154.75000000000
625.831000000000	3154.45000000000
631.694000000000	3154.24000000000
637.485000000000	3154.15000000000
643.537000000000	3154.05000000000
649.412000000000	3153.87000000000
655.246000000000	3153.69000000000
661.199000000000	3153.55000000000
667.025000000000	3153.36000000000
672.892000000000	3153.23000000000
678.728000000000	3153.10000000000
684.628000000000	3152.91000000000
691.492000000000	3152.66000000000
697.401000000000	3152.33000000000
703.171000000000	3152.23000000000
709.025000000000	3152.17000000000
714.771000000000	3151.99000000000
720.549000000000	3151.84000000000
726.301000000000	3151.69000000000
732.074000000000	3151.52000000000
737.880000000000	3151.40000000000
743.661000000000	3151.22000000000
749.424000000000	3150.75000000000
755.170000000000	3150.45000000000
760.923000000000	3150.30000000000
766.778000000000	3149.95000000000
772.541000000000	3149.54000000000
778.325000000000	3149.28000000000
784.126000000000	3149.14000000000
789.920000000000	3148.94000000000
795.712000000000	3148.86000000000
801.494000000000	3148.64000000000
807.222000000000	3148.38000000000
812.980000000000	3148.02000000000
818.741000000000	3147.71000000000
824.566000000000	3147.51000000000
830.347000000000	3147.40000000000
836.153000000000	3147.31000000000
842.018000000000	3147.22000000000
847.835000000000	3147.10000000000
854.203000000000	3147.06000000000
860.419000000000	3146.83000000000
866.561000000000	3146.56000000000
872.776000000000	3146.28000000000
879.043000000000	3146.14000000000
885.348000000000	3145.89000000000
891.144000000000	3145.79000000000
896.930000000000	3145.74000000000
902.716000000000	3145.62000000000
908.587000000000	3145.50000000000
914.352000000000	3145.36000000000
920.165000000000	3145.08000000000
925.923000000000	3144.93000000000
931.756000000000	3144.78000000000
937.535000000000	3144.58000000000
943.295000000000	3144.38000000000
949.090000000000	3144.10000000000
954.877000000000	3143.90000000000
961.191000000000	3143.74000000000
966.970000000000	3143.55000000000
985.462000000000	3143.35000000000
991.372000000000	3143.05000000000
997.235000000000	3142.70000000000
1003.09000000000	3142.47000000000
1008.91800000000	3142.24000000000
1014.72400000000	3142.05000000000
1020.61100000000	3141.95000000000
1026.47800000000	3141.87000000000
1032.40500000000	3141.70000000000
1038.25100000000	3141.60000000000
1044.37800000000	3141.48000000000
1050.21600000000	3141.34000000000
1056.06800000000	3141.13000000000
1061.92700000000	3141
1067.78300000000	3140.84000000000
1073.62100000000	3140.62000000000
1079.46500000000	3140.37000000000
1085.22000000000	3140.27000000000
1091.01800000000	3140.13000000000
1096.80000000000	3140.05000000000
1102.69200000000	3139.93000000000
1108.52200000000	3139.75000000000
1114.33800000000	3139.54000000000
1120.15800000000	3139.40000000000
1125.93500000000	3139.25000000000
1131.73500000000	3139.14000000000
1137.59000000000	3138.92000000000
1143.41400000000	3138.83000000000
1149.25000000000	3138.80000000000
1155.08500000000	3138.78000000000
1160.90300000000	3138.77000000000
1166.69200000000	3138.73000000000
1172.41100000000	3138.61000000000
1178.27800000000	3138.49000000000
1184.09200000000	3138.31000000000
1189.90600000000	3138.23000000000
1195.78700000000	3138.17000000000
1201.64600000000	3137.94000000000
1207.84700000000	3137.83000000000
1214.08600000000	3137.75000000000
1220.18400000000	3137.69000000000
1225.96300000000	3137.61000000000
1231.77600000000	3137.52000000000
1237.60700000000	3137.30000000000
1243.44700000000	3137.10000000000
1249.29300000000	3136.95000000000
1255.13000000000	3136.79000000000
1261.28500000000	3136.74000000000
1267.12300000000	3136.73000000000
1272.98900000000	3136.67000000000
1278.85100000000	3136.39000000000
1285.30600000000	3136.25000000000
1291.24300000000	3136.18000000000
1297.20800000000	3136.17000000000
1302.98200000000	3136.16000000000
1308.79600000000	3136.14000000000
1314.60600000000	3136.13000000000
1320.38500000000	3136.10000000000
1326.24000000000	3136.05000000000
1332.07000000000	3135.90000000000
1337.85000000000	3135.74000000000
1343.63700000000	3135.71000000000
1349.39300000000	3135.64000000000
1368.01900000000	3135.58000000000
1373.88900000000	3135.47000000000
1379.65900000000	3135.29000000000
1385.50100000000	3134.74000000000
1391.29800000000	3134.63000000000
1397.08500000000	3134.59000000000
1402.94800000000	3134.53000000000
1408.66200000000	3134.38000000000
1414.42700000000	3134.29000000000
1420.19200000000	3134.26000000000
1439.17200000000	3134.20000000000
1445.14700000000	3134.14000000000
1450.90900000000	3133.97000000000
1456.65900000000	3133.91000000000
1462.34900000000	3133.79000000000
1468.16900000000	3133.66000000000
1473.88400000000	3133.62000000000
1479.62500000000	3133.61000000000
1485.50200000000	3133.56000000000
1491.26600000000	3133.52000000000
1497.03400000000	3133.45000000000
1503.01200000000	3133.28000000000
1508.75800000000	3133.19000000000
1514.49700000000	3133.14000000000
1520.24200000000	3133.06000000000
1526.06300000000	3133
1531.81100000000	3132.96000000000
1537.79500000000	3132.83000000000
1543.54200000000	3132.64000000000
1549.41800000000	3132.60000000000
1555.32300000000	3132.58000000000
1561.24600000000	3132.56000000000
1567.04000000000	3132.53000000000
1572.94700000000	3132.50000000000
1578.71600000000	3132.39000000000
1585.10700000000	3132.34000000000
1590.87200000000	3132.23000000000
1596.66800000000	3131.98000000000
1602.46700000000	3131.69000000000
1608.23200000000	3131.50000000000
1613.96700000000	3131.32000000000
1619.78900000000	3131.29000000000
1625.67000000000	3131.26000000000
1631.42500000000	3131.24000000000
1637.22700000000	3131.17000000000
1643.05700000000	3131.05000000000
1648.81500000000	3131.02000000000
1654.59100000000	3130.91000000000
1660.35800000000	3130.85000000000
1666.32300000000	3130.83000000000
1672.09400000000	3130.81000000000
1677.86600000000	3130.78000000000
1683.63700000000	3130.78000000000
1689.46400000000	3130.78000000000
1695.17000000000	3130.78000000000
1700.92800000000	3130.70000000000
1706.72000000000	3130.62000000000
1712.57500000000	3130.46000000000
1718.49500000000	3130.42000000000
1724.26400000000	3130.40000000000
1730.01800000000	3130.36000000000
1735.76900000000	3130.33000000000
1741.51600000000	3130.32000000000
1747.35700000000	3130.26000000000
1753.11700000000	3130.16000000000
1758.85400000000	3129.95000000000
1764.60000000000	3129.89000000000
1770.47000000000	3129.88000000000
1776.35400000000	3129.82000000000
1782.10900000000	3129.64000000000
1787.87200000000	3129.45000000000
1793.61300000000	3129.36000000000
1799.36000000000	3129.34000000000
1805.17300000000	3129.34000000000
1810.96400000000	3129.30000000000
1816.74400000000	3129.24000000000
1822.69000000000	3129.21000000000
1828.48000000000	3129.21000000000
1834.23200000000	3129.20000000000
1839.98700000000	3129.13000000000
1845.75300000000	3129.13000000000
1851.52000000000	3129.11000000000
1857.30000000000	3129.04000000000
1863.11800000000	3129.01000000000
1868.86800000000	3129
1874.59900000000	3128.96000000000
1880.81600000000	3128.85000000000
1886.57300000000	3128.81000000000
1892.49600000000	3128.78000000000
1898.16400000000	3128.74000000000
1903.83500000000	3128.71000000000
1909.47100000000	3128.70000000000
1915.11300000000	3128.59000000000
1920.63300000000	3128.50000000000
1926.22500000000	3128.33000000000
1931.80300000000	3128.26000000000
1937.40200000000	3128.20000000000
1942.92700000000	3128.10000000000
1973.28700000000	3128.07000000000
1978.82700000000	3128.07000000000
1984.43200000000	3128.07000000000
1989.96000000000	3128.07000000000
};
\end{axis}

\end{tikzpicture}

%% file: dubrovnik_iter_float.tex
%
%
%
\begin{tikzpicture}

\pgfplotsset{compat=newest} 

\definecolor{color0}{rgb}{1, 1, 1}

\tikzstyle{every node}=[font=\footnotesize]

\definecolor{colCholesky}{RGB}{0, 255, 0}
\definecolor{colQRkit}{RGB}{255, 0, 0}
\definecolor{colQRChol}{RGB}{0, 0, 255}
\definecolor{colSPQR}{RGB}{0, 255, 255}
\definecolor{colSSBA}{RGB}{255, 255, 0}
\definecolor{colMoreQR}{RGB}{200, 0, 255}

\begin{axis}[
xlabel={Iterations},
ylabel={Energy},
xmin=0, xmax=150,
ymin=3000, ymax=4500,
width=\figurewidth,
height=\figureheight,
at={(0\figurewidth,0\figureheight)},
xmajorgrids,
x grid style={lightgray},
ymajorgrids,
y grid style={lightgray},
axis line style={black},
axis background/.style={fill=color0},
legend style={at={(0.97,0.97)}, anchor=north east},
legend cell align={left},
legend entries={{Cholesky},{QRkit},{QRkit + Cholesky},{SSBA},{Mor{\'e} QRkit}}
]
%
%
\addplot [line width=1, colCholesky]
table {%
1	4529.42000000000
2	4423.65000000000
3	4353.31000000000
4	4296.09000000000
5	4246.01000000000
6	4187.86000000000
7	4128.22000000000
8	4073.12000000000
9	4018.39000000000
10	3958.82000000000
11	3900
12	3845.31000000000
13	3791.09000000000
14	3742.79000000000
15	3702.40000000000
16	3665.64000000000
17	3630.97000000000
18	3597.41000000000
19	3564.79000000000
20	3534.18000000000
21	3506.04000000000
22	3481.03000000000
23	3457.08000000000
24	3434.62000000000
25	3414.85000000000
26	3398.25000000000
27	3383.88000000000
28	3369.72000000000
29	3356.28000000000
30	3343.43000000000
31	3331.33000000000
32	3321.77000000000
33	3313.84000000000
34	3307.30000000000
35	3301.58000000000
36	3296.39000000000
37	3291.98000000000
38	3287.36000000000
39	3282.41000000000
40	3278.26000000000
41	3273.83000000000
42	3269.77000000000
43	3266.90000000000
44	3264.03000000000
45	3260.57000000000
46	3257.92000000000
47	3255.24000000000
48	3253.56000000000
49	3251.11000000000
50	3249.43000000000
51	3247.97000000000
52	3245.70000000000
53	3244.21000000000
54	3242.46000000000
55	3240.97000000000
56	3238.99000000000
57	3237.83000000000
58	3236.41000000000
59	3235.18000000000
60	3233.80000000000
61	3232.64000000000
62	3231.44000000000
63	3230.34000000000
64	3228.76000000000
65	3227.95000000000
66	3226.76000000000
67	3225.81000000000
68	3224.80000000000
69	3223.91000000000
70	3223.02000000000
71	3221.77000000000
72	3221.26000000000
73	3220.18000000000
74	3219.80000000000
75	3219.13000000000
76	3218.87000000000
77	3217.85000000000
78	3217.34000000000
79	3217.05000000000
80	3216.42000000000
81	3215.87000000000
82	3215.24000000000
83	3214.45000000000
84	3213.74000000000
85	3213.44000000000
86	3213.15000000000
87	3212.73000000000
88	3212.45000000000
89	3212.18000000000
90	3211.58000000000
91	3211.32000000000
92	3211
93	3210.80000000000
94	3209.92000000000
95	3209.85000000000
96	3209.66000000000
97	3209.37000000000
98	3209.20000000000
99	3208.42000000000
100	3207.95000000000
101	3207.56000000000
102	3207.04000000000
103	3206.60000000000
104	3206.22000000000
105	3205.84000000000
106	3205.56000000000
107	3205.32000000000
108	3205.05000000000
109	3204.70000000000
110	3204.27000000000
111	3204.20000000000
112	3203.79000000000
113	3203.53000000000
114	3203.39000000000
115	3203.32000000000
116	3203.10000000000
117	3202.85000000000
118	3202.58000000000
119	3202.50000000000
120	3202.47000000000
121	3202.40000000000
122	3202.26000000000
123	3201.93000000000
124	3201.77000000000
125	3201.76000000000
126	3201.47000000000
127	3201.36000000000
128	3201.25000000000
129	3201.09000000000
130	3200.94000000000
131	3200.88000000000
132	3200.76000000000
133	3200.68000000000
134	3200.68000000000
135	3200.33000000000
136	3200.29000000000
137	3200.02000000000
138	3199.99000000000
139	3199.93000000000
140	3199.86000000000
141	3199.77000000000
142	3199.48000000000
143	3199.33000000000
144	3199.27000000000
145	3199.11000000000
146	3198.99000000000
147	3198.87000000000
148	3198.69000000000
149	3198.68000000000
150	3198.37000000000
151	3198.12000000000
152	3197.93000000000
153	3197.78000000000
154	3197.72000000000
155	3197.48000000000
156	3197.44000000000
157	3197.30000000000
158	3197.22000000000
159	3197.17000000000
160	3197.11000000000
161	3197.06000000000
162	3197.04000000000
163	3196.96000000000
164	3196.84000000000
165	3196.63000000000
166	3196.55000000000
167	3196.41000000000
168	3196.20000000000
169	3196.10000000000
170	3195.92000000000
171	3195.84000000000
172	3195.77000000000
173	3195.67000000000
174	3195.60000000000
175	3195.52000000000
176	3195.44000000000
177	3195.30000000000
178	3195.23000000000
179	3195.10000000000
180	3194.99000000000
181	3194.97000000000
182	3194.83000000000
183	3194.82000000000
184	3194.79000000000
185	3194.61000000000
186	3194.49000000000
187	3194.16000000000
188	3194.13000000000
189	3194.05000000000
190	3193.80000000000
191	3193.48000000000
192	3193.34000000000
193	3193.20000000000
194	3193.02000000000
195	3193.01000000000
196	3192.85000000000
197	3192.83000000000
198	3192.73000000000
199	3192.68000000000
200	3192.54000000000
201	3192.53000000000
202	3192.31000000000
203	3192.20000000000
204	3192.14000000000
205	3191.82000000000
206	3191.32000000000
207	3190.94000000000
208	3190.44000000000
209	3190.23000000000
210	3190.09000000000
211	3189.82000000000
212	3189.47000000000
213	3189.27000000000
214	3188.91000000000
215	3188.56000000000
216	3188.30000000000
217	3187.84000000000
218	3187.72000000000
219	3187.51000000000
220	3187.22000000000
221	3186.89000000000
222	3186.81000000000
223	3185.96000000000
224	3185.82000000000
225	3185.54000000000
226	3185.45000000000
227	3185.32000000000
228	3185.17000000000
229	3185.06000000000
230	3185.06000000000
231	3184.56000000000
232	3184.47000000000
233	3184.34000000000
234	3184.23000000000
235	3184.20000000000
236	3184.18000000000
237	3184.17000000000
238	3184.11000000000
239	3184.03000000000
240	3183.98000000000
241	3183.73000000000
242	3183.70000000000
243	3183.40000000000
244	3183.28000000000
245	3183.22000000000
246	3182.98000000000
247	3182.93000000000
248	3182.88000000000
249	3182.78000000000
250	3182.73000000000
251	3182.59000000000
252	3182.47000000000
253	3182.33000000000
254	3182.18000000000
255	3182.13000000000
256	3181.98000000000
257	3181.92000000000
258	3181.83000000000
259	3181.79000000000
260	3181.66000000000
261	3181.62000000000
262	3181.54000000000
263	3181.35000000000
264	3181.29000000000
265	3181.28000000000
266	3181.23000000000
267	3181.16000000000
268	3181.07000000000
269	3181.01000000000
270	3180.95000000000
271	3180.92000000000
272	3180.90000000000
273	3180.74000000000
274	3180.72000000000
275	3180.65000000000
276	3180.44000000000
277	3180.42000000000
278	3180.33000000000
279	3180.16000000000
280	3179.97000000000
281	3179.87000000000
282	3179.73000000000
283	3179.68000000000
284	3179.62000000000
285	3179.52000000000
286	3179.46000000000
287	3179.38000000000
288	3179.34000000000
289	3179.30000000000
290	3179.25000000000
291	3179.22000000000
292	3179.02000000000
293	3178.85000000000
294	3178.73000000000
295	3178.56000000000
296	3178.42000000000
297	3178.33000000000
298	3178.32000000000
299	3178.27000000000
300	3178.24000000000
301	3178.16000000000
302	3178.16000000000
303	3178.12000000000
304	3178.07000000000
305	3178
306	3177.90000000000
307	3177.82000000000
308	3177.82000000000
309	3177.72000000000
310	3177.69000000000
311	3177.61000000000
312	3177.58000000000
313	3177.56000000000
314	3177.53000000000
315	3177.53000000000
316	3177.51000000000
317	3177.45000000000
318	3177.43000000000
319	3177.36000000000
320	3177.34000000000
321	3177.31000000000
322	3177.26000000000
323	3177.23000000000
324	3177.23000000000
325	3177.19000000000
326	3177
327	3176.96000000000
328	3176.86000000000
329	3176.83000000000
330	3176.78000000000
331	3176.69000000000
332	3176.60000000000
333	3176.57000000000
334	3176.34000000000
335	3176.25000000000
336	3176.19000000000
337	3176.02000000000
338	3176.01000000000
339	3175.98000000000
340	3175.80000000000
341	3175.64000000000
342	3175.56000000000
343	3175.49000000000
344	3175.47000000000
345	3175.41000000000
346	3175.30000000000
347	3175.14000000000
348	3175.08000000000
349	3174.99000000000
350	3174.83000000000
351	3174.58000000000
352	3174.42000000000
353	3174.15000000000
354	3174.08000000000
355	3173.97000000000
356	3173.79000000000
357	3173.75000000000
358	3173.67000000000
359	3173.67000000000
360	3173.66000000000
361	3173.62000000000
362	3173.53000000000
363	3173.50000000000
364	3173.49000000000
365	3173.44000000000
366	3173.32000000000
367	3173.22000000000
368	3172.94000000000
369	3172.88000000000
370	3172.85000000000
371	3172.71000000000
372	3172.35000000000
373	3172.14000000000
374	3171.90000000000
375	3171.89000000000
376	3171.79000000000
377	3171.74000000000
378	3171.68000000000
379	3171.37000000000
380	3171.08000000000
381	3171.06000000000
382	3170.95000000000
383	3170.90000000000
384	3170.90000000000
385	3170.74000000000
386	3170.55000000000
387	3170.45000000000
388	3170.43000000000
389	3170.38000000000
390	3170.36000000000
391	3170.34000000000
392	3170.24000000000
393	3170.08000000000
394	3170.01000000000
395	3169.90000000000
396	3169.73000000000
397	3169.72000000000
398	3169.60000000000
399	3169.54000000000
400	3169.51000000000
401	3169.45000000000
402	3169.36000000000
403	3169.30000000000
404	3169.16000000000
405	3169.10000000000
406	3168.90000000000
407	3168.81000000000
408	3168.78000000000
409	3168.40000000000
410	3168.31000000000
411	3168.27000000000
412	3168.25000000000
413	3168.17000000000
414	3168.15000000000
415	3168.04000000000
416	3167.98000000000
417	3167.92000000000
418	3167.84000000000
419	3167.76000000000
420	3167.74000000000
421	3167.66000000000
422	3167.60000000000
423	3167.54000000000
424	3167.50000000000
425	3167.46000000000
426	3167.39000000000
427	3167.37000000000
428	3167.32000000000
429	3167.25000000000
430	3167.24000000000
431	3167.23000000000
432	3167.21000000000
433	3167.14000000000
434	3167.12000000000
435	3167.09000000000
436	3167.09000000000
437	3166.91000000000
438	3166.89000000000
439	3166.89000000000
440	3166.74000000000
441	3166.63000000000
442	3166.62000000000
443	3166.51000000000
444	3166.51000000000
};
\addplot [line width=1, colQRkit]
table {%
 1	4529.42000000000
2	4419.14000000000
3	4350.64000000000
4	4293.80000000000
5	4253.31000000000
6	4177.51000000000
7	4117.36000000000
8	4058.96000000000
9	4005.22000000000
10	3932.07000000000
11	3864.74000000000
12	3804.32000000000
13	3755.61000000000
14	3695.38000000000
15	3658.61000000000
16	3602.43000000000
17	3566.71000000000
18	3524.96000000000
19	3499.71000000000
20	3476.81000000000
21	3450.08000000000
22	3408.74000000000
23	3379.27000000000
24	3376.16000000000
25	3370.47000000000
26	3360.13000000000
27	3311.87000000000
28	3308.25000000000
29	3285.89000000000
30	3257.81000000000
31	3238.19000000000
32	3237.89000000000
33	3214.92000000000
34	3211.48000000000
35	3197.38000000000
36	3191.24000000000
37	3190.35000000000
38	3180.58000000000
39	3172.32000000000
40	3170.44000000000
41	3164.61000000000
42	3160.53000000000
43	3156.71000000000
44	3154.08000000000
45	3151.25000000000
46	3150.60000000000
47	3149.42000000000
48	3147.50000000000
49	3144.71000000000
50	3142.92000000000
51	3140.12000000000
52	3138.44000000000
53	3137.42000000000
54	3135.46000000000
55	3134.45000000000
56	3132.90000000000
57	3132.23000000000
58	3131.15000000000
59	3130.23000000000
60	3129.64000000000
61	3129
62	3128.81000000000
63	3128.14000000000
64	3128
65	3126.71000000000
66	3126.03000000000
67	3125.12000000000
68	3124.66000000000
69	3124.21000000000
70	3123.76000000000
71	3123.55000000000
72	3123.12000000000
73	3121.93000000000
74	3121.55000000000
75	3121.30000000000
76	3120.28000000000
77	3119.71000000000
78	3119.36000000000
79	3118.09000000000
80	3117.49000000000
81	3117.17000000000
82	3116.13000000000
83	3115.56000000000
84	3114.73000000000
85	3114.42000000000
86	3114.02000000000
87	3113.52000000000
88	3112.96000000000
89	3112.80000000000
90	3112.77000000000
91	3112.02000000000
92	3111.63000000000
93	3111.35000000000
94	3110.94000000000
95	3110.50000000000
96	3110.20000000000
97	3110.08000000000
98	3110.02000000000
99	3109.80000000000
100	3109.51000000000
101	3109.30000000000
102	3109.26000000000
103	3108.93000000000
104	3108.72000000000
105	3108.46000000000
106	3108.27000000000
107	3108.10000000000
108	3107.93000000000
109	3107.79000000000
110	3107.14000000000
111	3106.67000000000
112	3106.17000000000
113	3106.12000000000
114	3104.95000000000
115	3104.95000000000
116	3103.49000000000
117	3102.97000000000
118	3102.84000000000
119	3102.35000000000
120	3102.30000000000
121	3101.33000000000
122	3100.75000000000
123	3100.61000000000
124	3100.21000000000
125	3099.74000000000
126	3099.21000000000
127	3098.91000000000
128	3098.87000000000
129	3098.36000000000
130	3098.15000000000
131	3097.73000000000
132	3097.33000000000
133	3095.98000000000
134	3095.01000000000
135	3094.69000000000
136	3094.46000000000
137	3094.25000000000
138	3094.22000000000
139	3094.11000000000
140	3094.02000000000
141	3093.99000000000
142	3093.85000000000
143	3093.74000000000
144	3093.68000000000
145	3093.66000000000
146	3093.65000000000
147	3093.62000000000
148	3093.47000000000
149	3093.36000000000
150	3093.28000000000
151	3093.22000000000
152	3093.16000000000
153	3093.14000000000
154	3093.10000000000
155	3093.10000000000
156	3093.03000000000
157	3092.92000000000
};
\addplot [line width=1, colQRChol]
table {%
1	4529.42000000000
2	4418.78000000000
3	4347.03000000000
4	4318.19000000000
5	4239.33000000000
6	4221.79000000000
7	4132.50000000000
8	4071.07000000000
9	4003.54000000000
10	3938.32000000000
11	3870.99000000000
12	3811.52000000000
13	3754.27000000000
14	3702.06000000000
15	3656.05000000000
16	3613.08000000000
17	3576.54000000000
18	3540.73000000000
19	3508.25000000000
};
\addplot [line width=1, colSSBA]
table {%
1	4529.41699200000
2	4413.05419900000
3	4333.51220700000
4	4328.76464800000
5	4216.72802700000
6	4135.75830100000
7	4076.15527300000
8	4019.31396500000
9	3962.88037100000
10	3909.35205100000
11	3860.09912100000
12	3810.58569300000
13	3763.75708000000
14	3720.35131800000
15	3680.22583000000
16	3643.53051800000
17	3611.99902300000
18	3580.79834000000
19	3549.94677700000
20	3521.74731400000
21	3495.52368200000
22	3472.06054700000
23	3450.72827100000
24	3432.91601600000
25	3416.79589800000
26	3401.24365200000
27	3386.07226600000
28	3372.63061500000
29	3360.59130900000
30	3349.72143600000
31	3340.07080100000
32	3330.83252000000
33	3322.12622100000
34	3314.19287100000
35	3307.11377000000
36	3300.83642600000
37	3295.13916000000
38	3289.95581100000
39	3284.87719700000
40	3280.51367200000
41	3276.38525400000
42	3272.89135700000
43	3269.11303700000
44	3266.02685500000
45	3263.18725600000
46	3260.17089800000
47	3257.54907200000
48	3255.12475600000
49	3252.41503900000
50	3250.64062500000
51	3249.25219700000
};
\addplot [line width=1, colMoreQR]
table {%
1	4529.42000000000
2	4433.02000000000
3	4382.56000000000
4	4340.65000000000
5	4298.48000000000
6	4253.23000000000
7	4221.80000000000
8	4161.27000000000
9	4112.85000000000
10	4051.34000000000
11	4000.09000000000
12	3992.34000000000
13	3902.04000000000
14	3859.78000000000
15	3741.86000000000
16	3686.17000000000
17	3684.54000000000
18	3577.31000000000
19	3540.69000000000
20	3513.16000000000
21	3473.19000000000
22	3434.77000000000
23	3406.78000000000
24	3379.05000000000
25	3378.47000000000
26	3337.58000000000
27	3327.01000000000
28	3306.22000000000
29	3288.70000000000
30	3268.06000000000
31	3254.61000000000
32	3251.10000000000
33	3232.87000000000
34	3218.17000000000
35	3213.37000000000
36	3206.54000000000
37	3193.88000000000
38	3186.03000000000
39	3176.94000000000
40	3171.26000000000
41	3165.90000000000
42	3162.64000000000
43	3157.38000000000
44	3154.44000000000
45	3151.13000000000
46	3148.03000000000
47	3144.46000000000
48	3142.78000000000
49	3140.29000000000
50	3138.85000000000
51	3134.89000000000
52	3134.02000000000
53	3128.59000000000
54	3126.13000000000
55	3125.11000000000
56	3123.62000000000
57	3122.62000000000
58	3122.13000000000
59	3120.86000000000
60	3120.83000000000
61	3118.65000000000
62	3117.59000000000
63	3117.07000000000
64	3116.87000000000
65	3116.71000000000
66	3115.26000000000
67	3114.23000000000
68	3112.39000000000
69	3111.61000000000
70	3110.73000000000
71	3110.23000000000
72	3108.80000000000
73	3108.26000000000
74	3107.86000000000
75	3107.46000000000
76	3107.05000000000
77	3106.68000000000
78	3106.45000000000
79	3106.37000000000
80	3105.83000000000
81	3105.15000000000
82	3104.52000000000
83	3104.09000000000
84	3104.04000000000
85	3103.82000000000
86	3103.71000000000
87	3103.09000000000
88	3103.08000000000
89	3103.01000000000
90	3102.67000000000
91	3102.52000000000
92	3102.32000000000
93	3102.18000000000
94	3101.97000000000
95	3101.81000000000
96	3101.66000000000
97	3101.58000000000
98	3101.46000000000
99	3101.39000000000
100	3101.21000000000
101	3101.06000000000
102	3100.90000000000
103	3100.64000000000
104	3100.45000000000
105	3100.24000000000
106	3100.01000000000
107	3099.86000000000
108	3099.80000000000
109	3099.79000000000
110	3099.65000000000
111	3099.59000000000
112	3099.53000000000
113	3099.50000000000
114	3099.46000000000
115	3099.32000000000
116	3099.25000000000
117	3099.23000000000
118	3099.18000000000
119	3099.18000000000
120	3099.14000000000
121	3098.99000000000
122	3098.94000000000
123	3098.92000000000
124	3098.90000000000
125	3098.77000000000
126	3098.75000000000
127	3098.74000000000
128	3098.67000000000
129	3098.63000000000
130	3098.60000000000
131	3098.49000000000
132	3098.49000000000
133	3098.42000000000
134	3098.38000000000
135	3098.37000000000
136	3098.32000000000
137	3098.28000000000
138	3098.25000000000
139	3098.18000000000
140	3098.16000000000
141	3098.11000000000
142	3098.08000000000
143	3098.04000000000
144	3098.02000000000
145	3098.01000000000
146	3098
147	3097.96000000000
148	3097.93000000000
149	3097.91000000000
150	3097.83000000000
151	3097.79000000000
152	3097.77000000000
153	3097.76000000000
154	3097.76000000000
155	3097.76000000000
156	3097.76000000000
157	3097.76000000000
158	3097.72000000000
159	3097.71000000000
160	3097.65000000000
161	3097.64000000000
162	3097.62000000000
163	3097.61000000000
164	3097.61000000000
165	3097.58000000000
166	3097.39000000000
167	3097.38000000000
168	3097.35000000000
169	3097.34000000000
170	3097.30000000000
171	3097.29000000000
172	3097.29000000000
173	3097.24000000000
174	3097.20000000000
175	3097.17000000000
176	3097.12000000000
177	3097.11000000000
178	3097.09000000000
179	3097.07000000000
180	3096.98000000000
181	3096.97000000000
182	3096.77000000000
183	3096.76000000000
184	3096.51000000000
185	3096.49000000000
186	3096.45000000000
187	3096.38000000000
188	3096.35000000000
189	3096.33000000000
190	3096.32000000000
191	3096.20000000000
192	3096.01000000000
193	3096
194	3095.96000000000
195	3095.95000000000
196	3095.93000000000
197	3095.92000000000
198	3095.79000000000
199	3095.73000000000
200	3095.71000000000
201	3095.71000000000
202	3095.70000000000
203	3095.67000000000
204	3095.62000000000
205	3095.62000000000
206	3095.60000000000
207	3095.56000000000
208	3095.44000000000
209	3095.36000000000
210	3095.29000000000
211	3095.17000000000
212	3095.01000000000
213	3094.97000000000
214	3094.94000000000
215	3094.91000000000
216	3094.90000000000
217	3094.78000000000
218	3094.77000000000
219	3094.76000000000
220	3094.74000000000
221	3094.73000000000
222	3094.73000000000
223	3094.73000000000
224	3094.73000000000
225	3094.73000000000
};
\end{axis}

\end{tikzpicture}

%% file: dubrovnik_float.tex
%
%
%
\begin{tikzpicture}

\pgfplotsset{compat=newest} 

\definecolor{color0}{rgb}{1, 1, 1}

\tikzstyle{every node}=[font=\footnotesize]

\definecolor{colCholesky}{RGB}{0, 255, 0}
\definecolor{colQRkit}{RGB}{255, 0, 0}
\definecolor{colQRChol}{RGB}{0, 0, 255}
\definecolor{colSPQR}{RGB}{0, 255, 255}
\definecolor{colSSBA}{RGB}{255, 255, 0}
\definecolor{colMoreQR}{RGB}{200, 0, 255}

\begin{axis}[
xmode=log,
log ticks with fixed point,
xtick={1,5,10,20,50,120,300,1000},
xlabel={Time [s]},
ylabel={Energy},
xmin=0, xmax=5e3,
ymin=3000, ymax=4500,
width=\figurewidth,
height=\figureheight,
at={(0\figurewidth,0\figureheight)},
xmajorgrids,
x grid style={lightgray},
ymajorgrids,
y grid style={lightgray},
axis line style={black},
axis background/.style={fill=color0},
]
%
%
\addplot [line width=1, colCholesky]
table {%
0.482000000000000	4529.42000000000
0.959000000000000	4423.65000000000
1.43400000000000	4353.31000000000
1.90600000000000	4296.09000000000
2.38000000000000	4246.01000000000
3.15400000000000	4187.86000000000
3.93200000000000	4128.22000000000
4.71000000000000	4073.12000000000
5.80200000000000	4018.39000000000
6.57200000000000	3958.82000000000
7.65500000000000	3900
8.43800000000000	3845.31000000000
10.3100000000000	3791.09000000000
10.7850000000000	3742.79000000000
11.2570000000000	3702.40000000000
11.7270000000000	3665.64000000000
12.5130000000000	3630.97000000000
13.5910000000000	3597.41000000000
14.3690000000000	3564.79000000000
14.8460000000000	3534.18000000000
15.9340000000000	3506.04000000000
17.0200000000000	3481.03000000000
17.8060000000000	3457.08000000000
18.8950000000000	3434.62000000000
19.6840000000000	3414.85000000000
20.4700000000000	3398.25000000000
21.2470000000000	3383.88000000000
22.3330000000000	3369.72000000000
23.4200000000000	3356.28000000000
24.1990000000000	3343.43000000000
25.2880000000000	3331.33000000000
26.0700000000000	3321.77000000000
27.1640000000000	3313.84000000000
27.6540000000000	3307.30000000000
28.8730000000000	3301.58000000000
29.3480000000000	3296.39000000000
30.1290000000000	3291.98000000000
31.5340000000000	3287.36000000000
32.0140000000000	3282.41000000000
32.7960000000000	3278.26000000000
34.2040000000000	3273.83000000000
34.6830000000000	3269.77000000000
35.4660000000000	3266.90000000000
36.5520000000000	3264.03000000000
37.0300000000000	3260.57000000000
38.1120000000000	3257.92000000000
38.6000000000000	3255.24000000000
39.9980000000000	3253.56000000000
40.4780000000000	3251.11000000000
41.2660000000000	3249.43000000000
42.3600000000000	3247.97000000000
42.8420000000000	3245.70000000000
43.9440000000000	3244.21000000000
44.7250000000000	3242.46000000000
45.8220000000000	3240.97000000000
46.6080000000000	3238.99000000000
47.7120000000000	3237.83000000000
48.4970000000000	3236.41000000000
49.5750000000000	3235.18000000000
50.3550000000000	3233.80000000000
51.1380000000000	3232.64000000000
51.9190000000000	3231.44000000000
53.0040000000000	3230.34000000000
53.7820000000000	3228.76000000000
55.1840000000000	3227.95000000000
55.6560000000000	3226.76000000000
56.1290000000000	3225.81000000000
57.2250000000000	3224.80000000000
58.0120000000000	3223.91000000000
59.0990000000000	3223.02000000000
59.5740000000000	3221.77000000000
60.6670000000000	3221.26000000000
61.1470000000000	3220.18000000000
61.9300000000000	3219.80000000000
62.7150000000000	3219.13000000000
63.4990000000000	3218.87000000000
64.2820000000000	3217.85000000000
65.2330000000000	3217.34000000000
65.7120000000000	3217.05000000000
66.4910000000000	3216.42000000000
66.9650000000000	3215.87000000000
67.7530000000000	3215.24000000000
68.5310000000000	3214.45000000000
69.3140000000000	3213.74000000000
69.7900000000000	3213.44000000000
70.5670000000000	3213.15000000000
71.3470000000000	3212.73000000000
72.6070000000000	3212.45000000000
73.3870000000000	3212.18000000000
74.1710000000000	3211.58000000000
74.9590000000000	3211.32000000000
75.4380000000000	3211
75.9170000000000	3210.80000000000
76.7040000000000	3209.92000000000
77.1800000000000	3209.85000000000
77.6530000000000	3209.66000000000
78.1310000000000	3209.37000000000
78.9150000000000	3209.20000000000
79.7050000000000	3208.42000000000
80.9680000000000	3207.95000000000
81.4430000000000	3207.56000000000
82.7100000000000	3207.04000000000
83.4980000000000	3206.60000000000
84.5890000000000	3206.22000000000
85.0630000000000	3205.84000000000
85.8480000000000	3205.56000000000
86.3330000000000	3205.32000000000
87.1170000000000	3205.05000000000
87.5930000000000	3204.70000000000
88.3810000000000	3204.27000000000
89.1550000000000	3204.20000000000
89.6320000000000	3203.79000000000
90.4100000000000	3203.53000000000
90.8940000000000	3203.39000000000
91.3720000000000	3203.32000000000
92.1530000000000	3203.10000000000
93.4170000000000	3202.85000000000
93.8920000000000	3202.58000000000
94.3700000000000	3202.50000000000
94.8460000000000	3202.47000000000
95.3170000000000	3202.40000000000
95.7919999999999	3202.26000000000
96.5690000000000	3201.93000000000
97.0449999999999	3201.77000000000
97.5269999999999	3201.76000000000
98.7889999999999	3201.47000000000
99.2659999999999	3201.36000000000
100.354000000000	3201.25000000000
101.436000000000	3201.09000000000
101.910000000000	3200.94000000000
103.166000000000	3200.88000000000
103.645000000000	3200.76000000000
104.126000000000	3200.68000000000
104.601000000000	3200.68000000000
105.554000000000	3200.33000000000
106.032000000000	3200.29000000000
106.982000000000	3200.02000000000
107.458000000000	3199.99000000000
107.937000000000	3199.93000000000
108.415000000000	3199.86000000000
109.198000000000	3199.77000000000
110.147000000000	3199.48000000000
110.620000000000	3199.33000000000
111.405000000000	3199.27000000000
112.664000000000	3199.11000000000
113.618000000000	3198.99000000000
114.880000000000	3198.87000000000
115.354000000000	3198.69000000000
115.830000000000	3198.68000000000
116.921000000000	3198.37000000000
117.396000000000	3198.12000000000
118.179000000000	3197.93000000000
118.968000000000	3197.78000000000
119.751000000000	3197.72000000000
121.018000000000	3197.48000000000
121.490000000000	3197.44000000000
122.750000000000	3197.30000000000
124.181000000000	3197.22000000000
124.654000000000	3197.17000000000
125.127000000000	3197.11000000000
126.549000000000	3197.06000000000
127.029000000000	3197.04000000000
127.812000000000	3196.96000000000
128.901000000000	3196.84000000000
129.861000000000	3196.63000000000
130.332000000000	3196.55000000000
130.809000000000	3196.41000000000
132.068000000000	3196.20000000000
132.543000000000	3196.10000000000
133.326000000000	3195.92000000000
134.165000000000	3195.84000000000
135.415000000000	3195.77000000000
136.680000000000	3195.67000000000
137.155000000000	3195.60000000000
137.627000000000	3195.52000000000
138.720000000000	3195.44000000000
139.200000000000	3195.30000000000
140.288000000000	3195.23000000000
140.770000000000	3195.10000000000
141.724000000000	3194.99000000000
142.208000000000	3194.97000000000
142.684000000000	3194.83000000000
143.160000000000	3194.82000000000
143.641000000000	3194.79000000000
144.122000000000	3194.61000000000
144.906000000000	3194.49000000000
145.859000000000	3194.16000000000
146.336000000000	3194.13000000000
146.915000000000	3194.05000000000
147.692000000000	3193.80000000000
148.967000000000	3193.48000000000
149.445000000000	3193.34000000000
150.528000000000	3193.20000000000
151.484000000000	3193.02000000000
151.958000000000	3193.01000000000
152.740000000000	3192.85000000000
153.515000000000	3192.83000000000
154.301000000000	3192.73000000000
154.777000000000	3192.68000000000
155.730000000000	3192.54000000000
156.208000000000	3192.53000000000
157.161000000000	3192.31000000000
157.637000000000	3192.20000000000
158.109000000000	3192.14000000000
158.892000000000	3191.82000000000
159.975000000000	3191.32000000000
161.056000000000	3190.94000000000
161.531000000000	3190.44000000000
162.002000000000	3190.23000000000
162.477000000000	3190.09000000000
162.951000000000	3189.82000000000
164.036000000000	3189.47000000000
164.811000000000	3189.27000000000
165.287000000000	3188.91000000000
166.472000000000	3188.56000000000
167.255000000000	3188.30000000000
168.033000000000	3187.84000000000
168.510000000000	3187.72000000000
168.986000000000	3187.51000000000
169.459000000000	3187.22000000000
169.936000000000	3186.89000000000
170.712000000000	3186.81000000000
171.490000000000	3185.96000000000
172.742000000000	3185.82000000000
173.520000000000	3185.54000000000
174.300000000000	3185.45000000000
175.080000000000	3185.32000000000
175.550000000000	3185.17000000000
176.635000000000	3185.06000000000
177.114000000000	3185.06000000000
177.892000000000	3184.56000000000
178.367000000000	3184.47000000000
179.453000000000	3184.34000000000
181.245000000000	3184.23000000000
181.721000000000	3184.20000000000
182.676000000000	3184.18000000000
183.154000000000	3184.17000000000
183.940000000000	3184.11000000000
184.416000000000	3184.03000000000
184.896000000000	3183.98000000000
185.982000000000	3183.73000000000
186.457000000000	3183.70000000000
186.934000000000	3183.40000000000
188.025000000000	3183.28000000000
188.809000000000	3183.22000000000
189.897000000000	3182.98000000000
191.326000000000	3182.93000000000
192.282000000000	3182.88000000000
193.234000000000	3182.78000000000
193.713000000000	3182.73000000000
194.193000000000	3182.59000000000
194.668000000000	3182.47000000000
195.901000000000	3182.33000000000
196.376000000000	3182.18000000000
197.154000000000	3182.13000000000
198.568000000000	3181.98000000000
199.047000000000	3181.92000000000
199.834000000000	3181.83000000000
200.723000000000	3181.79000000000
201.991000000000	3181.66000000000
202.467000000000	3181.62000000000
203.253000000000	3181.54000000000
203.730000000000	3181.35000000000
204.690000000000	3181.29000000000
205.166000000000	3181.28000000000
205.646000000000	3181.23000000000
206.425000000000	3181.16000000000
207.206000000000	3181.07000000000
208.943000000000	3181.01000000000
209.417000000000	3180.95000000000
210.201000000000	3180.92000000000
210.677000000000	3180.90000000000
211.790000000000	3180.74000000000
212.271000000000	3180.72000000000
212.750000000000	3180.65000000000
213.224000000000	3180.44000000000
214.180000000000	3180.42000000000
214.655000000000	3180.33000000000
215.448000000000	3180.16000000000
216.551000000000	3179.97000000000
217.331000000000	3179.87000000000
218.108000000000	3179.73000000000
219.205000000000	3179.68000000000
219.679000000000	3179.62000000000
221.418000000000	3179.52000000000
221.894000000000	3179.46000000000
223.155000000000	3179.38000000000
224.110000000000	3179.34000000000
224.589000000000	3179.30000000000
225.060000000000	3179.25000000000
225.535000000000	3179.22000000000
226.317000000000	3179.02000000000
227.103000000000	3178.85000000000
227.886000000000	3178.73000000000
229.148000000000	3178.56000000000
230.573000000000	3178.42000000000
231.044000000000	3178.33000000000
231.586000000000	3178.32000000000
232.067000000000	3178.27000000000
232.544000000000	3178.24000000000
233.331000000000	3178.16000000000
233.806000000000	3178.16000000000
234.282000000000	3178.12000000000
234.758000000000	3178.07000000000
235.236000000000	3178
236.015000000000	3177.90000000000
236.489000000000	3177.82000000000
236.967000000000	3177.82000000000
238.411000000000	3177.72000000000
238.944000000000	3177.69000000000
239.900000000000	3177.61000000000
240.681000000000	3177.58000000000
241.461000000000	3177.56000000000
241.938000000000	3177.53000000000
242.412000000000	3177.53000000000
242.890000000000	3177.51000000000
243.841000000000	3177.45000000000
244.317000000000	3177.43000000000
245.266000000000	3177.36000000000
246.219000000000	3177.34000000000
246.699000000000	3177.31000000000
247.961000000000	3177.26000000000
248.435000000000	3177.23000000000
249.390000000000	3177.23000000000
249.863000000000	3177.19000000000
250.641000000000	3177
251.119000000000	3176.96000000000
251.895000000000	3176.86000000000
252.973000000000	3176.83000000000
253.447000000000	3176.78000000000
254.398000000000	3176.69000000000
255.352000000000	3176.60000000000
255.826000000000	3176.57000000000
256.614000000000	3176.34000000000
257.568000000000	3176.25000000000
258.351000000000	3176.19000000000
258.828000000000	3176.02000000000
259.304000000000	3176.01000000000
259.778000000000	3175.98000000000
260.861000000000	3175.80000000000
261.649000000000	3175.64000000000
262.907000000000	3175.56000000000
263.381000000000	3175.49000000000
263.856000000000	3175.47000000000
264.331000000000	3175.41000000000
265.421000000000	3175.30000000000
266.198000000000	3175.14000000000
266.983000000000	3175.08000000000
268.370000000000	3174.99000000000
268.845000000000	3174.83000000000
269.320000000000	3174.58000000000
270.101000000000	3174.42000000000
271.195000000000	3174.15000000000
271.669000000000	3174.08000000000
272.757000000000	3173.97000000000
273.710000000000	3173.79000000000
274.491000000000	3173.75000000000
275.925000000000	3173.67000000000
276.401000000000	3173.67000000000
276.883000000000	3173.66000000000
277.360000000000	3173.62000000000
277.837000000000	3173.53000000000
278.313000000000	3173.50000000000
278.789000000000	3173.49000000000
279.263000000000	3173.44000000000
280.219000000000	3173.32000000000
281.312000000000	3173.22000000000
282.103000000000	3172.94000000000
282.575000000000	3172.88000000000
283.049000000000	3172.85000000000
283.832000000000	3172.71000000000
285.094000000000	3172.35000000000
285.572000000000	3172.14000000000
286.833000000000	3171.90000000000
287.312000000000	3171.89000000000
288.740000000000	3171.79000000000
289.521000000000	3171.74000000000
289.999000000000	3171.68000000000
291.094000000000	3171.37000000000
291.876000000000	3171.08000000000
292.351000000000	3171.06000000000
292.824000000000	3170.95000000000
293.609000000000	3170.90000000000
294.089000000000	3170.90000000000
294.565000000000	3170.74000000000
295.043000000000	3170.55000000000
295.825000000000	3170.45000000000
296.605000000000	3170.43000000000
298.032000000000	3170.38000000000
298.985000000000	3170.36000000000
299.461000000000	3170.34000000000
300.723000000000	3170.24000000000
301.509000000000	3170.08000000000
301.990000000000	3170.01000000000
303.245000000000	3169.90000000000
303.720000000000	3169.73000000000
304.196000000000	3169.72000000000
305.633000000000	3169.60000000000
307.068000000000	3169.54000000000
307.546000000000	3169.51000000000
308.018000000000	3169.45000000000
309.111000000000	3169.36000000000
309.595000000000	3169.30000000000
310.380000000000	3169.16000000000
310.855000000000	3169.10000000000
312.115000000000	3168.90000000000
313.071000000000	3168.81000000000
313.550000000000	3168.78000000000
314.812000000000	3168.40000000000
315.288000000000	3168.31000000000
315.773000000000	3168.27000000000
316.269000000000	3168.25000000000
316.836000000000	3168.17000000000
317.351000000000	3168.15000000000
318.219000000000	3168.04000000000
319.844000000000	3167.98000000000
320.473000000000	3167.92000000000
321.427000000000	3167.84000000000
322.309000000000	3167.76000000000
322.788000000000	3167.74000000000
324.038000000000	3167.66000000000
324.574000000000	3167.60000000000
325.046000000000	3167.54000000000
326.449000000000	3167.50000000000
326.923000000000	3167.46000000000
328.669000000000	3167.39000000000
329.144000000000	3167.37000000000
329.683000000000	3167.32000000000
330.949000000000	3167.25000000000
331.482000000000	3167.24000000000
332.436000000000	3167.23000000000
332.974000000000	3167.21000000000
333.918000000000	3167.14000000000
334.451000000000	3167.12000000000
335.397000000000	3167.09000000000
335.870000000000	3167.09000000000
336.349000000000	3166.91000000000
336.825000000000	3166.89000000000
337.298000000000	3166.89000000000
338.561000000000	3166.74000000000
339.338000000000	3166.63000000000
339.874000000000	3166.62000000000
343.537000000000	3166.51000000000
347.542000000000	3166.51000000000
};
\addplot [line width=1, colQRkit]
table {%
 2.59600000000000	4529.42000000000
4.88600000000000	4419.14000000000
7.13300000000000	4350.64000000000
9.43800000000000	4293.80000000000
11.6880000000000	4253.31000000000
13.9320000000000	4177.51000000000
16.1920000000000	4117.36000000000
18.4970000000000	4058.96000000000
25.3430000000000	4005.22000000000
27.5540000000000	3932.07000000000
34.3190000000000	3864.74000000000
36.6330000000000	3804.32000000000
43.5700000000000	3755.61000000000
48.0920000000000	3695.38000000000
54.9580000000000	3658.61000000000
59.4690000000000	3602.43000000000
66.6660000000000	3566.71000000000
71.2200000000000	3524.96000000000
75.9130000000000	3499.71000000000
78.3220000000000	3476.81000000000
82.9240000000000	3450.08000000000
89.7990000000000	3408.74000000000
91.9830000000000	3379.27000000000
94.2730000000000	3376.16000000000
96.6350000000000	3370.47000000000
98.9850000000000	3360.13000000000
103.687000000000	3311.87000000000
106.050000000000	3308.25000000000
108.448000000000	3285.89000000000
113.071000000000	3257.81000000000
117.737000000000	3238.19000000000
120.011000000000	3237.89000000000
124.536000000000	3214.92000000000
129.354000000000	3211.48000000000
133.949000000000	3197.38000000000
138.625000000000	3191.24000000000
140.895000000000	3190.35000000000
147.794000000000	3180.58000000000
153.169000000000	3172.32000000000
158.089000000000	3170.44000000000
162.818000000000	3164.61000000000
169.828000000000	3160.53000000000
174.463000000000	3156.71000000000
181.557000000000	3154.08000000000
183.758000000000	3151.25000000000
186.046000000000	3150.60000000000
188.398000000000	3149.42000000000
193.075000000000	3147.50000000000
197.881000000000	3144.71000000000
204.883000000000	3142.92000000000
209.495000000000	3140.12000000000
214.198000000000	3138.44000000000
218.746000000000	3137.42000000000
223.435000000000	3135.46000000000
230.498000000000	3134.45000000000
235.220000000000	3132.90000000000
237.568000000000	3132.23000000000
244.834000000000	3131.15000000000
249.524000000000	3130.23000000000
256.420000000000	3129.64000000000
258.772000000000	3129
263.491000000000	3128.81000000000
265.818000000000	3128.14000000000
268.191000000000	3128
272.895000000000	3126.71000000000
279.969000000000	3126.03000000000
287.829000000000	3125.12000000000
290.096000000000	3124.66000000000
294.869000000000	3124.21000000000
299.483000000000	3123.76000000000
301.863000000000	3123.55000000000
304.295000000000	3123.12000000000
308.912000000000	3121.93000000000
311.290000000000	3121.55000000000
313.705000000000	3121.30000000000
315.989000000000	3120.28000000000
318.293000000000	3119.71000000000
323.055000000000	3119.36000000000
327.704000000000	3118.09000000000
330.125000000000	3117.49000000000
332.411000000000	3117.17000000000
339.437000000000	3116.13000000000
341.739000000000	3115.56000000000
346.502000000000	3114.73000000000
348.767000000000	3114.42000000000
355.843000000000	3114.02000000000
360.613000000000	3113.52000000000
365.218000000000	3112.96000000000
367.594000000000	3112.80000000000
369.998000000000	3112.77000000000
377.041000000000	3112.02000000000
381.680000000000	3111.63000000000
384.098000000000	3111.35000000000
386.363000000000	3110.94000000000
388.676000000000	3110.50000000000
393.435000000000	3110.20000000000
398.080000000000	3110.08000000000
402.809000000000	3110.02000000000
407.558000000000	3109.80000000000
412.201000000000	3109.51000000000
414.563000000000	3109.30000000000
416.973000000000	3109.26000000000
419.252000000000	3108.93000000000
421.565000000000	3108.72000000000
426.407000000000	3108.46000000000
429.103000000000	3108.27000000000
431.376000000000	3108.10000000000
433.677000000000	3107.93000000000
436.028000000000	3107.79000000000
440.771000000000	3107.14000000000
443.102000000000	3106.67000000000
445.478000000000	3106.17000000000
447.923000000000	3106.12000000000
452.567000000000	3104.95000000000
454.962000000000	3104.95000000000
459.684000000000	3103.49000000000
464.440000000000	3102.97000000000
466.713000000000	3102.84000000000
469.011000000000	3102.35000000000
471.372000000000	3102.30000000000
478.456000000000	3101.33000000000
480.895000000000	3100.75000000000
485.528000000000	3100.61000000000
487.889000000000	3100.21000000000
494.998000000000	3099.74000000000
499.558000000000	3099.21000000000
501.930000000000	3098.91000000000
504.338000000000	3098.87000000000
509.030000000000	3098.36000000000
511.397000000000	3098.15000000000
513.806000000000	3097.73000000000
516.075000000000	3097.33000000000
523.131000000000	3095.98000000000
527.891000000000	3095.01000000000
534.952000000000	3094.69000000000
537.238000000000	3094.46000000000
541.906000000000	3094.25000000000
544.304000000000	3094.22000000000
546.578000000000	3094.11000000000
548.901000000000	3094.02000000000
551.273000000000	3093.99000000000
553.713000000000	3093.85000000000
560.754000000000	3093.74000000000
563.043000000000	3093.68000000000
567.710000000000	3093.66000000000
570.111000000000	3093.65000000000
572.396000000000	3093.62000000000
577.166000000000	3093.47000000000
579.437000000000	3093.36000000000
581.732000000000	3093.28000000000
586.501000000000	3093.22000000000
591.138000000000	3093.16000000000
593.514000000000	3093.14000000000
598.221000000000	3093.10000000000
600.587000000000	3093.10000000000
605.284000000000	3093.03000000000
626.202000000000	3092.92000000000
};
\addplot [line width=1, colQRChol]
table {%
4.05700000000000	4529.42000000000
6.75500000000000	4418.78000000000
9.34400000000000	4347.03000000000
11.9620000000000	4318.19000000000
17.3220000000000	4239.33000000000
20.1060000000000	4221.79000000000
22.8840000000000	4132.50000000000
25.4710000000000	4071.07000000000
28.1950000000000	4003.54000000000
30.7520000000000	3938.32000000000
33.3150000000000	3870.99000000000
35.9070000000000	3811.52000000000
38.4050000000000	3754.27000000000
41.0280000000000	3702.06000000000
43.8950000000000	3656.05000000000
46.5140000000000	3613.08000000000
49.0310000000000	3576.54000000000
51.6760000000000	3540.73000000000
75.1450000000000	3508.25000000000
};
\addplot [line width=1, colMoreQR]
table {%
4.11500000000000	4529.42000000000
7.99200000000000	4433.02000000000
12.3860000000000	4382.56000000000
16.0320000000000	4340.65000000000
19.8510000000000	4298.48000000000
23.5090000000000	4253.23000000000
27.3360000000000	4221.80000000000
30.9820000000000	4161.27000000000
34.8020000000000	4112.85000000000
38.6180000000000	4051.34000000000
42.4640000000000	4000.09000000000
46.3540000000000	3992.34000000000
50.0400000000000	3902.04000000000
53.8950000000000	3859.78000000000
57.5540000000000	3741.86000000000
65.7120000000000	3686.17000000000
69.3780000000000	3684.54000000000
75.4540000000000	3577.31000000000
81.4870000000000	3540.69000000000
87.4690000000000	3513.16000000000
93.4500000000000	3473.19000000000
101.570000000000	3434.77000000000
109.713000000000	3406.78000000000
113.520000000000	3379.05000000000
117.143000000000	3378.47000000000
123.118000000000	3337.58000000000
129.115000000000	3327.01000000000
135.060000000000	3306.22000000000
143.194000000000	3288.70000000000
149.162000000000	3268.06000000000
155.158000000000	3254.61000000000
158.836000000000	3251.10000000000
166.989000000000	3232.87000000000
170.682000000000	3218.17000000000
176.691000000000	3213.37000000000
180.525000000000	3206.54000000000
188.534000000000	3193.88000000000
196.866000000000	3186.03000000000
200.560000000000	3176.94000000000
208.714000000000	3171.26000000000
215.048000000000	3165.90000000000
221.057000000000	3162.64000000000
227.613000000000	3157.38000000000
233.756000000000	3154.44000000000
242.051000000000	3151.13000000000
250.184000000000	3148.03000000000
254.156000000000	3144.46000000000
258.074000000000	3142.78000000000
262.048000000000	3140.29000000000
266.017000000000	3138.85000000000
273.704000000000	3134.89000000000
277.612000000000	3134.02000000000
286.269000000000	3128.59000000000
292.446000000000	3126.13000000000
298.803000000000	3125.11000000000
307.436000000000	3123.62000000000
311.403000000000	3122.62000000000
315.193000000000	3122.13000000000
321.505000000000	3120.86000000000
325.435000000000	3120.83000000000
334.226000000000	3118.65000000000
340.453000000000	3117.59000000000
344.378000000000	3117.07000000000
348.310000000000	3116.87000000000
352.256000000000	3116.71000000000
356.197000000000	3115.26000000000
362.616000000000	3114.23000000000
368.816000000000	3112.39000000000
372.761000000000	3111.61000000000
376.734000000000	3110.73000000000
380.545000000000	3110.23000000000
389.334000000000	3108.80000000000
393.321000000000	3108.26000000000
399.583000000000	3107.86000000000
405.777000000000	3107.46000000000
409.727000000000	3107.05000000000
415.910000000000	3106.68000000000
420.277000000000	3106.45000000000
424.258000000000	3106.37000000000
428.249000000000	3105.83000000000
432.191000000000	3105.15000000000
436.166000000000	3104.52000000000
447.732000000000	3104.09000000000
452.021000000000	3104.04000000000
456.242000000000	3103.82000000000
460.475000000000	3103.71000000000
474.600000000000	3103.09000000000
481.802000000000	3103.08000000000
486.241000000000	3103.01000000000
490.623000000000	3102.67000000000
496.941000000000	3102.52000000000
500.879000000000	3102.32000000000
504.841000000000	3102.18000000000
511.225000000000	3101.97000000000
515.118000000000	3101.81000000000
519.100000000000	3101.66000000000
522.998000000000	3101.58000000000
526.950000000000	3101.46000000000
530.935000000000	3101.39000000000
537.556000000000	3101.21000000000
541.723000000000	3101.06000000000
545.633000000000	3100.90000000000
549.551000000000	3100.64000000000
553.476000000000	3100.45000000000
557.397000000000	3100.24000000000
563.779000000000	3100.01000000000
567.654000000000	3099.86000000000
571.769000000000	3099.80000000000
575.764000000000	3099.79000000000
579.651000000000	3099.65000000000
583.578000000000	3099.59000000000
587.574000000000	3099.53000000000
591.459000000000	3099.50000000000
595.560000000000	3099.46000000000
599.487000000000	3099.32000000000
603.376000000000	3099.25000000000
607.320000000000	3099.23000000000
611.240000000000	3099.18000000000
615.145000000000	3099.18000000000
619.095000000000	3099.14000000000
625.391000000000	3098.99000000000
629.355000000000	3098.94000000000
633.312000000000	3098.92000000000
637.242000000000	3098.90000000000
641.216000000000	3098.77000000000
647.496000000000	3098.75000000000
651.413000000000	3098.74000000000
657.687000000000	3098.67000000000
666.246000000000	3098.63000000000
670.138000000000	3098.60000000000
674.085000000000	3098.49000000000
678.049000000000	3098.49000000000
686.580000000000	3098.42000000000
692.937000000000	3098.38000000000
696.820000000000	3098.37000000000
702.929000000000	3098.32000000000
706.854000000000	3098.28000000000
710.801000000000	3098.25000000000
714.691000000000	3098.18000000000
718.621000000000	3098.16000000000
722.546000000000	3098.11000000000
726.575000000000	3098.08000000000
732.666000000000	3098.04000000000
736.685000000000	3098.02000000000
740.684000000000	3098.01000000000
744.560000000000	3098
748.732000000000	3097.96000000000
752.675000000000	3097.93000000000
756.707000000000	3097.91000000000
760.609000000000	3097.83000000000
766.767000000000	3097.79000000000
772.909000000000	3097.77000000000
777.008000000000	3097.76000000000
780.950000000000	3097.76000000000
784.862000000000	3097.76000000000
788.928000000001	3097.76000000000
792.841000000001	3097.76000000000
796.770000000000	3097.72000000000
800.708000000000	3097.71000000000
806.872000000001	3097.65000000000
810.800000000001	3097.64000000000
814.746000000001	3097.62000000000
818.665000000001	3097.61000000000
822.616000000001	3097.61000000000
826.542000000001	3097.58000000000
832.751000000001	3097.39000000000
837.083000000001	3097.38000000000
840.961000000001	3097.35000000000
844.902000000001	3097.34000000000
848.836000000001	3097.30000000000
852.752000000001	3097.29000000000
856.680000000001	3097.29000000000
860.595000000001	3097.24000000000
864.540000000001	3097.20000000000
868.474000000001	3097.17000000000
872.369000000001	3097.12000000000
876.314000000001	3097.11000000000
880.248000000001	3097.09000000000
884.156000000001	3097.07000000000
888.111000000001	3096.98000000000
892.076000000001	3096.97000000000
896.178000000001	3096.77000000000
900.066000000001	3096.76000000000
906.425000000001	3096.51000000000
910.646000000001	3096.49000000000
915.177000000001	3096.45000000000
919.923000000001	3096.38000000000
927.073000000001	3096.35000000000
931.363000000001	3096.33000000000
936.600000000001	3096.32000000000
942.336000000001	3096.20000000000
947.288000000001	3096.01000000000
954.122000000001	3096
958.240000000001	3095.96000000000
962.122000000001	3095.95000000000
968.288000000001	3095.93000000000
972.157000000001	3095.92000000000
976.056000000001	3095.79000000000
979.942000000001	3095.73000000000
986.130000000001	3095.71000000000
990.000000000001	3095.71000000000
993.880000000001	3095.70000000000
997.796000000001	3095.67000000000
1001.68500000000	3095.62000000000
1005.60100000000	3095.62000000000
1009.51100000000	3095.60000000000
1013.40400000000	3095.56000000000
1017.30300000000	3095.44000000000
1021.17900000000	3095.36000000000
1025.22500000000	3095.29000000000
1029.09700000000	3095.17000000000
1035.31300000000	3095.01000000000
1039.24100000000	3094.97000000000
1043.10400000000	3094.94000000000
1047.05400000000	3094.91000000000
1051.05900000000	3094.90000000000
1057.16000000000	3094.78000000000
1061.01800000000	3094.77000000000
1064.90600000000	3094.76000000000
1071.10400000000	3094.74000000000
1075.05000000000	3094.73000000000
1078.89800000000	3094.73000000000
1084.98600000000	3094.73000000000
1102.62800000000	3094.73000000000
1117.81600000000	3094.73000000000
};
\end{axis}

\end{tikzpicture}